\documentclass[12pt,a4paper]{amsart}
\usepackage{amsmath,amsthm,amsfonts,amssymb,mathrsfs}
\usepackage{subfigure} 
\usepackage{graphicx}
\usepackage{color}
\usepackage{caption}
\usepackage[arrow, matrix, curve]{xy}
\usepackage{tikz}
\usepackage{ulem}
\usetikzlibrary{patterns}

\usepackage{comment}

\newcommand\shorter[1]{}
\makeatletter
\renewcommand{\@secnumfont}{\bfseries}
\makeatother
\newcommand\R{\mathbb{R}}

\newcommand\vardefODE{x \in \overline{\text{$I$}}, t>0,}
\newcommand\vardefPDE{x \in \text{$I$}, t>0,}
\newcommand\vardefODEhat{x \in \overline{\text{$I$}}, \hat{t}>0,}
\newcommand\vardefPDEhat{x \in \text{$I$}, \hat{t}>0,}
\newcommand{\eb}{b}
\newcommand{\MB}{B}

\newcommand{\del}[2]{\frac{\partial #1}{\partial #2}}
\newcommand{\dell}[2]{\frac{\partial^2 #1}{\partial #2^2}}

\def\ueps{u_1^{\delta}}
\def\veps{u_2^{\delta}}
\def\weps{v^{\delta}}

\newtheorem{theorem}{Theorem}
\newtheorem{cor}[theorem]{Corollary}

\newtheorem{lemma}[theorem]{Lemma}
\newtheorem{sublemma}[theorem]{Sublemma}
\newtheorem{prop}[theorem]{Proposition}

\newtheorem{definition}[theorem]{Definition}
\newtheorem{assumption}[theorem]{Assumption}

\numberwithin{equation}{section}
\numberwithin{theorem}{section}
\numberwithin{figure}{section}

\theoremstyle{remark}
\newtheorem{rem}[theorem]{Remark}

\usepackage{delarray}
\usepackage{enumerate}

\title[]{Stable patterns with jump discontinuity in systems with Turing instability and hysteresis}
\author{Steffen H\"arting}
\author{Anna Marciniak-Czochra}
\author{Izumi Takagi}

\address[Steffen H\"arting]{Institute of Applied Mathematics and BIOQUANT, Heidelberg University, Im Neuenheimer Feld 294, 69120 Heidelberg, Germany}
\address[Anna Marciniak-Czochra]{Institute of Applied Mathematics, IWR and BIOQUANT, Heidelberg University, Im Neuenheimer Feld 294, 69120 Heidelberg, Germany}
\address[Izumi Takagi]{Mathematical Institute, Tohoku University, Sendai, 980-8578, Japan}

\begin{document}

\begin{abstract}

Classical models of pattern formation are based on diffusion-driven
instability (DDI) of constant stationary solutions of reaction-diffusion
equations, which leads to emergence of stable, regular Turing patterns formed
around that equilibrium. In this paper we show that coupling
reaction-diffusion equations with ordinary differential equations (ODE) may
lead to a novel pattern formation phenomenon in which DDI causes
destabilization of both constant solutions and Turing patterns. Bistability
and hysteresis effects in the null sets of model nonlinearities yield formation
of {\it far from the equilibrium} patterns with jump discontinuity.  We derive conditions for 
stability of stationary solutions with jump discontinuity in a suitable topology which
allows us to include the discontinuity points and leads to the definition of
$(\varepsilon_0,A)$-stability. Additionally, we provide conditions on stability
of patterns in a quasi-stationary model reduction. The analysis is illustrated on the example of
three-component model of receptor-ligand binding. The proposed model
provides an example of a mechanism of {\it de novo} formation of {\it far from
the equilibrium} patterns in reaction-diffusion-ODE models involving co-existence of
DDI and hysteresis.

\end{abstract}
\maketitle
\section{Introduction}
 Since the seminal paper of Alan Turing  \cite{Turing},  mathematical models of biological pattern formation have been constructed by using systems of reaction-diffusion 
 equations exhibiting diffusion-driven instability (DDI). DDI is related to a local behavior of solutions of a reaction-diffusion system in the neighborhood of a constant stationary solution that is destabilized through diffusion.  It may lead to emergence of stable continuous and spatially periodic structures around the destabilized constant equilibrium. The Turing concept became a paradigm for pattern formation and led to development of numerous theoretical models, though its biological verification has remained elusive \cite{Akam, G-M, Murray}. 

The models are based on the idea that cells differentiate according to positional information which is supplied to them by diffusing biochemical morphogens. Different morphogen concentrations are able to activate transcription of distinct target genes and thus lead to cell differentiation. However, both regulatory and signaling molecules (ligands) act by binding and activating receptor molecules, which are located in the cell membrane \cite{Lauffenburger,Mbook}. This observation leads to a hypothesis that the positional value of a cell may be determined by the density of bound receptors which do not diffuse. Taking dynamics of receptors into account leads to systems coupling reaction-diffusion equations with ordinary differential equations (ODE) \cite{Klika,Marciniak03,MC06}, called also receptor-based models. Such systems can be obtained as a homogenization limit of the models describing coupling of cell-localized processes with cell-to-cell communication via diffusion in a cell assembly \cite{MC12,MCP}.  Nonlinear interactions of diffusive and non-diffusive components  arise also from modeling of interactions between cellular or intracellular processes. Reaction-diffusion-ODE models have recently been employed in various biological contexts, see e.g., \cite{Hock,Klika,Marciniak03,MCK07,Pham,USOO}. 

Although receptor-based models may exhibit Turing instability as shown in \cite{Marciniak03} and discussed more recently on several examples from mathematical biology in \cite{Klika}, they are different from classical Turing-type models. Since DDI is induced by self-enhancement in the subsystem with smallest diffusion \cite{Sakamoto12}, in reaction-diffusion-ODE models instability is induced by growth properties of the ODE subsystem. Therefore, coupling of reaction-diffusion equations to ODEs may lead to DDI which does not arise if the non-diffusive components are neglected. 

To understand the role of non-diffusive components in pattern formation process, we focus on systems involving a single reaction-diffusion equation coupled to ODEs on one-dimensional domain $\overline{\text{$I$}}$. In general, equations of such models can be represented by the following initial-boundary value problem,
\begin{equation}
\begin{aligned}\label{sys:deg}
\del{u}{t} &= f(u,v),    &\vardefODE\\
\del{v}{t} &= D_v \dell{v}{x} + g(u,v),\qquad  &\vardefPDE\\
(u(0),v(0))&= (u_0,v_0) \in C(\overline{\text{$I$}})\times (C^2(\text{$I$}) \cap C(\overline{\text{$I$}})),
\end{aligned}
\end{equation}
with homogeneous Neumann boundary conditions for $v$.

It is an interesting case, since a scalar reaction-diffusion equation cannot exhibit stable spatially heterogeneous patterns \cite{CaHo}. Coupling it to an ODE fulfilling autocatalysis condition, i.e. $\partial f/\partial u \, >0$ at the constant equilibrium, leads to DDI.  However, in this case all Turing patterns are unstable, i.e. the same mechanism which destabilizes constant solutions, destabilizes also all continuous spatially heterogeneous stationary solutions  \cite{MCKS12,MCKS13}.  This instability result  holds also for patterns with jump discontinuity in case of a special class of nonlinearities \cite{MCKS12,MCKS13}.  The question then arises as to which patterns, if any, can be exhibited in such models.

As shown in \cite{MCKS12}, it may happen that there exist no stable stationary patterns  and the emerging spatially heterogeneous structures have a dynamical character. Simulations of different models of this form indicate formation of  dynamical, multimodal and apparently irregular structures, the shape of which depends strongly on initial conditions \cite{HM13,MCK07,Pham}.   On the other hand, reaction-diffusion-ODE models may give rise to discontinuous patterns due to the hysteresis effect in the model nonlinearities, i.e. when the equation $f(u,v)=0$ has multiple quasi-stationary solutions $v_{i}=H_{i}(u)$.  Diffusion tries to average different states and may lead to formation of far from equilibrium structures. Hysteresis yields emergence of stationary solutions with a jump discontinuity, which may be monotone, periodic or irregular.  The hysteresis-driven pattern formation has been investigated analytically in systems with two stable constant equilibria, i.e. without DDI  \cite{MC06,MCNT,KMCProc}. 
In such a case patterns depend strongly on initial conditions and require large initial heterogeneity. 

Consequently, the next question is whether a coupling of the two mechanisms, DDI and hysteresis, can provide a mechanism of {\it de novo} formation of stable  {\it far from equilibrium} patterns. By {\it de novo} pattern formation we understand emergence of spatially heterogeneous structures starting from small perturbation of the constant stationary solution. According to our knowledge, systems with DDI and hysteresis have not been studied so far.

Motivated by these observations, we propose a new receptor-based model exhibiting DDI  and hysteresis which lead to {\it de novo} formation of stable patterns with jump-type discontinuities. The model is a modification of the simplest model describing receptor-ligand binding dynamics developed in  \cite{MCP, Marciniak03}. The modification accounts for saturation in production terms, which provides uniform boundedness of the model solutions. We consider two- and three-component systems, i.e., 
\begin{itemize}
 \item a model consisting of one ordinary differential equation and one reaction-diffusion equation (1ODE+1RDE),
 \item a model consisting of two ordinary differential equations and one reaction-diffusion equation (2ODE+1RDE),
\end{itemize}
While the principle of DDI is similar to Turing's original idea, coupling to an ODE is crucial for stability of steady states in this paper. One interpretation of the observed phenomenon is that non-diffusive species, such as receptors, can facilitate formation of patterns with sharp transitions.

In this paper, we derive conditions for  stability of such patterns.  Since stationary solutions for nondiffusing variables have jump discontinuities, linearized stability analysis requires careful treatment. To cope with this difficulty, we consider a special topology which allows us to include the discontinuity points and yields a notion of $(\varepsilon_0,A)$-stability. Our approach follows the one of Weinberger proposed in refs.\ \cite{Weinberger,Weinberger:87} to study discontinuous solutions in a model of density dependent diffusion. A novelty of our paper lies in the stability analysis of patterns with jump discontinuity in systems involving several ODEs which may be reduced to a subsystem by using quasi-steady state approximation. We show that under suitable conditions on the nonlinearities, stability of stationary solutions is preserved under quasi-steady state reduction.

The paper is divided into five sections:
In Section \ref{sec:stab} the main results of this paper are formulated. We consider systems of type \eqref{sys:deg} with scalar-valued $v$ and $u$ having values in either $\mathbb{R}$ or $\mathbb{R}^2$ and provide conditions under which steady states with jump discontinuity are stable in the topology introduced in ref. \cite{Weinberger}. Then, we extend the stability results to systems with a small parameter and the corresponding quasi-stationary approximation.
In Section \ref{sec:ex}, we present a specific receptor-based model of type \eqref{sys:deg} exhibiting both pattern formation mechanisms, i.e.\ DDI and hysteresis. We show that the model satisfies the theory developed in this paper. Model analysis is supplemented by numerical simulations showing {\it de novo} formation of {\it far from equilibrium} patterns exhibiting jump-discontinuities of which the location depends both on the size of diffusion and on the initial data.  Section \ref{sec:mainproof} presents a detailed proof of the stability theorems stated in Section \ref{sec:stab}.
In Section \ref{sec:exproof}, we prove all statements from Section \ref{sec:ex} related to the analysis of DDI initiated formation of hysteresis-driven patterns with jump discontinuity.

\section{Main results. Stability conditions for steady states with jump discontinuity}\label{sec:stab}
We state our results on the stability of steady states of \eqref{sys:deg} with finitely many jump discontinuities under the following assumptions.
\begin{assumption}\label{assump} \hfill \break  \begin{enumerate}
  \item[{\rm (1)}] $v$ is a scalar, while $u$ is either a scalar or a vector consisting of two components.
  \item[{\rm (2)}] There exist constants $\underline{U}$, $\overline{U}$, $ \underline{V}$, $\overline{V}$, $\underline{u}$, $\overline{u}$, $\underline{v}$, $ \overline{v} \in \R$ (respectively $\underline{U}=(\underline{U}_1, \underline{U}_2), \, \overline{U}=(\overline{U}_1, \overline{U}_2)$,\ $\underline{u}=(\underline{u}_1, \underline{u}_2),\, \overline{u} =(\overline{u}_1, \overline{u}_2) \in \R^2$) such that {\rm (i)} $\underline{U} < \underline{u} < \overline{u} < \overline{U}$, $ \underline{V} < \underline{v} < \overline{v} < \overline{V}$ and {\rm (ii)}
 \begin{equation*}
\begin{aligned}
& -\infty < \underline{u} \leq \liminf_{t \to+\infty} \min_{x\in \overline{I}} u(t,x) \leq \limsup_{t\to+\infty} \max_{x\in \overline{I}} u(t,x) \leq \overline{u} < +\infty, \\
& -\infty < \underline{v} \leq \liminf_{t \to+\infty} \min_{x\in \overline{I}} v(t,x) \leq \limsup_{t\to+\infty} \max_{x\in \overline{I}} v(t,x) \leq \overline{v} < +\infty,
\end{aligned}
\end{equation*}
provided that initial function $(u_0, v_0)$ satisfies $\underline{U} \leq u_0 \leq \overline{U}$ and $\underline{V} \leq v_0 \leq \overline{V}$.
The inequalities are satisfied componentwise in the case $u=(u_1, u_2)$.
 \item[{\rm (3)}]  $f$ and $g$ are twice continuously differentiable on $[\underline{U}, \overline{U}] \times [\underline{V}, \overline{V}]$ if $u$ is a scalar, and on $[\underline{U}_1, \overline{U}_1] \times [\underline{U}_2, \overline{U}_2] \times  [\underline{V}, \overline{V}]$ if $u$ is a vector.
 \item[{\rm (4)}] Initial functions $(u_0, v_0)$ are in $C(\overline{I}) \times (C^2(I) \cap C^1(\overline{I}))$ (respectively $C(\overline{I})^2 \times (C^2(I) \cap C^1(\overline{I}))$, and satisfy $v_0^\prime(0)=v_0^\prime(l)=0$, where $I=(0, l)$.
 \end{enumerate}
\end{assumption} 

Assumption \ref{assump} is satisfied if there exists an invariant rectangle, see \cite{R84,Smoller}.

\subsection{The topology}$\,$ \\
In this section, we introduce $(\varepsilon_0,A)$-stability following the idea of Weinberger in \cite{Weinberger} and provide conditions for stability of steady states of system \eqref{sys:deg} in this topology.

For a bounded interval $\text{$I$}=(0,l)$, let $BV(\overline{I})$ denote the space of all functions of bounded variation defined on $\overline{I}$. We define a neighborhood basis for $\tilde{u} \in BV(\overline{\text{$I$}})$ with values in $[\underline{u},\overline{u}]$, where $\underline{u},\overline{u}$ are given by Assumption \ref{assump}, as
\begin{equation*}
\begin{aligned}
N_{\varepsilon, \underline{u}, \overline{u}}(\tilde{u})=\{u \in & BV(\overline{\text{$I$}}, [\underline{u},\overline{u}]) \,  | \\ &\text{there exists } {R \subset \text{$I$}}\mbox{ such that } \|u-\tilde{u}\|_{L^{\infty}(R)}^2<\varepsilon^2 \mbox{ and } \operatorname{meas}(\text{$I$}\setminus R)<\varepsilon^4 \}.
\end{aligned}
\end{equation*}
We are particularly interested in the situation where the initial function $u_0$ is {\it close} to a steady state $\tilde u$ with finitely many jump-type discontinuities, but $u_0$ is still continuous on $\overline{\text{$I$}}$. If we assume $v_0 \in C^2(\text{$I$})\cap C(\overline{\text{$I$}})$, then the solution  $(u(t,x), v(t,x))$ of the initial-boundary value problem is continuous for all $t>0$. If the steady state $\tilde u$ is `stable', then $u(t,x)$ is expected to converge towards a function with jump discontinuity in space. Therefore, the uniform norm is not appropriate to measure the closeness, see Figure 2.1. It turns out that $N_{\varepsilon,\underline{u},\overline{u}}(\tilde u)$ provides us a reasonable topology for this purpose.
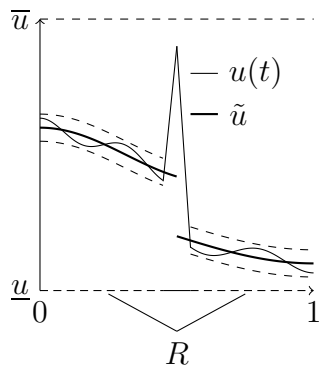
\begin{figure}[ht]
  \centering
  \begin{tikzpicture}[domain=0:1,scale=1.8]
    \draw[very thin,color=gray] (0.0,0.0) grid (0,0);
    \draw[-][densely dashed] (0,0) node[left] {$\underline{u}$} node[below] {$0$} --(0.9,0);
    \draw[->][densely dashed] (1.1,0) --(2,0) node[below]{$1$};
    \draw[-] (0.9,0)--(1.1,0);
    \draw[dashed] (0,2) -- (2,2);
    \draw[-][thin] (0.5,-0.025) -- (1,-0.3);
    \draw[-][thin] (1.5,-0.025) -- (1,-0.3);
    \draw[] (1,-0.3) node[below] {$R$};
    \draw[->] (0,0) -- (0,2) node[left] {$\overline{u}$};
    \draw[samples=100,domain=0.0:1, thick] plot(\x,   {cos(deg(0.5*\x*5))*0.2    +1});
    \draw[samples=100,domain=0.0:0.9, dashed] plot(\x,   {cos(deg(0.5*\x*5))*0.2    +1.1});
    \draw[samples=100,domain=0.0:0.9, dashed] plot(\x,   {cos(deg(0.5*\x*5))*0.2    +0.9});
    \draw[samples=100,domain=0.0:0.9] plot(\x,   {cos(deg(0.5*\x*5))*0.2    +1+0.07*cos(deg(0.5*\x*5*4))});
       
    \draw[samples=100,domain=1:2, thick] plot(\x,   {cos(deg((1-0.5)*\x*3.14))*0.2    +0.4});
    \draw[samples=100,domain=1.1:2, dashed] plot(\x,   {cos(deg((1-0.5)*\x*3.14))*0.2    +0.5});
    \draw[samples=100,domain=1.1:2, dashed] plot(\x,   {cos(deg((1-0.5)*\x*3.14))*0.2    +0.3});
    \draw[samples=100,domain=1.1:2] plot(\x,   {cos(deg((1-0.5)*\x*3.14))*0.2    +0.4 + 0.07*cos(deg(0.5*\x*5*3.14))});
        \draw[-] (0.9,0.805) -- (1,1.8);
        \draw[-] (1,1.8) -- (1.1,0.315);
        \draw[] (1.1,1.6) -- (1.3,1.6) node[right] {$u(t)$};
        \draw[thick] (1.1,1.3) -- (1.3,1.3) node[right] {$\tilde{u}$};
  \end{tikzpicture}
   \caption{Illustration of the topology applied to problem \eqref{sys:deg} for scalar $u$. Model \eqref{sys:deg} exhibits steady states with jump discontinuity and global existence of classical solutions. $\tilde{u}$ represents a steady state while $u(t,x)$ represents a solution for some $t$.}
  \label{fig:cutDisc}
\end{figure}
\noindent Following \cite{Weinberger}, $(\varepsilon_0,A)$-stability in this topology is defined as
\begin{definition}[$(\varepsilon_0,A)$-stability]\label{def-stab}
A stationary solution $(\tilde{u},\tilde{v})$ of system \eqref{sys:deg} is said to be $(\varepsilon_0,A)$-stable for positive constants $\varepsilon_0$ and $A$ if the initial functions $(u_0,v_0)$ satisfy \begin{equation}
 \|u_0 - \tilde{u}\|_{L^{\infty}(R)}^2+\|v_0-\tilde{v}\|_{H^1(\text{$I$})}^2 < \varepsilon^2
\end{equation}
for some $R \subset I$ with 
 $\operatorname{meas}(\text{$I$} \setminus R)<\varepsilon^4$,
and for some $\varepsilon \in (0,\varepsilon_0)$, then \begin{equation}
\|u(t,\cdot)-\tilde{u}\|_{L^{\infty}(R)}^2+ \|v(,\cdot)-\tilde{v}\|_{H^1(\text{$I$})}^2 < A \varepsilon^2.
\end{equation}
for all $t>0$. Here $H^1(I) = \{ u \in L^2(I) \mid u^\prime \in L^2(I)\}$ and $\Vert u\Vert_{H^1(I)} = \Vert u \Vert_{L^2(I)} + \Vert u^\prime \Vert_{L^2(I)} $.
\end{definition}

\subsection{Stability conditions} $\,$ \\
Definition \ref{def-stab} allows us to give conditions for the $(\varepsilon_0,A)$-stability of steady states with jump discontinuity. First, we analyze a system of one ordinary differential equation coupled to one reaction-diffusion equation.

\begin{theorem}[Stability for a scalar-valued $u$] \label{thm:stabC} Let $\text{$I$} =(0, l)$ be a bounded interval in $\mathbb{R}$. Under Assumption \ref{assump}, consider the system of equations
\begin{equation}
\left\{\ 
\begin{aligned}
 \del{u}{t} &= f(u,v),                           &\text{for}\ \vardefODE\\
 \del{v}{t} &= D \dell{v}{x} + g(u,v), \qquad &\text{for}\ \vardefPDE \\
 \del{v}{x} &(t,x) = 0 & \text{at}\quad x =0,\,l, \quad t>0.
\end{aligned}
\right.
\end{equation}
Let $(\tilde{u},\tilde{v})$ be a steady state with finitely many discontinuities of $\tilde{u}$.
Denote the Jacobian matrix of the kinetic system at the steady state by 
\begin{equation}
\MB(x)=\begin{pmatrix} \displaystyle \del{f}{u}(\tilde u(x), \tilde v(x)) & \displaystyle \del{f}{v}(\tilde u(x), \tilde v(x)) \vphantom{\int_M}\\
\displaystyle \del{g}{u}(\tilde u(x), \tilde v(x)) & \displaystyle \del{g}{v}(\tilde u(x), \tilde v(x)) \vphantom{\int^M}  \end{pmatrix}
\end{equation}
and assume that the following inequalities hold for all $x \in I$:
\begin{align}
 \del{f}{u}(\tilde u(x), \tilde v(x))&\leq -c_1 < 0, \label{A2.5}\\
 \del{g}{v}(\tilde u(x), \tilde v(x))&\leq -c_1< 0, \label{A2.6} \\
 \det \MB(x)& > 0, \label{A2.7}
\end{align}
where $c_1$ is a positive constant independent of $x$.
Then $(\tilde{u},\tilde{v})$ is $(\varepsilon_0,A)$-stable for a pair $(\varepsilon_0,A)$ with $0 < \varepsilon_0,A<\infty$.
\end{theorem}

The proof of Theorem \ref{thm:stabC} is similar to that of stability for a particular model with cross-diffusion in \cite{Weinberger}. We shall present it as a corollary to more general Theorem \ref{thm:stabC-3c-gen} in Section \ref{sec:mainproof}. Theorem \ref{thm:stabC} shows that the steady states of the model are stable, if both species are self-inhibitory and cross-inhibition/cross-proliferation is asymmetric. Even if each species inhibits all species, a steady state can be stable if the cross-inhibition is sufficiently small compared to the self-inhibition. However, in this case, external influx is necessary. \par

Next, we consider a model comprising of two ordinary differential equations and one reaction-diffusion equation, i.e., $u=(u_1, u_2)$: \begin{equation}\label{full}
\begin{aligned}
 \del{u_1}{t} &= f_1(u_1,u_2,v),                        &\vardefODE\\
 \del{u_2}{t} &= f_2(u_1,u_2,v),                 &\vardefODE\\
 \del{v}{t} &= D \dell{v}{x} + g(u_1,u_2,v),\qquad&\vardefPDE\\
 \del{v}{x}(t,0)&=\del{v}{x}(t,1)=0.
\end{aligned}
\end{equation}

We are interested in the stability question on a given steady state $(\tilde u_1(x), \tilde u_2(x), \tilde v(x))$ of \eqref{full} with finitely many jump discontinuities. It is convenient to put
\begin{equation}
\begin{aligned}
& a_{11}(x) = \del{f_1}{u_1}(\tilde u_1, \tilde u_2, \tilde v),
\quad a_{12}(x) = \del{f_1}{u_2}(\tilde u_1, \tilde u_2, \tilde v),
\quad a_{13}(x) = \del{f_1}{v}(\tilde u_1, \tilde u_2, \tilde v), \cr
& a_{21}(x) = \del{f_2}{u_1}(\tilde u_1, \tilde u_2, \tilde v),
\quad a_{22}(x) = \del{f_2}{u_2}(\tilde u_1, \tilde u_2, \tilde v),
\quad a_{23}(x) = \del{f_2}{v}(\tilde u_1, \tilde u_2, \tilde v), \cr
& a_{31}(x) = \del{g}{u_1}(\tilde u_1, \tilde u_2, \tilde v),
\quad a_{32}(x) = \del{g}{u_2}(\tilde u_1, \tilde u_2, \tilde v),
\quad a_{33}(x) = \del{g}{v}(\tilde u_1, \tilde u_2, \tilde v),
\end{aligned}
\end{equation}
and
\begin{equation} \label{eqn:JacobianA}
A(x) = (a_{ij}(x))_{1\leq i,\, j \leq 3}.
\end{equation}

\begin{theorem}[Stability for vector-valued $u$]\label{thm:stabC-3c-gen} Let $\text{$I$} \subset \mathbb{R}$ be a bounded interval.
Under Assumption \ref{assump}, consider a system of type \eqref{full}.
Let $(\tilde{u_1},\tilde{u_2},\tilde{v})$ be a steady state with finitely many discontinuities of $(\tilde{u_1},\tilde{u_2})$ in $x$.
Denote the Jacobian matrix of the kinetic system at the steady state by $A(x)$ as in \eqref{eqn:JacobianA}.
 Let $A_{ij}$ denote the matrix resulting from omitting the $i$-th row and the $j$-th column, i.e. $A_{ij}=(a_{kl})_{k\neq i, l\neq j}$.
Assume that $A(x)$ and $A_{ij}(x)$ satisfy
\begin{align}
& {\rm tr}\, A(x) < 0, \quad {\rm tr}\, A(x) \sum_{j=1}^3 \det A_{jj}(x) < \det A(x) < 0, & \label{A2.11} \\
& {\rm tr}\, A_{33}(x) \sum_{j=1}^3 \det A_{jj}(x) \leq \det A(x) + {\rm tr}\, A(x)\det A_{33}(x) , & \label{A2.12} \\
& 0 < 3 \det A_{33}(x) \leq {\rm tr}\, A(x) \, {\rm tr}\, A_{33}(x) + \sum_{j=1}^3 \det A_{jj}(x), & \label{A2.13} \\
& a_{33}(x)  \leq  -3\kappa <0 \quad \hbox{and} \quad \det A_{33}(x) \geq -3\kappa\, {\rm tr}\, A_{33}(x) \geq 18\kappa^2 & \label{A2.14}
\end{align}
for all $x \in \overline{I}$, where $\kappa$ is a positive constant.Then the steady state $(\tilde u_1, \tilde u_2, \tilde v) $ is $(\varepsilon_0, A)$-stable for some positive constants $\varepsilon_0$ and $A$. 
\end{theorem}

\noindent
{\bf Remark.} Condition \eqref{A2.11} is the Routh-Hurwitz condition for the matrix $A(x)$, i.e., all eigenvalues of $A(x)$ have negative real part if and only if \eqref{A2.11} holds. Note that it implies $\sum_{j=1}^3 \det A_{jj}(x) > 0$.

\medskip

The above theorem deals with general three-component systems, and it is non-trivial to check the set of conditions \eqref{A2.11}--\eqref{A2.14}, practically speaking.  However, in the case where one species, say $u_2$, reacts very rapidly, the situation becomes more manageable. Therefore,
we consider the following system with a small positive parameter $\delta$:
\begin{align}
\left\{ \ 
\begin{aligned}
 \del{\ueps}{t}  & = f_1(\ueps, \veps, \weps) &\text{for}\ \vardefPDE \\
\delta \del{\veps}{t}  & = f_2(\ueps, \veps, \weps) &\text{for}\ \vardefPDE \\
 \del{\weps}{t}  & = D \dell{v}{x} + g(\ueps, \veps, \weps) &\text{for}\ \vardefPDE\\
 \del{\weps}{x}&(t, 0) = \del{v}{x}(t, l) = 0 &\text{for}\ t>0.
\end{aligned}
\right.
\label{A2.15}
\end{align}

\begin{cor}[Stability for vector-valued $u$ with a small parameter] \label{thm:stabC-3c}  Let $ I =(0, l) \subset \R$ be a bounded interval. Under Assumption 2.1, consider a system of type {\rm \eqref{A2.15}}. Let $(\tilde u_1, \tilde u_2, \tilde v)$ be a steady state with finitely many discontinuities of $(\tilde u_1, \tilde u_2)$ in $x$.  Let $A(x)$ be the matrix defined by {\rm \eqref{eqn:JacobianA}}.
Let $A_{ij}$ denote the matrix resulting from omitting the $i$-th row and the $j$-th column, i.e., $A_{ij}=(a_{kl})_{k\neq i, l\neq j}$. Assume that the following inequalities are satisfied for all $x\in \overline{I}$:
\begin{align}
& a_{22}(x)  \leq - k_1 < 0, & \label{A2.16} \\
& a_{33}(x)  \leq - k_2 < 0, & \label{A2.17} \\
& \det A_{11}(x)  >0, & \label{A2.18}\\
& \det A_{33}(x)  > -k_3 a_{22}(x) >0, & \label{A2.19} \\
& \det A(x)  < 0, & \label{A2.20} 
\end{align}
where $k_1,\,k_2,\,k_3$ are positive constants independent of $x$.
Then there exists a positive constant $\delta^*$ such that  the steady state $(\tilde u_1, \tilde u_2, \tilde v) $ is $(\varepsilon_0, A)$-stable for some positive constants $\varepsilon_0$ and $A$, provided that $0<\delta<\delta^*$.
\end{cor}

\subsection{Quasi-steady state approximation} $\,$ \\
By the quasi-steady state approximation of \eqref{A2.15}, we understand system \eqref{A2.15} for $\delta=0$, i.e., the following system of differential-algebraic equations:
\begin{equation} \label{A2.21}\left\{\ 
\begin{aligned}
 & \del{u_1^0}{t} = f_1(u_1^0,u_2^0,v^0),              &\vardefODE\\
& 0 = f_2(u_1^0,u_2^0,v^0),                             &\vardefODE\\
 &\del{v^0}{t} = D \dell{v^0}{x} + g(u_1^0,u_2^0,v^0),\qquad &\vardefPDE\\
 &\del{v^0}{x}(t, 0)=\del{v^0}{x}(t, l)=0 & \text{for} \ t>0.
\end{aligned}
\right.
\end{equation}
Assuming that $f_2(u_1^0,u_2^0,v^0)=0$ has an isolated solution $u_2^0= u_2^*(u_1^0,v^0)$, we can re-write \eqref{A2.21} as
\begin{equation}\label{A2.22} \left\{ \ 
\begin{aligned}
 & \del{u_1^0}{t} = f_1(u_1^0,u_2^*(u_1^0,v^0),v^0),                             &\vardefODE\\
& \del{v^0}{t} = D \dell{v^0}{x} + g(u_1^0,u_2^*(u_1^0,u_2^0),v^0),\qquad &\vardefPDE\\
 & \del{v^0}{x}(t, 0)=\del{v^0}{x}(t,l)=0 & \text{for}\ t>0.
\end{aligned}
\right.
\end{equation}

It follows immediately that any steady state of the reduced system \eqref{A2.22} gives rise to a steady state of \eqref{A2.21}. Hence, we want to investigate the following questions: 
\begin{enumerate}
 \item Does system  \eqref{A2.15} exhibit DDI if its quasi-steady state reduction \eqref{A2.22} does?
 \item Are spatially inhomogeneous steady states of model \eqref{A2.15} $(\varepsilon_0,A)$-stable if they are $(\varepsilon_0,A)$-stable steady states of model \eqref{A2.22}?
\end{enumerate}

The first question is addressed in Section 4, see Proposition \ref{propDDIImply}. For $\partial f_2/\partial u_2 \, \leq -c <0$, it turns out that DDI of the reduced system implies DDI of the unreduced system for sufficiently small $\delta$. Here, we address the second question on the `transfer' of stability:

\begin{theorem}[Stability transfer]\label{lem:stabtrans} Consider a steady state $(\tilde{u}_1,\tilde{u}_2,\tilde{v})$ of the system
\begin{equation} \label{A2.23} \begin{aligned}
 \del{\ueps}{t} =& f_1(\ueps,\veps,\weps),                               &\vardefODE\\
 \delta \del{\veps}{t} =& f_2(\ueps,\veps,\weps),                        &\vardefODE\\
 \del{\weps}{t} =& D \dell{\weps}{x} + g(\ueps,\veps,\weps),\qquad &\vardefPDE\\
 (\ueps(0,x),\veps(0,x),&\weps(0,x)) \in C(\overline{\text{$I$}})^2 \times (C^2(\text{$I$})\cap C(\overline{\text{$I$}})), 
\end{aligned}
\end{equation}
with homogeneous Neumann boundary conditions for $v$
and its quasi-steady state reduction \begin{equation} \label{A2.24}  \begin{aligned}
 \del{u_1^0}{t} =& f_1(u_1^0,u_2^*(u_1^0,v^0),v^0),           &\vardefODE\\
 \del{v^0}{t} =& D \dell{v^0}{x} + g(u_1^0,u_2^*(u_1^0,v^0),v^0),\qquad &x \in \text{$I$}, t>0.
\end{aligned}
\end{equation}
Denote the Jacobian matrix of the kinetic system of \eqref{A2.23} with $\delta=1$ at $(\tilde{u}_1,\tilde{u}_2,\tilde{v})$ by $A=(a_{ij})$ and denote the Jacobian matrix of the kinetic system of \eqref{A2.24} at $(\tilde{u}_1,\tilde{v})$ by $\MB=(\eb_{ij})$. Suppose that
\begin{enumerate}
 \item[{\rm (1)}] system \eqref{A2.24} satisfies the conditions of Theorem \ref{thm:stabC},
 \item[{\rm (2)}] system \eqref{A2.23} satisfies Assumption \ref{assump},
 \item[{\rm (3)}] $a_{22},a_{33}\leq -c < 0$,
\end{enumerate}
where $c$ is a positive constant independent of $x$.
Then, there exists a $\delta^*>0$ such that $(\tilde{u}_1,v^*(\tilde{u}_1,\tilde{v}),\tilde{v})$ is an $(\epsilon_0,A)$-stable steady state of system \eqref{A2.23} for all $\delta<\delta^*$.\end{theorem}

\section{Application to a system exhibiting DDI and hysteresis}\label{sec:ex}
In this section, we present a model exhibiting the coexistence of \textit{DDI} and \textit{hysteresis}. First, we propose a model consisting of two ordinary differential equations coupled to one reaction-diffusion equation. We reduce this model using quasi-steady state approximation. Then, we rescale the reduced model and show that it satisfies the conditions of Theorem \ref{thm:stabC}. Finally, based on Theorem \ref{lem:stabtrans}, we show that the stability of steady states of the reduced model implies stability of corresponding patterns in the original model.\par

Our example takes the form
\begin{equation}\label{3er-sys}
\begin{aligned}
 \del{\ueps}{\hat{t}} &= -\nu_1 \ueps - \beta \ueps \weps + \theta_1 \frac{(\ueps)^2}{1+\kappa (\ueps)^2} + \alpha \veps, &\vardefODEhat\\
 \delta \del{\veps}{\hat{t}} &= -\nu_2 \veps + \beta \ueps \weps - \alpha \veps, &\vardefODEhat\\
 \del{\weps}{\hat{t}} &= \gamma \dell{\weps}{x} -\nu_3 \weps - \beta \ueps \weps + \theta_2 \frac{(\ueps)^2}{1+\kappa (\ueps)^2} + \alpha \veps,\qquad &\vardefPDEhat
\end{aligned}
\end{equation}
supplemented with homogeneous Neumann boundary conditions and initial functions $(\ueps(0),\veps(0),\weps(0)) \in C(\overline{\text{$I$}})^2 \times (C^2(\text{$I$})\cap C(\overline{\text{$I$}}))$. Here, $\alpha$, $\beta$, $\gamma$, $\theta_1$, $\theta_2$, $\nu_1$, $\nu_2$, $\nu_3$ are positive constants. Without loss of generality, we consider $\text{$I$}=(0,1)$ in order to simplify analysis.

The model is a modification of the receptor-based model proposed in \cite{Marciniak03,MCP}. Numerical investigations for $\kappa=0$ show spike solutions blowing up in infinite time, see \cite{HM13}. Introducing $\kappa > 0$ reflects saturation of {\it de novo} production. The modification assures uniform boundedness of the model solutions.

Variables  $u_1$ and $u_2$ describe cell surface receptors at free and occupied states, respectively. The state of receptors is being changed through binding to diffusive ligands, denoted by $v$, and a reverse process of dissociation of bound receptors. Free receptors and ligands are produced {\it de novo}, what is described by a function with saturation effect for $\kappa>0$.

For $\delta=0$, \eqref{3er-sys} reduces to  a system of differential-algebraic equations. In this case, the second equation of \eqref{3er-sys} can be solved uniquely for $u_2^0$ leading to a system of one ordinary differential equation coupled to one reaction-diffusion equation. This reduction leads to simplification of analytic investigation of the model. We show that the limit $\delta \rightarrow 0$ is regular in the sense of stability-preservation of steady states. 

The quasi-steady state approximation, {i.e. $\delta = 0$}, reads
\begin{equation}
\begin{aligned}
 \del{u_1^0}{\hat{t}} &= -\nu_1 u_1^0 - \beta \frac{\nu_2}{\nu_2+\alpha} u_1^0 v^0 + \theta_1 \frac{(u_1^0)^2}{1+\kappa (u_1^0)^2}, &\vardefODEhat\\
 \del{v^0}{\hat{t}} &= \gamma \dell{v^0}{x} -\nu_3 v^0 - \beta \frac{\nu_2}{\nu_2+\alpha} u_1^0 v^0 + \theta_2 \frac{(u_1^0)^2}{1+\kappa (u_1^0)^2}, \qquad &\vardefPDEhat
\end{aligned}
\end{equation}
After rescaling $u_1^0,v^0,\hat{t}$, we obtain
\begin{align}
 \del{u}{t} &= - u - uv + m_1 \frac{u^2}{1+k u^2}, &\vardefODE \label{eq1}\\
 \del{v}{t} &= D\dell{v}{x} -\mu_3 v - uv + m_2 \frac{u^2}{1+k u^2},       &\vardefPDE \label{eq2}\\
 &\del{v}{x}(t, 0)=\del{v}{x}(t,1)=0,\label{bc}\\
 (u(0,x),v(0,x))&=(u_0(x),v_0(x)) \in C(\overline{\text{$I$}})\times (C^2(\text{$I$})\cap C(\overline{\text{$I$}})),\label{ic}
\end{align}
with positive initial data.

To simplify notation, we define 
\begin{equation*}
\begin{aligned}
 f_r(u,v) &:= -u - uv + m_1 \frac{u^2}{1+ku^2}, \\
 g_r(u,v) &:= -\mu_3 v - uv + m_2 \frac{u^2}{1+ku^2}. 
\end{aligned}
\end{equation*}
\noindent
We show that the model has a spatially homogeneous steady state that is stable to spatially homogeneous perturbation and is destabilized due to autocatalysis in dynamics of non-diffusive component. Due to the nonlinear nature of the reaction terms, there exists a `far from equilibrium' regime, where the species become self-inhibitory, stabilizing steady states with jump discontinuity. 

Our first step is to show that models \eqref{eq1}-\eqref{bc} and \eqref{3er-sys} satisfy Assumption \ref{assump}, i.e., the solutions are uniformly bounded for nonnegative initial functions.  Numerical simulations of the model solutions are presented in Figure \ref{fig:illu1}. We formulate the results in the remainder of this section and defer their proofs to Section \ref{sec:exproof}.

\subsection{Global existence and uniform boundedness} $\,$ \\
For model \eqref{eq1}-\eqref{ic}, global existence of a solution in $C^1((0,\infty);C(\overline{\text{$I$}})\times (C^2(\text{$I$})\cap C(\overline{\text{$I$}})))$ is assured within the framework of invariant rectangles, see, e.g., \cite{Smoller}. The solution is uniformly bounded, i.e., \begin{lemma}[Uniform boundednes of solutions of system \eqref{eq1}-\eqref{bc}]\label{lem:trivstab} {\rm (i)} Assume that $m_1 \geq 2 \sqrt{k}$. Then, for positive initial functions, the set
\begin{equation*}
 \left(0, \frac{1}{2k} \big( {m_1} + \sqrt{ ( {m_1} )^2 - 4k}\, \big) \right)
\times \left( 0, \frac{m_2}{k \mu_3}\right)
\end{equation*}
is an attracting invariant rectangle of model \eqref{eq1}-\eqref{bc}. {\rm (ii)} Assume that $m_1<2\sqrt{k}$. For initial functions $(u_0(x), v_0(x))$ satisfying $0<u_0(x)$, $0<v_0(x)< m_2/(k\mu_3)$, it holds that $u(t,x)$ is monotone decreasing in $t \in (0, \infty)$ and
\begin{equation}
(u(t,x), v(t,x))\rightarrow (0, 0) \qquad \text{uniformly on $\overline{I}$, as $t \to+\infty$}.
\end{equation}
\end{lemma}

Also for model \eqref{3er-sys}, global existence of a solution in $C^1((0,\infty);C(\overline{\text{$I$}})^2\times (C^2(\text{$I$})\cap C(\overline{\text{$I$}})))$ follows from the theory of invariant rectangles.
\begin{lemma}[Uniform boundednes of solutions of system \eqref{3er-sys}]\label{lem:3-ub} For positive initial functions, the set
\begin{align*} 
 &\left(0,\frac{1}{\kappa}\left( \frac{\theta_1}{2 \nu_1} + \sqrt{\left(\frac{\theta_1}{2\nu_1}\right)^2-\kappa}\right)\right)
\times
\left(0,\frac{\theta_1}{\kappa \delta \min(\nu_1,\nu_2)}\right)  \\
& \qquad \times 
\left(0,\frac{1}{\nu_3\kappa}\left( \theta_2 + \frac{\alpha \theta_1}{\delta \alpha \min(\nu_1,\nu_2)} \right)\right)
\end{align*}
is an attracting invariant rectangle of model \eqref{3er-sys}. Here, in the case $\theta_1<2\nu_1\sqrt{\kappa}$, the expression  $\theta_1/(2\nu_1) + \sqrt{(\theta_1/(2\nu_1))^2-\kappa}$ is interpreted as an arbitrary positive number. For $\theta_1<2\nu_1 \sqrt{\kappa}$, it holds that
\begin{equation*}
(\ueps(t,x),\veps(t,x),\weps(t,x))\rightarrow (0,0,0) \qquad \text{uniformly on $\overline{I}$, as $t\to+\infty$}.
\end{equation*}
\end{lemma}

\subsection{Existence of steady states} $\,$ \\
By definition of quasi-stationary reduction and linear rescaling, there exists a one-to-one mapping between the sets of steady states of model \eqref{eq1}-\eqref{bc} and of model \eqref{3er-sys}. Therefore, we may focus on existence of steady states of model \eqref{eq1}-\eqref{bc}. 
A function $(u,v)$ is a steady state of model \eqref{eq1}-\eqref{bc} only if 
\begin{equation} \del{u}{t}=f_r(u,v)=0, \end{equation}
is satisfied. Note that  $u(v):=\{ u \in \mathbb{R}_{\geq 0} \mid f_r(u,v)=0 \}$ may consist of more than one element. We characterize different possible branches of solutions for given $v$:

\begin{lemma} \label{stst:f=0} For given $v \in \mathbb{R}_{\geq 0}$,
\begin{equation}
 f_r(u,v)=0
\end{equation}
has three solutions
\begin{align}
 u_0(v)  &= 0,\\
 u_-(v)  &= \frac{1}{k} \left(\frac{m_1}{2(1+v)}-\sqrt{  \left(\frac{m_1}{2(1+v)}\right)^2  -k} \right),\\
 u_+(v)  &=\frac{1}{k} \left(\frac{m_1}{2(1+v)}+\sqrt{  \left(\frac{m_1}{2(1+v)}\right)^2  -k} \right). 
\end{align}
\end{lemma}

For $v \geq m_1/(2 \sqrt{k})-1$, $u_-(v)$ and $u_+(v)$ are not real-valued. We therefore define
\begin{equation}
 v_r := \frac{m_1}{2 \sqrt{k}}-1,
\end{equation}
which is the largest $v\in \mathbb{R}$ such that $u_{\pm}(v)$ is {\bf r}eal-valued. Note that $v_r>0$ is equivalent to $m_1>2 \sqrt{k}$. Therefore, $v_r\leq 0$ implies $(u(t,x),v(t,x))\rightarrow (0,0)$ uniformly as $t \to \infty$, due to Lemma \ref{lem:trivstab}.

We show in Lemma \ref{lem:entries} that the steady states have different stability properties depending on the branch. 
\begin{definition}
We say that a steady state $(u,v)$ is of \textit{type $u_-$}(respectively $u_+$, $u_0$) at $x \in \overline{\text{$I$}}$ if $u(x)=u_-(v(x))$ (respectively $u(x)=u_+(v(x))$, $u(x)=u_0(v(x))=0$).
\end{definition} \par
The existence of steady states can be summarized in the following lemma, illustrated in Figure
\ref{fig:homStSt}:
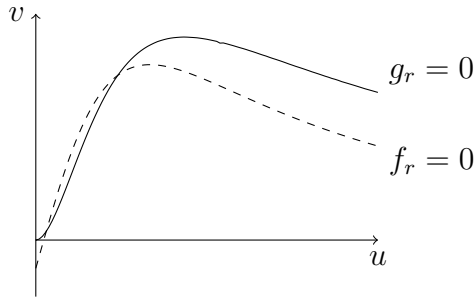
\begin{figure}[ht]
  \centering
  \begin{tikzpicture}[domain=0:3,scale=1.5]
    \draw[very thin,color=gray] (0.0,0.0) grid (0,0);
    \draw[->] (0,0) -- (3,0) node[below] {$u$};
    \draw[->] (0,-0.5) -- (0,2) node[left] {$v$};
    \draw[samples=105,domain=0:3, dashed] plot(\x,   {0.25*(-1+1.44*10*\x/(1+0.01*100*\x*\x))});
    \draw[samples=105,domain=0:3] plot(\x,   {0.25*2/(4.51061+10*\x) * 10*\x * 10*\x / (1+0.01*100*\x*\x)});
    \draw (3,1.5) node[right] {$g_r=0$};
    \draw (3,0.7) node[right] {$f_r=0$};
\end{tikzpicture}
   \caption{Plot of the nullclines of $f_r(u,v)=-(1+v)u+m_1(u^2/(1+ku^2))$ for $u\neq 0$ and $g_r(u,v)=-(\mu_3+u)v+m_2 (u^2/(1+ku^2))$.}
  \label{fig:homStSt}
\end{figure}

\begin{lemma}[Existence of spatially homogeneous steady states]\label{lem:homStSt} $\,$\\
\begin{enumerate}[{\rm (1)}]
\item For arbitrary positive parameters, the trivial solution
$(u,v)=(0,0)$ is a spatially homogeneous steady state of model \eqref{eq1}-\eqref{bc}.\\
\item Let $m_1<m_2$ and 
\begin{equation*} \mu_3 > \frac{1}{m_1}\left(\frac{2 m_2 - m_1}{m_1} + 2 \sqrt{\left(\frac{m_2}{m_1}\right)^2-\frac{m_2}{m_1}}\right).\end{equation*}
Then, there exists a positive $k_1^*$ such that for all $k<k_1^*$, system \eqref{eq1}-\eqref{bc} has exactly two positive homogeneous steady states. Moreover, if $k<\min(k_1^*,((m_2-m_1)/(m_1 \mu_3))^2)$, both nontrivial positive homogeneous steady states are of type $u_-$.
\end{enumerate}
\end{lemma}
We show in Lemma \ref{lem:entries} that the conditions of Theorem \ref{thm:stabC} are never satisfied if $u$ is of type $u_-$ at some $x \in \overline{\text{$I$}}$, but always satisfied if $u$ is of type $u_0$ or $u_+$ for all $x\in \overline{\text{$I$}}$. Hence, the latter type of solutions is constructed.

\begin{lemma}[Existence of steady states with jump discontinuities]\label{lem:disc}  Assume $2 \sqrt{k}<m_1<m_2$. Then, there exist infinitely many steady states $(\tilde{u},\tilde{v}) \in (L^{\infty}(\text{$I$})\cap BV(\overline{\text{$I$}})) \times (H^2(\text{$I$})\cap C^1(\overline{\text{$I$}}))$ which are of type $u_+$ or $u_0$ for all $x \in \overline{\text{$I$}}$.\\
\end{lemma}

\subsection{Stability of steady states for the quasi-steady state approximation} $\,$ \\
Following the same principle as in Section \ref{sec:stab}, we prove stability of steady states of the quasi-steady state approximation \eqref{eq1}-\eqref{bc}. In the next subsection, we conclude stability of steady states of model \eqref{3er-sys}. 
In order to apply Theorem \ref{thm:stabC}, we investigate how the entries of Jacobian matrix at a steady state depend on the branch:
\begin{lemma}[Characterization of Jacobian matrix for model \eqref{eq1}-\eqref{bc}]\label{lem:entries} Assume $2 \sqrt{k}<m_1<m_2$. For a pair of nonnegative functions $(u(x), v(x))$, denote the Jacobian matrix of the kinetic system of \eqref{eq1}-\eqref{eq2} by
\begin{equation}
 \MB(x)=\MB(u(x),v(x))=\begin{pmatrix} \del{f_r}{u}(x) & \del{f_r}{v}(x) \\ \del{g_r}{u}(x) & \del{g_r}{v}(x) \end{pmatrix} =
\begin{pmatrix} \eb_{11}(x) & \eb_{12}(x) \\ \eb_{21}(x) & \eb_{22}(x) \end{pmatrix}.
\end{equation}
Then, it holds for $0 \leq v(x) < m_1/(2\sqrt{k})-1=:v_r$ that
\begin{equation}
\begin{aligned}
 \eb_{11}(x) \left\{\begin{array}{cl} 
               >0, & \text{if}\ u(x)=u_-(v(x)) \\ 
               <0, & \text{if}\ u(x)=u_+(v(x)) \\
               <0, & \text{if}\ u(x)=u_0
            \end{array}\right. ,
\quad 
 \eb_{12}(x) \left\{\begin{array}{cl} 
               <0, & \text{if}\ u(x)=u_-(v(x)) \\ 
               <0, & \text{if}\ u(x)=u_+(v(x)) \\
               =0, & \text{if}\ u(x)=u_0
            \end{array}\right. , \\
 \eb_{21}(x) \left\{\begin{array}{cl} 
               >0, & \text{if}\ u(x)=u_-(v(x)) \\ 
               <0, & \text{if}\ u(x)=u_0
            \end{array}\right. ,
\quad 
  \eb_{22}(x) \left\{\begin{array}{cl} 
               <0, & \text{if}\ u(x)=u_-(v(x)) \\ 
               <0, & \text{if}\ u(x)=u_+(v(x)) \\
               <0, & \text{if}\ u(x)=u_0
            \end{array}\right. .
\end{aligned}
\end{equation}
Moreover, there exists $0<v^*<v_r$ such that, for $u(x)=u_+(v(x))$, it holds that
\begin{equation}
 \eb_{21}(x) \left\{\begin{array}{cl} >0 & \text{if} \ v(x) \in (v^*,v_r), \\ =0 & \text{if} \ v(x)=v^*.\end{array}\right.
\end{equation}
For an isolated solution $u(v)$ of $f_r(u,v)=0$ (i.e., $\partial f_r/ \partial u \cdot  \partial f_r/ \partial v \neq 0 $ at $(u(v),v)$), it holds that
\begin{equation}
 \operatorname{det}(\MB(u(v),v))=-\eb_{12} \left(\frac{d}{dv}u(v)\right)^{-1} \frac{d}{dv}g_r(u(v),v),
\end{equation}
and, in particular, for $u(x)=u_+(v(x))$,
\begin{equation}
  \operatorname{det}(\MB(u_+(v),v))>0.
\end{equation}
\end{lemma} 

Applying Theorem \ref{thm:stabC} yields stability of the steady states constructed in Lemma \ref{lem:disc}.
\begin{cor}[Stability of solutions with jump discontinuity of model \eqref{eq1}-\eqref{bc}]\label{lem:stabRD} $\,$ Assume $D>0$ and $2 \sqrt{k}<m_1<m_2$. Then,
there exist infinitely many $(\varepsilon_0,A)$-stable steady states with jump discontinuity of model \eqref{eq1}-\eqref{bc}, which are for every $x \in \overline{I}$ either of type $u_0$ or $u_+$. 
\end{cor}
To obtain DDI, it is left to prove that at least one of the spatially homogeneous steady states of type $u_-$ is stable under spatially constant perturbations, i.e.~a stable steady state of the kinetic system.
Instability of all spatially homogeneous steady states of type $u_-$ follows from \cite{MCKS13} due to autocatalysis of the non-diffusive species.
Existence of exactly two homogeneous steady states of type $u_-$ results from Lemma \ref{lem:homStSt}.

\begin{lemma}[Diffusion-Driven Instability in model \eqref{eq1}-\eqref{bc}]\label{lem:stabKin}  Let $m_1<\min(m_2,\sqrt{m_2})$. If there exist two homogeneous steady states $(u_-(v_-),v_-)$, $(u_-(v_+),v_+)$ with $v_-<v_+$, then there exists $k_2^*>0$ such that for all $k<k_2^*$
\begin{enumerate}[{\rm (1)}]
\setlength{\leftskip}{1.5pc}
 \item $(u_-(v_-),v_-)$ is unstable to spatially homogeneous perturbations,
 \item $(u_-(v_+),v_+)$ is stable to spatially homogeneous perturbations.
\end{enumerate}
 Furthermore, all spatially homogeneous steady states of type $u_-$ are unstable for $k \geq 0$ and $D>0$.
\end{lemma}
We summarize our results in the following 
\begin{theorem}[Coexistence of DDI and Hysteresis for model \eqref{eq1}-\eqref{bc}] Let $k_1^*$ and $k_2^*$ be defined as in Lemmas \ref{lem:homStSt} and \ref{lem:stabKin}. Under conditions
\begin{equation}
\begin{aligned}
 m_1,m_2,k,\mu_3 &> 0,\\
 m_1   &< \min(m_2,\sqrt{m_2}),\\
 \mu_3 &> \frac{1}{m_1} \left(\frac{2m_2-m_1}{m_1}+2 \sqrt{\left(\frac{m_2}{m_1}\right)^2 - \frac{m_2}{m_1}}\right),\\
\end{aligned}
\end{equation}
there exists a $k^*$ in the interval
\begin{equation}
 0<k^* \leq \min\left(k_1^*,k_2^*,\left(\frac{m_2-m_1}{m_1 \mu_3}\right)^2, \frac{m_1^2}{4}\right)
\end{equation}
such that, for all $0<k<k^*$, the following $(1)$--$(5)$ hold true:
\begin{enumerate}[{\rm (1)}]
\setlength{\leftskip}{1pc}
 \item system \eqref{eq1}-\eqref{bc} has exactly two strictly positive spatially homogeneous steady states $(u_-(v_-), v_-)$ and $(u_-(v_+), v_+)$ with $v_-<v_+$.
 \item $(u_-(v_-), v_-)$ is an unstable equilibrium of the kinetic system for \eqref{eq1}-\eqref{bc}, hence is an unstable steady state of system \eqref{eq1}-\eqref{bc}.  \item $(u_-(v_+), v_+)$ is a stable steady state of the kinetic system of \eqref{eq1}-\eqref{bc} and an unstable steady state of \eqref{eq1}-\eqref{bc}.
 \item $(0,0)$ is a stable steady state of the kinetic system for \eqref{eq1}-\eqref{bc} and of the original system \eqref{eq1}-\eqref{bc}.
 \item system \eqref{eq1}-\eqref{bc} has infinitely many $(\varepsilon_0,A)$-stable steady states with jump discontinuity, which are at $x\in \overline{\text{$I$}}$ of type $u_+$ or $u_0$.
\end{enumerate}
\end{theorem}
A numerically obtained solution, using the finite element library \textit{deal.ii} \cite{deal.ii}, is presented in Figure \ref{fig:illu1}. 
We observe destabilization of the spatially homogeneous steady state and convergence towards a pattern with jump discontinuity. Compartment $u$ is discretized with cell-wise constant
finite elements and $v$ is discretized with cell-wise linear, globally continuous finite elements. The formulation is the weak formulation.
The time-stepping scheme is Crank-Nicholson.\\ 
\begin{minipage}[t]{36em}
\begin{tabular}{p{18em}p{18em}}
   \begingroup
  \makeatletter
  \providecommand\color[2][]{    \GenericError{(gnuplot) \space\space\space\@spaces}{      Package color not loaded in conjunction with
      terminal option `colourtext'    }{See the gnuplot documentation for explanation.    }{Either use 'blacktext' in gnuplot or load the package
      color.sty in LaTeX.}    \renewcommand\color[2][]{}  }  \providecommand\includegraphics[2][]{    \GenericError{(gnuplot) \space\space\space\@spaces}{      Package graphicx or graphics not loaded    }{See the gnuplot documentation for explanation.    }{The gnuplot epslatex terminal needs graphicx.sty or graphics.sty.}    \renewcommand\includegraphics[2][]{}  }  \providecommand\rotatebox[2]{#2}  \@ifundefined{ifGPcolor}{    \newif\ifGPcolor
    \GPcolorfalse
  }{}  \@ifundefined{ifGPblacktext}{    \newif\ifGPblacktext
    \GPblacktexttrue
  }{}    \let\gplgaddtomacro\g@addto@macro
    \gdef\gplbacktext{}  \gdef\gplfronttext{}  \makeatother
  \ifGPblacktext
        \def\colorrgb#1{}    \def\colorgray#1{}  \else
        \ifGPcolor
      \def\colorrgb#1{\color[rgb]{#1}}      \def\colorgray#1{\color[gray]{#1}}      \expandafter\def\csname LTw\endcsname{\color{white}}      \expandafter\def\csname LTb\endcsname{\color{black}}      \expandafter\def\csname LTa\endcsname{\color{black}}      \expandafter\def\csname LT0\endcsname{\color[rgb]{1,0,0}}      \expandafter\def\csname LT1\endcsname{\color[rgb]{0,1,0}}      \expandafter\def\csname LT2\endcsname{\color[rgb]{0,0,1}}      \expandafter\def\csname LT3\endcsname{\color[rgb]{1,0,1}}      \expandafter\def\csname LT4\endcsname{\color[rgb]{0,1,1}}      \expandafter\def\csname LT5\endcsname{\color[rgb]{1,1,0}}      \expandafter\def\csname LT6\endcsname{\color[rgb]{0,0,0}}      \expandafter\def\csname LT7\endcsname{\color[rgb]{1,0.3,0}}      \expandafter\def\csname LT8\endcsname{\color[rgb]{0.5,0.5,0.5}}    \else
            \def\colorrgb#1{\color{black}}      \def\colorgray#1{\color[gray]{#1}}      \expandafter\def\csname LTw\endcsname{\color{white}}      \expandafter\def\csname LTb\endcsname{\color{black}}      \expandafter\def\csname LTa\endcsname{\color{black}}      \expandafter\def\csname LT0\endcsname{\color{black}}      \expandafter\def\csname LT1\endcsname{\color{black}}      \expandafter\def\csname LT2\endcsname{\color{black}}      \expandafter\def\csname LT3\endcsname{\color{black}}      \expandafter\def\csname LT4\endcsname{\color{black}}      \expandafter\def\csname LT5\endcsname{\color{black}}      \expandafter\def\csname LT6\endcsname{\color{black}}      \expandafter\def\csname LT7\endcsname{\color{black}}      \expandafter\def\csname LT8\endcsname{\color{black}}    \fi
  \fi
  \setlength{\unitlength}{0.0500bp}  \begin{picture}(5040.00,3528.00)    \gplgaddtomacro\gplbacktext{    }    \gplgaddtomacro\gplfronttext{      \csname LTb\endcsname      \put(594,1142){\makebox(0,0)[r]{\strut{} 0}}      \put(2248,643){\makebox(0,0)[r]{\strut{} 14}}      \put(1191,819){\makebox(0,0){\strut{}$t$}}      \put(2613,588){\makebox(0,0){\strut{} 0}}      \put(4385,1173){\makebox(0,0){\strut{} 1}}      \put(3849,819){\makebox(0,0){\strut{}$x$}}      \put(622,1656){\makebox(0,0)[r]{\strut{} 0}}      \put(622,2528){\makebox(0,0)[r]{\strut{} 16}}      \put(891,819){\vector(3,-1){500}}
    }    \gplbacktext
    \put(0,0){\includegraphics{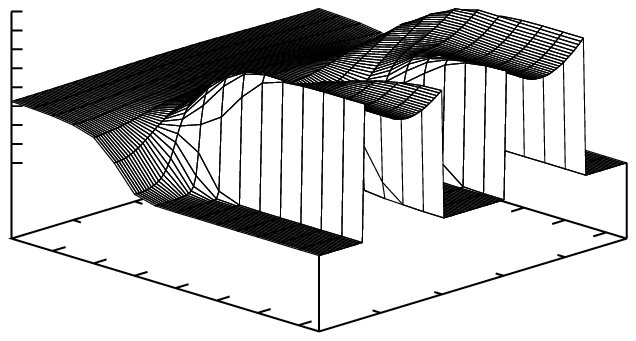}}    \gplfronttext
  \end{picture}\endgroup
 &
   \begingroup
  \makeatletter
  \providecommand\color[2][]{    \GenericError{(gnuplot) \space\space\space\@spaces}{      Package color not loaded in conjunction with
      terminal option `colourtext'    }{See the gnuplot documentation for explanation.    }{Either use 'blacktext' in gnuplot or load the package
      color.sty in LaTeX.}    \renewcommand\color[2][]{}  }  \providecommand\includegraphics[2][]{    \GenericError{(gnuplot) \space\space\space\@spaces}{      Package graphicx or graphics not loaded    }{See the gnuplot documentation for explanation.    }{The gnuplot epslatex terminal needs graphicx.sty or graphics.sty.}    \renewcommand\includegraphics[2][]{}  }  \providecommand\rotatebox[2]{#2}  \@ifundefined{ifGPcolor}{    \newif\ifGPcolor
    \GPcolorfalse
  }{}  \@ifundefined{ifGPblacktext}{    \newif\ifGPblacktext
    \GPblacktexttrue
  }{}    \let\gplgaddtomacro\g@addto@macro
    \gdef\gplbacktext{}  \gdef\gplfronttext{}  \makeatother
  \ifGPblacktext
        \def\colorrgb#1{}    \def\colorgray#1{}  \else
        \ifGPcolor
      \def\colorrgb#1{\color[rgb]{#1}}      \def\colorgray#1{\color[gray]{#1}}      \expandafter\def\csname LTw\endcsname{\color{white}}      \expandafter\def\csname LTb\endcsname{\color{black}}      \expandafter\def\csname LTa\endcsname{\color{black}}      \expandafter\def\csname LT0\endcsname{\color[rgb]{1,0,0}}      \expandafter\def\csname LT1\endcsname{\color[rgb]{0,1,0}}      \expandafter\def\csname LT2\endcsname{\color[rgb]{0,0,1}}      \expandafter\def\csname LT3\endcsname{\color[rgb]{1,0,1}}      \expandafter\def\csname LT4\endcsname{\color[rgb]{0,1,1}}      \expandafter\def\csname LT5\endcsname{\color[rgb]{1,1,0}}      \expandafter\def\csname LT6\endcsname{\color[rgb]{0,0,0}}      \expandafter\def\csname LT7\endcsname{\color[rgb]{1,0.3,0}}      \expandafter\def\csname LT8\endcsname{\color[rgb]{0.5,0.5,0.5}}    \else
            \def\colorrgb#1{\color{black}}      \def\colorgray#1{\color[gray]{#1}}      \expandafter\def\csname LTw\endcsname{\color{white}}      \expandafter\def\csname LTb\endcsname{\color{black}}      \expandafter\def\csname LTa\endcsname{\color{black}}      \expandafter\def\csname LT0\endcsname{\color{black}}      \expandafter\def\csname LT1\endcsname{\color{black}}      \expandafter\def\csname LT2\endcsname{\color{black}}      \expandafter\def\csname LT3\endcsname{\color{black}}      \expandafter\def\csname LT4\endcsname{\color{black}}      \expandafter\def\csname LT5\endcsname{\color{black}}      \expandafter\def\csname LT6\endcsname{\color{black}}      \expandafter\def\csname LT7\endcsname{\color{black}}      \expandafter\def\csname LT8\endcsname{\color{black}}    \fi
  \fi
  \setlength{\unitlength}{0.0500bp}  \begin{picture}(5040.00,3528.00)    \gplgaddtomacro\gplbacktext{    }    \gplgaddtomacro\gplfronttext{      \csname LTb\endcsname      \put(594,1142){\makebox(0,0)[r]{\strut{} 0}}      \put(2011,714){\makebox(0,0)[r]{\strut{} 12}}      \put(1191,819){\makebox(0,0){\strut{}$t$}}      \put(2613,588){\makebox(0,0){\strut{} 0}}      \put(4385,1173){\makebox(0,0){\strut{} 1}}      \put(3849,819){\makebox(0,0){\strut{}$x$}}      \put(622,1656){\makebox(0,0)[r]{\strut{} 4.8}}      \put(622,2528){\makebox(0,0)[r]{\strut{} 6.2}}      \put(891,819){\vector(3,-1){500}}    }    \gplbacktext
    \put(0,0){\includegraphics{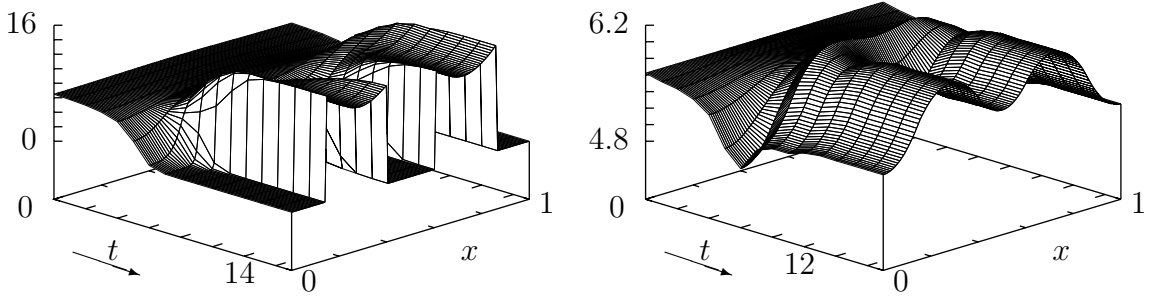}}    \gplfronttext
  \end{picture}\endgroup
 \end{tabular}
\captionof{figure}{Numerically obtained solution to model \eqref{eq1}-\eqref{bc} for parameters $m_1 =1.44, m_2 = 2,\mu_3 \approx 4.1, k=0.01, D=1$. We observe convergence towards a steady state with jump discontinuity. Left: Non-diffusive component $u$. Right: Diffusive component $v$.}
\label{fig:illu1}
\end{minipage}\\

\subsection{Stability of steady states for the unreduced model} $\,$ \\
In this section we show that model \eqref{3er-sys} satisfies the conditions of the `stability transfer' Theorem \ref{lem:stabtrans}.
\begin{lemma}\label{lem:3er-sys}
Consider model \eqref{3er-sys}. Define $u_2^*(u_1,v):= \beta u_1 v / (\nu_2+\alpha)$. Under conditions
\begin{align}
 0 &< \nu_1, \nu_3, \theta_1, \theta_2, \nu_2, \alpha, \beta, \gamma,\\
 \theta_1 &< \min\left(\theta_2, \sqrt{\theta_2 \beta \frac{\nu_2}{\alpha+\nu_2}}\right),\\
 \nu_3 &> 2 \nu_1 \beta \frac{\nu_2}{\alpha+\nu_2} \frac{1}{\theta_1} \left(\frac{2 \theta_2-\theta_1}{\theta_1}+2 \sqrt{\left(\frac{\theta_2}{\theta_1}\right)^2-\frac{\theta_2}{\theta_1}}\right),
\end{align}
there exists a constant $\kappa^*$ in the interval
\begin{equation} 
  0<\kappa^*<\min\left( \left( \beta \frac{\theta_2-\theta_1}{\nu_1} \frac{\nu_2}{\alpha+\nu_2}\right)^2, \left(\frac{\theta_1}{2\nu_1}\right)^2\right)
\end{equation} 
and a positive constant $\delta^*$ such that for all $0 \leq \delta<\delta^*$ and for all $0<\kappa<\kappa^*$ the following statements hold true:
\begin{enumerate}[{\rm (1)}]
\setlength{\leftskip}{1pc}
 \item system \eqref{3er-sys} has exactly two strictly positive spatially homogeneous steady states, $(u_-(v_-), u_2^*(u_-(v_-),v_-), v_-)$ and $(u_-(v_+),u_2^*(u_-(v_+),v_+), v_+)$ with $v_-<v_+$.
 \item $(u_-(v_-), u_2^*(u_-(v_-),v_-), v_-)$ is an unstable steady state of the kinetic system for \eqref{3er-sys} and of the original system \eqref{3er-sys}.
 \item $(u_-(v_+), u_2^*(u_-(v_+) ,v_+),v_+)$ is a stable steady state of the kinetic system of \eqref{3er-sys} and an unstable steady state of \eqref{3er-sys}.
 \item $(0,0,0)$ is a stable steady state of the kinetic system for \eqref{3er-sys} and of the original system \eqref{3er-sys} .
 \item system \eqref{3er-sys} has infinitely many $(\varepsilon_0,A)$-stable steady states with jump discontinuity, with component $u_1$ being at $x\in \overline{\text{$I$}}$ of type $u_+$ or $u_0$.
\end{enumerate}
\end{lemma}

\subsection{Numerical results: Effects of initial conditions and diffusion coefficient.} $\,$\\
In the previous subsection, we showed that a spatially homogeneous steady state is destabilized and that infinitely many, discontinuous
steady states are stable, both due to introduction of diffusion. In this subsection, we show some numerical results showing that the arising pattern depends on both, 
the initial conditions and the diffusion coefficient. In figure \ref{fig:illu2}, numerical approximation of solutions to
problem \eqref{eq1}-\eqref{bc} are shown. We observe that for the same initial conditions and parameters of the kinetic system, 
the number of jump-type discontinuities varies according to the size of diffusion coefficient $D$. A trend towards a higher number for smaller
diffusion-coefficient can be observed and the shape of the arising pattern is different from the shape of the initial conditions.
For large diffusion coefficient, a `plateau' arises close to the spatial position on which the initial condition of $u$ assumes 
its maximal value. Here, initial conditions seem to determine the position of `plateaus' as can be observed in figure \ref{fig:illu3}.
Note that a local maximum is positioned at $x=0$ in figure \ref{fig:illu3}.
\\
\begin{minipage}[t]{36em}
\begin{tabular}{p{18em}p{18em}}
   \begingroup
  \makeatletter
  \providecommand\color[2][]{    \GenericError{(gnuplot) \space\space\space\@spaces}{      Package color not loaded in conjunction with
      terminal option `colourtext'    }{See the gnuplot documentation for explanation.    }{Either use 'blacktext' in gnuplot or load the package
      color.sty in LaTeX.}    \renewcommand\color[2][]{}  }  \providecommand\includegraphics[2][]{    \GenericError{(gnuplot) \space\space\space\@spaces}{      Package graphicx or graphics not loaded    }{See the gnuplot documentation for explanation.    }{The gnuplot epslatex terminal needs graphicx.sty or graphics.sty.}    \renewcommand\includegraphics[2][]{}  }  \providecommand\rotatebox[2]{#2}  \@ifundefined{ifGPcolor}{    \newif\ifGPcolor
    \GPcolorfalse
  }{}  \@ifundefined{ifGPblacktext}{    \newif\ifGPblacktext
    \GPblacktexttrue
  }{}    \let\gplgaddtomacro\g@addto@macro
    \gdef\gplbacktext{}  \gdef\gplfronttext{}  \makeatother
  \ifGPblacktext
        \def\colorrgb#1{}    \def\colorgray#1{}  \else
        \ifGPcolor
      \def\colorrgb#1{\color[rgb]{#1}}      \def\colorgray#1{\color[gray]{#1}}      \expandafter\def\csname LTw\endcsname{\color{white}}      \expandafter\def\csname LTb\endcsname{\color{black}}      \expandafter\def\csname LTa\endcsname{\color{black}}      \expandafter\def\csname LT0\endcsname{\color[rgb]{1,0,0}}      \expandafter\def\csname LT1\endcsname{\color[rgb]{0,1,0}}      \expandafter\def\csname LT2\endcsname{\color[rgb]{0,0,1}}      \expandafter\def\csname LT3\endcsname{\color[rgb]{1,0,1}}      \expandafter\def\csname LT4\endcsname{\color[rgb]{0,1,1}}      \expandafter\def\csname LT5\endcsname{\color[rgb]{1,1,0}}      \expandafter\def\csname LT6\endcsname{\color[rgb]{0,0,0}}      \expandafter\def\csname LT7\endcsname{\color[rgb]{1,0.3,0}}      \expandafter\def\csname LT8\endcsname{\color[rgb]{0.5,0.5,0.5}}    \else
            \def\colorrgb#1{\color{black}}      \def\colorgray#1{\color[gray]{#1}}      \expandafter\def\csname LTw\endcsname{\color{white}}      \expandafter\def\csname LTb\endcsname{\color{black}}      \expandafter\def\csname LTa\endcsname{\color{black}}      \expandafter\def\csname LT0\endcsname{\color{black}}      \expandafter\def\csname LT1\endcsname{\color{black}}      \expandafter\def\csname LT2\endcsname{\color{black}}      \expandafter\def\csname LT3\endcsname{\color{black}}      \expandafter\def\csname LT4\endcsname{\color{black}}      \expandafter\def\csname LT5\endcsname{\color{black}}      \expandafter\def\csname LT6\endcsname{\color{black}}      \expandafter\def\csname LT7\endcsname{\color{black}}      \expandafter\def\csname LT8\endcsname{\color{black}}    \fi
  \fi
  \setlength{\unitlength}{0.0425bp}  \begin{picture}(5760.00,4032.00)    \gplgaddtomacro\gplbacktext{    }    \gplgaddtomacro\gplfronttext{      \csname LTb\endcsname      \put(4087,3308){\makebox(0,0)[r]{\strut{}shape of $u_0$}}      \csname LTb\endcsname      \put(667,1269){\makebox(0,0)[r]{\strut{} 0}}      \put(2729,593){\makebox(0,0)[r]{\strut{} 10}}      \put(1334,870){\makebox(0,0){\strut{}$t$}}      \put(2973,563){\makebox(0,0){\strut{} 0}}      \put(5035,1199){\makebox(0,0){\strut{} 1}}      \put(4426,870){\makebox(0,0){\strut{}$x$}}      \put(692,2040){\makebox(0,0)[r]{\strut{} 0}}      \put(692,2732){\makebox(0,0)[r]{\strut{} 20}}      \put(891,819){\vector(3,-1){500}}    }    \gplbacktext
    \put(0,0){\includegraphics[scale=0.85]{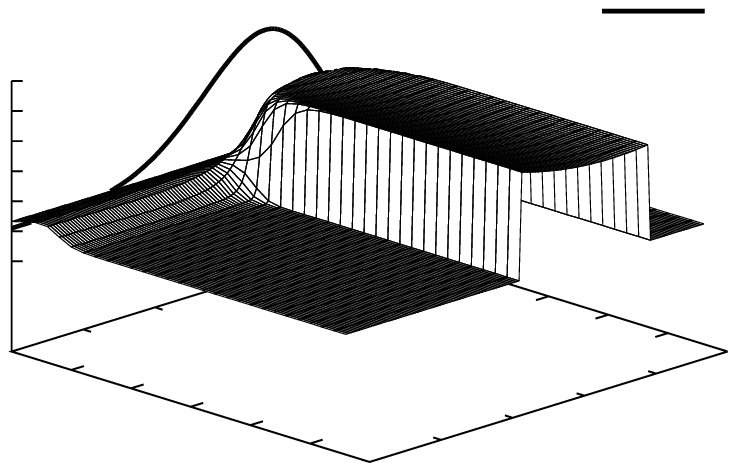}}    \gplfronttext
  \end{picture}\endgroup
 &
   \begingroup
  \makeatletter
  \providecommand\color[2][]{    \GenericError{(gnuplot) \space\space\space\@spaces}{      Package color not loaded in conjunction with
      terminal option `colourtext'    }{See the gnuplot documentation for explanation.    }{Either use 'blacktext' in gnuplot or load the package
      color.sty in LaTeX.}    \renewcommand\color[2][]{}  }  \providecommand\includegraphics[2][]{    \GenericError{(gnuplot) \space\space\space\@spaces}{      Package graphicx or graphics not loaded    }{See the gnuplot documentation for explanation.    }{The gnuplot epslatex terminal needs graphicx.sty or graphics.sty.}    \renewcommand\includegraphics[2][]{}  }  \providecommand\rotatebox[2]{#2}  \@ifundefined{ifGPcolor}{    \newif\ifGPcolor
    \GPcolorfalse
  }{}  \@ifundefined{ifGPblacktext}{    \newif\ifGPblacktext
    \GPblacktexttrue
  }{}    \let\gplgaddtomacro\g@addto@macro
    \gdef\gplbacktext{}  \gdef\gplfronttext{}  \makeatother
  \ifGPblacktext
        \def\colorrgb#1{}    \def\colorgray#1{}  \else
        \ifGPcolor
      \def\colorrgb#1{\color[rgb]{#1}}      \def\colorgray#1{\color[gray]{#1}}      \expandafter\def\csname LTw\endcsname{\color{white}}      \expandafter\def\csname LTb\endcsname{\color{black}}      \expandafter\def\csname LTa\endcsname{\color{black}}      \expandafter\def\csname LT0\endcsname{\color[rgb]{1,0,0}}      \expandafter\def\csname LT1\endcsname{\color[rgb]{0,1,0}}      \expandafter\def\csname LT2\endcsname{\color[rgb]{0,0,1}}      \expandafter\def\csname LT3\endcsname{\color[rgb]{1,0,1}}      \expandafter\def\csname LT4\endcsname{\color[rgb]{0,1,1}}      \expandafter\def\csname LT5\endcsname{\color[rgb]{1,1,0}}      \expandafter\def\csname LT6\endcsname{\color[rgb]{0,0,0}}      \expandafter\def\csname LT7\endcsname{\color[rgb]{1,0.3,0}}      \expandafter\def\csname LT8\endcsname{\color[rgb]{0.5,0.5,0.5}}    \else
            \def\colorrgb#1{\color{black}}      \def\colorgray#1{\color[gray]{#1}}      \expandafter\def\csname LTw\endcsname{\color{white}}      \expandafter\def\csname LTb\endcsname{\color{black}}      \expandafter\def\csname LTa\endcsname{\color{black}}      \expandafter\def\csname LT0\endcsname{\color{black}}      \expandafter\def\csname LT1\endcsname{\color{black}}      \expandafter\def\csname LT2\endcsname{\color{black}}      \expandafter\def\csname LT3\endcsname{\color{black}}      \expandafter\def\csname LT4\endcsname{\color{black}}      \expandafter\def\csname LT5\endcsname{\color{black}}      \expandafter\def\csname LT6\endcsname{\color{black}}      \expandafter\def\csname LT7\endcsname{\color{black}}      \expandafter\def\csname LT8\endcsname{\color{black}}    \fi
  \fi
  \setlength{\unitlength}{0.0425bp}  \begin{picture}(5760.00,4032.00)    \gplgaddtomacro\gplbacktext{    }    \gplgaddtomacro\gplfronttext{      \csname LTb\endcsname      \put(4087,3308){\makebox(0,0)[r]{\strut{}shape of $u_0$}}      \csname LTb\endcsname      \put(667,1269){\makebox(0,0)[r]{\strut{} 0}}      \put(2729,593){\makebox(0,0)[r]{\strut{} 10}}      \put(1334,870){\makebox(0,0){\strut{}$t$}}      \put(2973,563){\makebox(0,0){\strut{} 0}}      \put(5035,1199){\makebox(0,0){\strut{} 1}}      \put(4426,870){\makebox(0,0){\strut{}$x$}}      \put(692,2040){\makebox(0,0)[r]{\strut{} 0}}      \put(692,2732){\makebox(0,0)[r]{\strut{} 20}}      \put(891,819){\vector(3,-1){500}}    }    \gplbacktext
    \put(0,0){\includegraphics[scale=0.85]{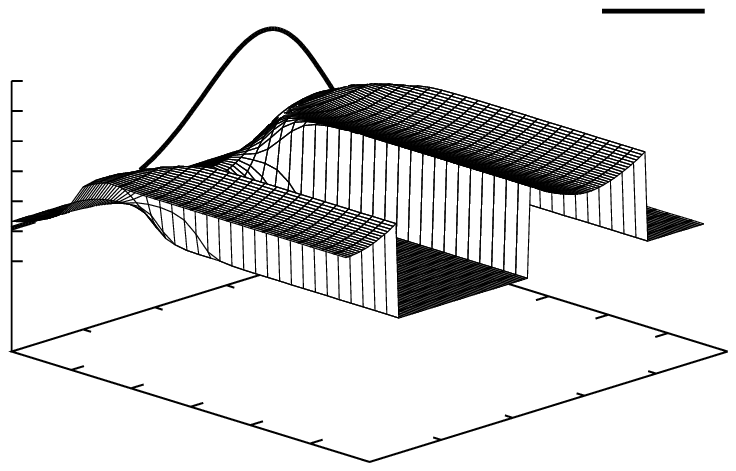}}    \gplfronttext
  \end{picture}\endgroup
 \\
   \begingroup
  \makeatletter
  \providecommand\color[2][]{    \GenericError{(gnuplot) \space\space\space\@spaces}{      Package color not loaded in conjunction with
      terminal option `colourtext'    }{See the gnuplot documentation for explanation.    }{Either use 'blacktext' in gnuplot or load the package
      color.sty in LaTeX.}    \renewcommand\color[2][]{}  }  \providecommand\includegraphics[2][]{    \GenericError{(gnuplot) \space\space\space\@spaces}{      Package graphicx or graphics not loaded    }{See the gnuplot documentation for explanation.    }{The gnuplot epslatex terminal needs graphicx.sty or graphics.sty.}    \renewcommand\includegraphics[2][]{}  }  \providecommand\rotatebox[2]{#2}  \@ifundefined{ifGPcolor}{    \newif\ifGPcolor
    \GPcolorfalse
  }{}  \@ifundefined{ifGPblacktext}{    \newif\ifGPblacktext
    \GPblacktexttrue
  }{}    \let\gplgaddtomacro\g@addto@macro
    \gdef\gplbacktext{}  \gdef\gplfronttext{}  \makeatother
  \ifGPblacktext
        \def\colorrgb#1{}    \def\colorgray#1{}  \else
        \ifGPcolor
      \def\colorrgb#1{\color[rgb]{#1}}      \def\colorgray#1{\color[gray]{#1}}      \expandafter\def\csname LTw\endcsname{\color{white}}      \expandafter\def\csname LTb\endcsname{\color{black}}      \expandafter\def\csname LTa\endcsname{\color{black}}      \expandafter\def\csname LT0\endcsname{\color[rgb]{1,0,0}}      \expandafter\def\csname LT1\endcsname{\color[rgb]{0,1,0}}      \expandafter\def\csname LT2\endcsname{\color[rgb]{0,0,1}}      \expandafter\def\csname LT3\endcsname{\color[rgb]{1,0,1}}      \expandafter\def\csname LT4\endcsname{\color[rgb]{0,1,1}}      \expandafter\def\csname LT5\endcsname{\color[rgb]{1,1,0}}      \expandafter\def\csname LT6\endcsname{\color[rgb]{0,0,0}}      \expandafter\def\csname LT7\endcsname{\color[rgb]{1,0.3,0}}      \expandafter\def\csname LT8\endcsname{\color[rgb]{0.5,0.5,0.5}}    \else
            \def\colorrgb#1{\color{black}}      \def\colorgray#1{\color[gray]{#1}}      \expandafter\def\csname LTw\endcsname{\color{white}}      \expandafter\def\csname LTb\endcsname{\color{black}}      \expandafter\def\csname LTa\endcsname{\color{black}}      \expandafter\def\csname LT0\endcsname{\color{black}}      \expandafter\def\csname LT1\endcsname{\color{black}}      \expandafter\def\csname LT2\endcsname{\color{black}}      \expandafter\def\csname LT3\endcsname{\color{black}}      \expandafter\def\csname LT4\endcsname{\color{black}}      \expandafter\def\csname LT5\endcsname{\color{black}}      \expandafter\def\csname LT6\endcsname{\color{black}}      \expandafter\def\csname LT7\endcsname{\color{black}}      \expandafter\def\csname LT8\endcsname{\color{black}}    \fi
  \fi
  \setlength{\unitlength}{0.0425bp}  \begin{picture}(5760.00,4032.00)    \gplgaddtomacro\gplbacktext{    }    \gplgaddtomacro\gplfronttext{      \csname LTb\endcsname      \put(4087,3308){\makebox(0,0)[r]{\strut{}shape of $u_0$}}      \csname LTb\endcsname      \put(667,1269){\makebox(0,0)[r]{\strut{} 0}}      \put(2729,593){\makebox(0,0)[r]{\strut{} 10}}      \put(1334,870){\makebox(0,0){\strut{}$t$}}      \put(2973,563){\makebox(0,0){\strut{} 0}}      \put(5035,1199){\makebox(0,0){\strut{} 1}}      \put(4426,870){\makebox(0,0){\strut{}$x$}}      \put(692,2040){\makebox(0,0)[r]{\strut{} 0}}      \put(692,2732){\makebox(0,0)[r]{\strut{} 20}}      \put(891,819){\vector(3,-1){500}}    }    \gplbacktext
    \put(0,0){\includegraphics[scale=0.85]{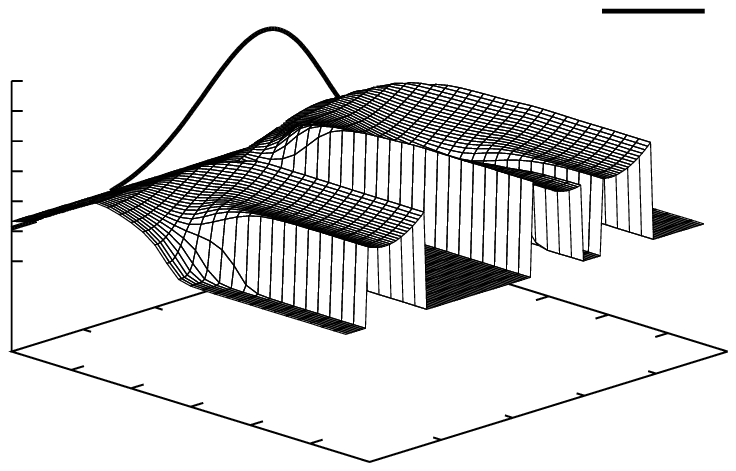}}    \gplfronttext
  \end{picture}\endgroup
 &
   \begingroup
  \makeatletter
  \providecommand\color[2][]{    \GenericError{(gnuplot) \space\space\space\@spaces}{      Package color not loaded in conjunction with
      terminal option `colourtext'    }{See the gnuplot documentation for explanation.    }{Either use 'blacktext' in gnuplot or load the package
      color.sty in LaTeX.}    \renewcommand\color[2][]{}  }  \providecommand\includegraphics[2][]{    \GenericError{(gnuplot) \space\space\space\@spaces}{      Package graphicx or graphics not loaded    }{See the gnuplot documentation for explanation.    }{The gnuplot epslatex terminal needs graphicx.sty or graphics.sty.}    \renewcommand\includegraphics[2][]{}  }  \providecommand\rotatebox[2]{#2}  \@ifundefined{ifGPcolor}{    \newif\ifGPcolor
    \GPcolorfalse
  }{}  \@ifundefined{ifGPblacktext}{    \newif\ifGPblacktext
    \GPblacktexttrue
  }{}    \let\gplgaddtomacro\g@addto@macro
    \gdef\gplbacktext{}  \gdef\gplfronttext{}  \makeatother
  \ifGPblacktext
        \def\colorrgb#1{}    \def\colorgray#1{}  \else
        \ifGPcolor
      \def\colorrgb#1{\color[rgb]{#1}}      \def\colorgray#1{\color[gray]{#1}}      \expandafter\def\csname LTw\endcsname{\color{white}}      \expandafter\def\csname LTb\endcsname{\color{black}}      \expandafter\def\csname LTa\endcsname{\color{black}}      \expandafter\def\csname LT0\endcsname{\color[rgb]{1,0,0}}      \expandafter\def\csname LT1\endcsname{\color[rgb]{0,1,0}}      \expandafter\def\csname LT2\endcsname{\color[rgb]{0,0,1}}      \expandafter\def\csname LT3\endcsname{\color[rgb]{1,0,1}}      \expandafter\def\csname LT4\endcsname{\color[rgb]{0,1,1}}      \expandafter\def\csname LT5\endcsname{\color[rgb]{1,1,0}}      \expandafter\def\csname LT6\endcsname{\color[rgb]{0,0,0}}      \expandafter\def\csname LT7\endcsname{\color[rgb]{1,0.3,0}}      \expandafter\def\csname LT8\endcsname{\color[rgb]{0.5,0.5,0.5}}    \else
            \def\colorrgb#1{\color{black}}      \def\colorgray#1{\color[gray]{#1}}      \expandafter\def\csname LTw\endcsname{\color{white}}      \expandafter\def\csname LTb\endcsname{\color{black}}      \expandafter\def\csname LTa\endcsname{\color{black}}      \expandafter\def\csname LT0\endcsname{\color{black}}      \expandafter\def\csname LT1\endcsname{\color{black}}      \expandafter\def\csname LT2\endcsname{\color{black}}      \expandafter\def\csname LT3\endcsname{\color{black}}      \expandafter\def\csname LT4\endcsname{\color{black}}      \expandafter\def\csname LT5\endcsname{\color{black}}      \expandafter\def\csname LT6\endcsname{\color{black}}      \expandafter\def\csname LT7\endcsname{\color{black}}      \expandafter\def\csname LT8\endcsname{\color{black}}    \fi
  \fi
  \setlength{\unitlength}{0.0425bp}  \begin{picture}(5760.00,4032.00)    \gplgaddtomacro\gplbacktext{    }    \gplgaddtomacro\gplfronttext{      \csname LTb\endcsname      \put(4087,3308){\makebox(0,0)[r]{\strut{}shape of $u_0$}}      \csname LTb\endcsname      \put(667,1269){\makebox(0,0)[r]{\strut{} 0}}      \put(2729,593){\makebox(0,0)[r]{\strut{} 10}}      \put(1334,870){\makebox(0,0){\strut{}$t$}}      \put(2973,563){\makebox(0,0){\strut{} 0}}      \put(5035,1199){\makebox(0,0){\strut{} 1}}      \put(4426,870){\makebox(0,0){\strut{}$x$}}      \put(692,2040){\makebox(0,0)[r]{\strut{} 0}}      \put(692,2732){\makebox(0,0)[r]{\strut{} 20}}      \put(891,819){\vector(3,-1){500}}    }    \gplbacktext
    \put(0,0){\includegraphics[scale=0.85]{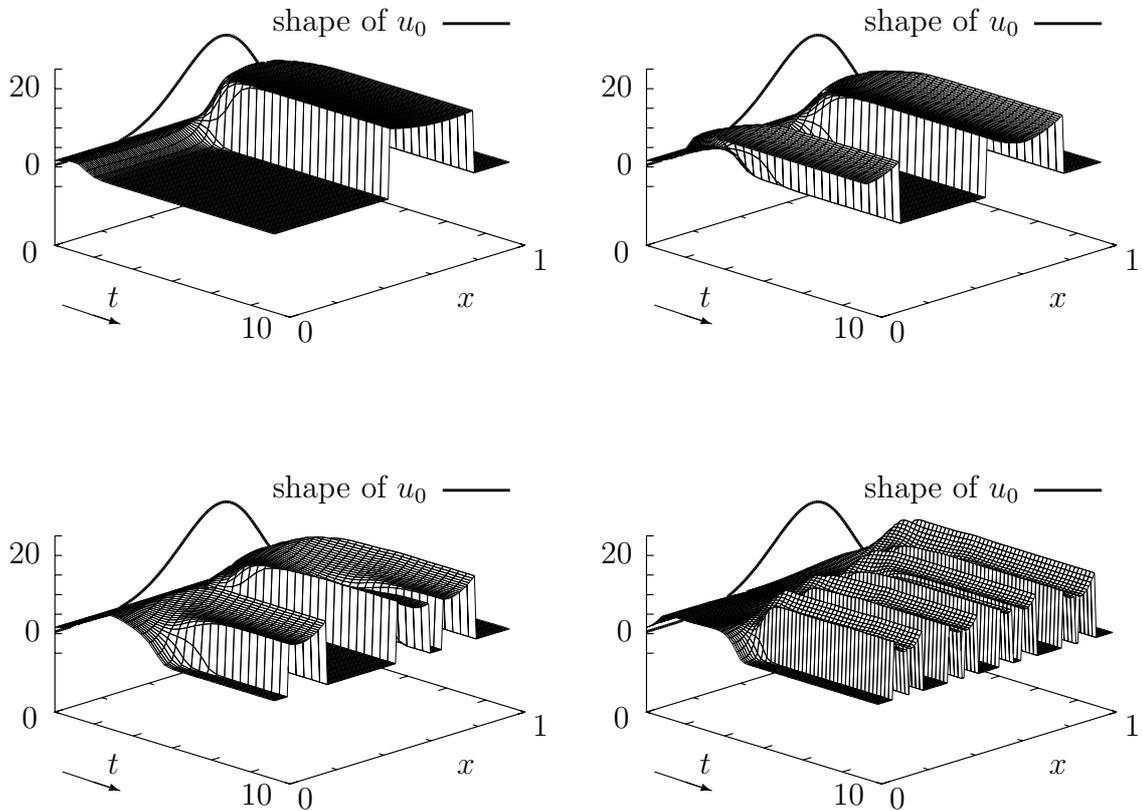}}    \gplfronttext
  \end{picture}\endgroup
 \end{tabular}
\captionof{figure}{Numerically obtained solution component $u$ of model \eqref{eq1}-\eqref{bc} for parameters $m_1 =1.44, m_2 = 2,\mu_3 \approx 4.1, k=0.01$ with varying diffusion coefficient: upper left: $D=5$, upper right: $D=1$, lower left: $D=0.5$, lower right: $D=0.1$.
We observe emergence of more jump-type discontinuities for smaller diffusion coefficient. Initial conditions are $u_0(x)=1.725-0.1\cos(2 \pi x^2), v_0(x)=2.48615$.}
\label{fig:illu2}
\end{minipage}\\
\begin{minipage}[t]{36em}
\center{\begingroup
  \makeatletter
  \providecommand\color[2][]{    \GenericError{(gnuplot) \space\space\space\@spaces}{      Package color not loaded in conjunction with
      terminal option `colourtext'    }{See the gnuplot documentation for explanation.    }{Either use 'blacktext' in gnuplot or load the package
      color.sty in LaTeX.}    \renewcommand\color[2][]{}  }  \providecommand\includegraphics[2][]{    \GenericError{(gnuplot) \space\space\space\@spaces}{      Package graphicx or graphics not loaded    }{See the gnuplot documentation for explanation.    }{The gnuplot epslatex terminal needs graphicx.sty or graphics.sty.}    \renewcommand\includegraphics[2][]{}  }  \providecommand\rotatebox[2]{#2}  \@ifundefined{ifGPcolor}{    \newif\ifGPcolor
    \GPcolorfalse
  }{}  \@ifundefined{ifGPblacktext}{    \newif\ifGPblacktext
    \GPblacktexttrue
  }{}    \let\gplgaddtomacro\g@addto@macro
    \gdef\gplbacktext{}  \gdef\gplfronttext{}  \makeatother
  \ifGPblacktext
        \def\colorrgb#1{}    \def\colorgray#1{}  \else
        \ifGPcolor
      \def\colorrgb#1{\color[rgb]{#1}}      \def\colorgray#1{\color[gray]{#1}}      \expandafter\def\csname LTw\endcsname{\color{white}}      \expandafter\def\csname LTb\endcsname{\color{black}}      \expandafter\def\csname LTa\endcsname{\color{black}}      \expandafter\def\csname LT0\endcsname{\color[rgb]{1,0,0}}      \expandafter\def\csname LT1\endcsname{\color[rgb]{0,1,0}}      \expandafter\def\csname LT2\endcsname{\color[rgb]{0,0,1}}      \expandafter\def\csname LT3\endcsname{\color[rgb]{1,0,1}}      \expandafter\def\csname LT4\endcsname{\color[rgb]{0,1,1}}      \expandafter\def\csname LT5\endcsname{\color[rgb]{1,1,0}}      \expandafter\def\csname LT6\endcsname{\color[rgb]{0,0,0}}      \expandafter\def\csname LT7\endcsname{\color[rgb]{1,0.3,0}}      \expandafter\def\csname LT8\endcsname{\color[rgb]{0.5,0.5,0.5}}    \else
            \def\colorrgb#1{\color{black}}      \def\colorgray#1{\color[gray]{#1}}      \expandafter\def\csname LTw\endcsname{\color{white}}      \expandafter\def\csname LTb\endcsname{\color{black}}      \expandafter\def\csname LTa\endcsname{\color{black}}      \expandafter\def\csname LT0\endcsname{\color{black}}      \expandafter\def\csname LT1\endcsname{\color{black}}      \expandafter\def\csname LT2\endcsname{\color{black}}      \expandafter\def\csname LT3\endcsname{\color{black}}      \expandafter\def\csname LT4\endcsname{\color{black}}      \expandafter\def\csname LT5\endcsname{\color{black}}      \expandafter\def\csname LT6\endcsname{\color{black}}      \expandafter\def\csname LT7\endcsname{\color{black}}      \expandafter\def\csname LT8\endcsname{\color{black}}    \fi
  \fi
  \setlength{\unitlength}{0.0425bp}  \begin{picture}(5760.00,4032.00)    \gplgaddtomacro\gplbacktext{    }    \gplgaddtomacro\gplfronttext{      \csname LTb\endcsname      \put(4187,3308){\makebox(0,0)[r]{\strut{}shape of $u_0$}}      \csname LTb\endcsname      \put(667,1269){\makebox(0,0)[r]{\strut{} 0}}      \put(2729,593){\makebox(0,0)[r]{\strut{} 10}}      \put(1334,870){\makebox(0,0){\strut{}$t$}}      \put(2973,563){\makebox(0,0){\strut{} 0}}      \put(5035,1199){\makebox(0,0){\strut{} 1}}      \put(4426,870){\makebox(0,0){\strut{}$x$}}      \put(891,819){\vector(3,-1){500}}      \put(692,1961){\makebox(0,0)[r]{\strut{} 0}}      \put(692,2905){\makebox(0,0)[r]{\strut{} 20}}    }    \gplbacktext
    \put(0,0){\includegraphics[scale=0.85]{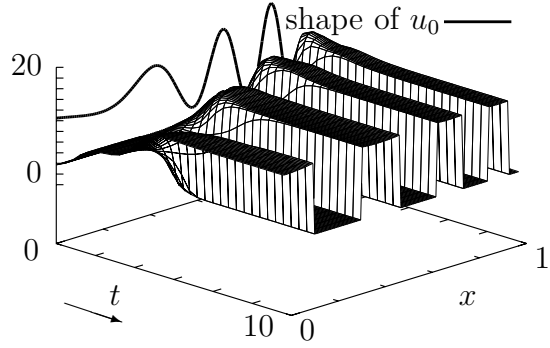}}    \gplfronttext
  \end{picture}\endgroup
 }
\captionof{figure}{Numerically obtained solution component $u$ of model \eqref{eq1}-\eqref{bc} for parameters $m_1 =1.44, m_2 = 2,\mu_3 \approx 4.1, k=0.01, D=5$.
We observe emergence of jump-type discontinuities around local maxima of the initial conditions. Initial conditions are $u_0(x)=1.725-0.1x^4\cos(8 \pi x^2), v_0(x)=2.48615$.}
\label{fig:illu3}
\end{minipage}\\

\newcommand{\pd}[2]{\frac{\partial #1}{\partial #2}}
\newcommand{\dpd}[2]{\frac{\partial^2 #1}{\partial #2^2}}
\newcommand{\od}[2]{\frac{d #1}{d #2}}
\newcommand{\dod}[2]{\frac{d^2 #1}{d #2^2}}

\section{Proof of the stability theorems} \label{sec:mainproof}
\subsection{Spectral structure of the linearized operator} $\,$\\
To handle Theorem \ref{thm:stabC} we introduce the following notation for simplicity:
\begin{align*}
& \eb_{11}(x) = \pd{f}{u}(\tilde u(x), \tilde v(x)), \quad \eb_{12}(x) = \pd{f}{v}(\tilde u(x), \tilde v(x)), \\
& \eb_{21}(x) = \pd{g}{u}(\tilde u(x), \tilde v(x)), \quad \eb_{22}(x) = \pd{g}{v}(\tilde u(x), \tilde v(x)).
\end{align*}
\begin{lemma}\label{L4.1} {Under the assumptions of Theorem \ref{thm:stabC}, consider the operator
\begin{equation}\label{A4.1}
\mathcal{L}_1 : \begin{pmatrix}\phi \\ \psi \end{pmatrix} \mapsto \begin{pmatrix} \eb_{11}(x) &  \eb_{12}(x)  \\ \eb_{21}(x) & \eb_{22}(x) + D \dod{}{x} \\\end{pmatrix} \begin{pmatrix}\phi \\ \psi \\ \end{pmatrix}
\end{equation}
on $(L^2(I))^2 $ with domain $L^2(I) \times (H^2(I) \cap \{ {\psi} \mid {\psi}^\prime=0 \ \hbox{at}\ x\in \partial I\} )$. Then $- \mathcal{L}_1$ is a sectorial operator on $(L^2(I))^2$, i.e., there exist positive constants $M_1$, $\kappa_0$ and $\omega \in (0, \pi/2)$ such that
\begin{equation}
\left\Vert \left(\lambda + \mathcal{L}_1 \right)^{-1} \right\Vert_{B((L^2(I))^2)} \leq \frac{M_1}{| \lambda - \kappa_0|+\kappa_0} \quad \hbox{for all } \ \lambda \in \{\lambda \in \mathbb{C} \mid |{\rm arg}\, (\lambda -\kappa_0) | > \omega \},
\label{A4.2}
\end{equation}
where ${B}((L^2(I))^2) $ denotes the Banach space of all bounded linear operators on $(L^2(I))^2$ equipped with the operator norm.
Hence $\mathcal{L}_1$ generates the analytic semigroup $e^{t \mathcal{L}_1}$ on $(L^2(I))^2$ and satisfies the estimate
\begin{equation}
\left\Vert e^{t \mathcal{L}_1} \begin{pmatrix}\phi \\ \psi \\\end{pmatrix}  \right\Vert_{(L^2(I))^2} \leq M e^{-\kappa_1 t} \left(\Vert \phi \Vert_{L^2(I)} + \Vert \psi \Vert_{L^2(I)} \right) \qquad \hbox{for all} \ t \geq 0.
\label{A4.3}
\end{equation}
Here $M$ and $\kappa_1$ are positive constants independent of $(\phi, \psi)$.
} 
\end{lemma}

\begin{lemma} \label{L4.2}
{Under the assumptions of Theorem \ref{thm:stabC-3c-gen}, consider the operator
\begin{equation}
\mathcal{ L}_2 : \begin{pmatrix} \phi_1 \\ \phi_2 \\ \psi \\ \end{pmatrix} \mapsto \begin{pmatrix} a_{11}(x) &  a_{12}(x)  & a_{13}(x) \\ a_{21}(x) & a_{22}(x) & a_{23}(x) \\
a_{31}(x) & a_{32}(x) & a_{33}(x) + D \dod{}{x} \\ \end{pmatrix} \begin{pmatrix} \phi_1 \\ \phi_2 \\ \psi \\ \end{pmatrix}
\label{A4.4}
\end{equation}
defined on $(L^2(I))^3$ with domain $(L^2(I))^2 \times (H^2(I) \cap \{\psi \mid \psi^\prime=0 \ \hbox{at}\ x\in \partial I\} )$. Then $- \mathcal{L}_2 $ is a sectorial operator on $(L^2(I))^3$, i.e., there exist positive constants $M_2$, $\kappa_0$ and $\omega \in (0, \pi/2)$ such that
\begin{equation}
\left\Vert \left(\lambda + \mathcal{L}_2 \right)^{-1} \right\Vert_{B((L^{2}(I))}^3) \leq \frac{M_2}{| \lambda - \kappa_0|+\kappa_0} \quad \hbox{for all } \ \lambda \in \{\lambda \in \mathbb{C} \mid |{\rm arg}\, (\lambda - \kappa_0) | > \omega \},
\label{A4.5}
\end{equation}
where ${B}((L^2(I))^3) $ denotes the Banach space of all bounded linear operators on $(L^2(I))^3$ equipped with the operator norm.
Hence $\mathcal{L}_2$ generates the analytic semigroup $e^{t \mathcal{L}_2}$ on $(L^2(I))^3 $ and satisfies the estimate
\begin{equation}
\left\Vert e^{t \mathcal{L}_2} \begin{pmatrix} \phi_1 \\ \phi_2 \\ \psi \\ \end{pmatrix}  \right\Vert_{(L^2(I))^3} \leq M e^{-\kappa_1 t} \left(\Vert \phi_1 \Vert_{L^2(I)} + \Vert \phi_2 \Vert_{L^2(I)} + \Vert \psi \Vert_{L^2(I)} \right) \quad \hbox{for all} \ t \geq 0.
\label{A4.6}
\end{equation}
Here $M$ and $\kappa_1$ are positive constants independent of $(\phi_1, \phi_2, \psi)$.
} 
\end{lemma}

Since Lemma \ref{L4.1} is obtained as a corollary to Lemma \ref{L4.2}, we prove only Lemma \ref{L4.2}. (Indeed,  choosing $\alpha>0$ sufficiently large, we see that the matrix
\begin{equation*}
\begin{pmatrix} -\alpha & 0 & 0 \\ 0 & \eb_{11}(x) & \eb_{12}(x) \\ 0 & \eb_{21}(x) & \eb_{22}(x) \\ \end{pmatrix}
\end{equation*}
satisfies conditions \eqref{A2.11}--\eqref{A2.14} of Theorem \ref{thm:stabC-3c-gen}.)

{\it Proof of Lemma \ref{L4.2}\/.}  Instead of $-\mathcal{L}_2$ we study $\mathcal {L}_2$; the estimate \eqref{A4.5} is obtained by replacing $\lambda$ with $-\lambda$ in \eqref{A4.888} below. Hence, we consider the nonhomogeneous equation
\begin{equation}
\mathcal{L}_2 \begin{pmatrix} \phi_1 \\ \phi_2 \\ \psi \\\end{pmatrix} = \lambda \begin{pmatrix}\phi_1 \\ \phi_2 \\ \psi \\ \end{pmatrix} + \begin{pmatrix}r_1 \\ r_2 \\ s \\ \end{pmatrix}
\label{A4.7}
\end{equation}
for $r_1,\, r_2, \, s \in L^2(I)$ under homogeneous Neumann boundary conditions on $\psi$. 

{\it Step 1.}\quad By assumptions \eqref{A2.11}--\eqref{A2.13} we see that ${\rm tr}\, A_{33}(x) < 0 < \det A_{33}(x)$, and hence both eigenvalues $\mu_1(x)$ and $\mu_2(x)$ of $A_{33}(x)$ lie in the left-half plane. In what follows we number $\mu_1$, $\mu_2$ so that ${\rm Re}\,\mu_1(x) \leq {\rm Re}\, \mu_2(x)$. More precisely, there exists a positive number $\mu^*$ such that the rectangle $\mathcal{R} = \{ \lambda \in \mathbb{C} \mid -\mu^* \leq {\rm Re}\,\lambda \leq - 3\kappa,\  |{\rm Im}\, \lambda | \leq \mu^*\} $ contains $\mu_1(x)$ and $\mu_2(x)$ for all $x\in \overline{I}$. Indeed, $\mu_j(x)$ are bounded because $a_{ij}(x)$ are bounded. It is elementary to check that $
{\rm Re}\, \mu_j(x) \leq \max\{{\rm tr}\, A_{33}(x)/2, \det A_{33}(x) / {\rm tr}\, A_{33}(x) \} $. By assumption \eqref{A2.14}, the right-hand side of this inequality does not exceed $-3\kappa$.
Therefore, in particular,
\begin{equation*}
\sigma (A_{33}(x)) \subset \mathcal{R} \subset \{ \lambda \in \mathbb{C} \mid {\rm Re}\, \lambda \leq -3 \kappa\}.
\end{equation*}
Hence, whenever ${\rm Re}\, \lambda > - 3\kappa$, we can solve
\begin{equation}
(A_{33}(x) - \lambda) \begin{pmatrix} \phi_1 \\ \phi_2\\ \end{pmatrix} = - \begin{pmatrix} a_{13}(x) \\ a_{23}(x) \\\end{pmatrix} \psi + \begin{pmatrix} r_1 \\ r_2 \\ \end{pmatrix}
\label{A4.10}
\end{equation}
for $(\phi_1, \phi_2)$:
\begin{equation}\label{A4.11}
\begin{aligned}
\phi_1 & = \bigg\{ \det \begin{pmatrix} r_1 & a_{12}(x) \\ r_2 & a_{22}(x) - \lambda \\ \end{pmatrix} - \det \begin{pmatrix} a_{13}(x) & a_{12}(x) \\ a_{23}(x) & a_{22}(x) - \lambda \\ \end{pmatrix} \psi \bigg\} \bigg/ \det (A_{33}(x) - \lambda), \\
\phi_2 & = \bigg\{ \det \begin{pmatrix} a_{11}(x) - \lambda & r_1 \\ a_{21}(x) & r_2 \\ \end{pmatrix} - \det \begin{pmatrix} a_{11}(x) - \lambda & a_{13}(x) \\ a_{21}(x) & a_{23}(x) \\ \end{pmatrix} \psi \bigg\} \bigg/ \det (A_{33}(x) - \lambda).
\end{aligned}
\end{equation}
Inserting this result into the third equation yields a scalar inhomogeneous Sturm-Liouville problem subject to homogeneous Neumann boundary conditions:
\begin{equation}
D \dod{\psi}{x} - q(x,\lambda) \psi = p_1(x,\lambda) r_1 + p_2(x,\lambda) r_2 + s,
\label{A4.12}
\end{equation}
where
\begin{equation}
\begin{aligned}
& p_1(x,\lambda) = - \frac{\det A_{13}(x) + \lambda a_{31}(x)}{\det (A_{33}(x)-\lambda)}, \quad
p_2(x,\lambda) = - \frac{\det A_{23}(x) - \lambda a_{32}(x)}{\det (A_{33}(x) - \lambda)}, \\
& q(x,\lambda) = - \frac{\det (A(x)-\lambda)}{\det (A_{33}(x)-\lambda)}.
\end{aligned}
\label{A4.13}
\end{equation}

{\it Step 2.} \quad
Let us define a sesquilinear form $\mathcal{B}(\psi_1, \psi_2) $ on $H^1(I)$ by
\begin{equation}
\mathcal{B}(\psi_1, \psi_2) = D \int \limits_I \od{\psi_1}{x} \overline{\od{\psi_2}{x}}\,dx + \int \limits_I q(x,\lambda) \psi_1 \, \overline{\psi_2}\,dx.
\label{A4.14}
\end{equation}

\begin{sublemma}\label{SL4.3}{There exist positive constants $\kappa_0$ and $\mathit{\Gamma}_1$ such that if $\lambda \in \mathbb{C}$ satisfies ${\rm Re}\,\lambda \geq -2\kappa_0$ or $|{\rm Im}\,\lambda| \geq \mathit{\Gamma}_1$, then the sesquilinear form $\mathcal{B}$ is bounded and  coercive. Namely, there exist positive constants $M_0$, $\gamma_0$ such that
\begin{align}
& |\mathcal{B}(\psi_1, \psi_2)| \leq M_0(1+|\lambda|) \Vert \psi_1 \Vert_{H^1(I)} \Vert \psi_2 \Vert_{H^1(I)} \quad \hbox{for all}\ \psi_1,\,\psi_2 \in H^1(I), & \label{A4.200}\\
& {\rm Re}\, \mathcal{B}(\psi,\psi) \geq \gamma_0 \Vert \psi \Vert_{H^1(I)}^2 \quad \hbox{for all}\ \psi \in H^1(I). & \label{A4.201}
\end{align}
Here, $\gamma_0$ depends on $\lambda$.
}
\end{sublemma}

{\it Proof of Sublemma \ref{SL4.3}.}\quad First we prove the boundedness of $\mathcal{B}$. 
It is convenient to introduce the following notation:
\begin{equation}
\left\{
\begin{aligned}
& D(x) = \det A(x), \quad T(x) = {\rm tr}\, A(x), \quad S(x) = \sum_{j=1}^3 \det A_{jj}(x), \\
& D_3(x)=\det A_{33}(x), \quad T_3(x) = {\rm tr}\, A_{33}(x), \\
& P(\lambda) = \det (A(x)-\lambda), \quad Q(\lambda)=\det(A_{33}(x)-\lambda).
\end{aligned}
\right.
\label{A4.16}
\end{equation}
We note that $P(\lambda)  = D(x)-S(x)\lambda+T(x)\lambda^2-\lambda^3 $ and $ Q(\lambda)=D_3(x) - T_3(x) \lambda + \lambda^2$, so that
\begin{equation}
|q(x,\lambda)| = |P(\lambda)/Q(\lambda)| = |\lambda| + O(1) \qquad \hbox{as}\quad |\lambda| \to +\infty
\label{A4.17}
\end{equation}
uniformly in $x \in \overline{I}$.  Moreover, zeros of $Q(\lambda)$ are both confined in the rectangle $\mathcal{R}$ defined in Step 1 for all $x \in \overline{I}$. Therefore, there exist positive constants $ \mathit{\Gamma}_0$ and $m_0$ such that
$$
|q(x,\lambda)| \leq m_0(|\lambda|+1) 
$$
for all $\lambda$ satisfying ${\rm Re}\, \lambda + 2\kappa \geq 0$ or $|{\rm Im}\, \lambda| \geq \mathit{\Gamma}_0$, which implies the desired inequality:
\begin{equation*}
|\mathcal{B}(\psi_1, \psi_2)| \leq D \Vert \psi_1^\prime \Vert_{L^2(I)} \Vert \psi_2^\prime \Vert_{L^2(I)} + m_0 (|\lambda|+1) \Vert \psi_1 \Vert_{L^2(I)} \Vert \psi_2 \Vert_{L^2(I)}.
\end{equation*}

Second, we prove the coerciveness of $\mathcal{B}$. 
Let us start with the case where $|\lambda|$ is sufficiently large. From
\begin{align}
P(\lambda)Q(\overline{\lambda}) + P(\overline{\lambda})Q(\lambda)
& = 2D(x)D_3(x) -( D_3(x)S(x)+D(x)T_3(x))(\lambda+\overline{\lambda}) & \label{A4.20}\\
& \quad +(D_3(x) T(x)+ D(x)) (\lambda^2+\overline{\lambda}^2) +2S(x)T_3(x) |\lambda|^2 \nonumber \\
& \quad -D_3(x) (\lambda^3+\overline{\lambda}^3) -(T(x)T_3(x)+S(x))|\lambda|^2(\lambda+\overline{\lambda}) \nonumber \\
& \qquad + T_3(x)|\lambda|^2(\lambda^2 + \overline{\lambda}^2) + 2T(x)|\lambda|^4 - |\lambda|^4(\lambda+\overline{\lambda}), \nonumber
\end{align}
it follows, due to $\lambda^2+\overline{\lambda}^2=(\lambda+\overline{\lambda})^2-2|\lambda|^2$ and $T(x)-T_{3}(x)=a_{33}(x)$, that
\begin{equation}
-{\rm Re}\,P(\lambda)Q(\overline{\lambda}) = |\lambda|^4 \left({\rm Re}\,\lambda - a_{33}(x) - 2T_3(x) ({\rm Re}\,\lambda)^2 / |\lambda|^2 + O(|\lambda|^{-1}) \right)
\label{A4.21}
\end{equation}
as $|\lambda| \to +\infty$, uniformly in $x \in \overline{I}$. On the other hand, we have
\begin{align}
|Q(\lambda)|^2 
& = D_{3}(x)^2 - 2 D_{3}(x) T_{3}(x)  + (T_{3}(x)^2 - 2 D_{3}(x)) |\lambda|^2  & \label{A4.22} \\
& \quad + 4D_{3}(x) ({\rm Re}\,\lambda)^2- 2 T_{3}(x) |\lambda|^2{\rm Re}\,\lambda +|\lambda|^4. \nonumber
\end{align}
Therefore, we see that
\begin{equation}
{\rm Re}\, q(x,\lambda) = {\rm Re}\, \lambda - a_{33}(x) - 2T_3(x) ({\rm Re}\,\lambda)^2 / |\lambda|^2 + O(|\lambda|^{-1}) \qquad \hbox{as}\quad |\lambda| \to +\infty,
\label{A4.23}
\end{equation}
uniformly in $x \in \overline{I}$.  Since ${\rm Re}\,\lambda - a_{33}(x) \geq {\rm Re}\, \lambda +3\kappa $ and $T_3(x) < 0$ by assumption \eqref{A2.14}, we conclude that there exists a positive constant $\mathit{\Gamma}_1 \geq \mathit{\Gamma}_0$ such that
\begin{equation}
{\rm Re}\, \mathcal{B}(\psi,\psi) \geq D \Vert \psi^\prime \Vert_{L^2(I)}^2 + ({\rm Re}\, \lambda + 3\kappa) \Vert \psi \Vert_{L^2(I)}^2 \geq \min\{D,\kappa/2 \} \Vert \psi \Vert_{H^1(I)}^2 
\label{A4.24}
\end{equation}
whenever ${\rm Re}\,\lambda \geq -2\kappa$ and $ |{\rm Im}\, \lambda|\geq \mathit{\Gamma}_1$.

Next we turn to the case where $|{\rm Im}\, \lambda|\leq \mathit{\Gamma}_1$. By direct computations we have
\begin{align*}
{\rm Re}\, \det (A(x)-\lambda) & = \det (A(x)- {\rm Re}\, \lambda) +({\rm Im}\,\lambda)^2(3 {\rm Re}\,\lambda - T(x)), \\
{\rm Im}\, \det (A(x)- \lambda) & = ({\rm Im}\,\lambda) \big( ({\rm Im}\,\lambda)^2 +(2T(x) - 3{\rm Re}\,\lambda) {\rm Re}\, \lambda - S(x) \big), \\
{\rm Re}\, \det (A_{33}(x)-\lambda) & = \det (A_{33}(x)-{\rm Re}\,\lambda) - ({\rm Im}\,\lambda)^2, \\
{\rm Im}\, \det (A_{33}(x)-\lambda) & = ({\rm Im}\,\lambda) (2 {\rm Re}\, \lambda - T_3(x)).
\end{align*}
Using these expressions, we obtain  
\begin{equation}
{\rm Re}\, q(x,\lambda) = \frac{1}{|\det (A_{33}(x)- \lambda)|^2} \sum_{j=0}^5 q_k(\eta) \xi^k \quad \hbox{for}\ \lambda=\xi+i\eta,
\label{A4.25}
\end{equation}
where we have defined as follows:
\begin{align*}
q_5(\eta) & = 1, \\
q_4(\eta) & = - T(x) - T_3(x), \\
q_3(\eta) & = D_3(x) +S(x) +T(x)T_3(x) +2\eta^2, \\
q_2(\eta) & = - \left ( D(x) + D_3(x) T(x) +S(x)T_3(x) + 2 T(x) \eta^2 \right) , \\
q_1(\eta) & = D_3(x) S(x) + D(x)T_3(x) + (T(x)T_3(x) + S(x) - 3 D_3(x))\eta^2 +\eta^4, \\
q_0(\eta) & = -D(x)D_3(x) + (D(x) + D_3(x)T(x) - S(x) T_3(x) ) \eta^2 - (T(x)-T_3(x)) \eta^4.
\end{align*}
Thanks to assumptions \eqref{A2.11}--\eqref{A2.14} of Theorem \ref{thm:stabC-3c-gen}, we see that
\begin{equation}
q_5(\eta) = 1, \quad q_4(\eta) > 0, \quad q_3(\eta) >0, \quad
 q_2(\eta) >0, \quad q_1(\eta) >0, \quad q_0(\eta) >0
\label{A4.26}
\end{equation}
for all $x \in \overline{I}$ and ${\rm Re}\,\lambda + 2\kappa \geq 0$. 
Observe that
\begin{align}
& \xi^5+q_4(\eta)\xi^4  =\xi^4(\xi-(T+T_3)); & \label{A4.27} \\
& q_3(\eta)\xi^3 + q_2(\eta)\xi^2  \geq \xi^2((D_3+TT_3 + S+2\eta^2)\xi  - (D+D_3T+ST_3)); & \label{A4.28} \\
& q_1(\eta)\xi + q_0(\eta) \geq \left\{ D_3 S + DT_3 +(T T_3 +S-3D_3) \eta^2 + \eta^4\right\} \xi - DD_3. & \label{A4.29}
\end{align}

Let $\kappa_0 $ be a positive constant defined by
\begin{align*}
-3\kappa_0 & = \max  \left\{ \sup_{x \in I} (T(x)+T_3(x)), \ \sup_{x \in I} \frac{D(x)+D_3(x)T(x)+S(x)T_3(x)}{D_3(x)+T(x)T_3(x)+S(x)+2\mathit{ \Gamma}_{1}^2}, \right. \\
& \qquad\qquad \left. \sup_{x \in I} \frac{D(x)D_3(x)}{D_3(x) S(x) + D(x)T_3(x)  + (T(x)T_3(x)+S(x)-3D_3(x))\mathit{\Gamma}_{1}^2 + \mathit{\Gamma}_{1}^4} \right\}.
\end{align*}
We may assume that $\kappa_0  \leq \kappa$.
If $\xi \geq -3\kappa_0$ and $| {\rm Im}\, \lambda | \leq \mathit{ \Gamma}_{1}$, then
\begin{align*}
& \xi -(T+T_3) \geq \xi + 3\kappa_0; \\
& (D_3+TT_3+S+2\eta^2)\xi - (D+D_3T+ST_3) \\
& \quad \geq  (D_3+TT_3+S+2\eta^2) \xi + 3 \kappa_0 (D_3+TT_3+S+2\mathit{ \Gamma}_{1}^2)\geq (D_3+TT_3+S)(\xi+3\kappa_0);\\
& \left\{(D_3S+DT_3) +(TT_3+S-3D_3)\eta^2 + \eta^4\right\}\xi -DD_3 \\
& \quad \geq \left\{(D_3S+DT_3) +(TT_3+S-3D_3)\eta^2 + \eta^4\right\}\xi \\
& \qquad + 3 \kappa_0 \left\{(D_3S+DT_3) +(TT_3+S-3D_3) \mathit{ \Gamma}_{1}^2 + \mathit{ \Gamma}_{1}^4\right\} \geq (D_3S+DT_3) (\xi+3\kappa_0).
\end{align*}
Therefore,
\begin{equation}
\begin{aligned} 
& \xi^5+q_4(\eta)\xi^4+q_3(\eta)\xi^3+q_2(\eta)\xi^2+q_1(\eta)\xi+q_0(\eta) \\
& \quad \geq (\xi+3\kappa_0) (\xi^4+(D_3(x)+T(x)T_3(x)+S(x))\xi^2+D_3(x)S(x)+D(x)T_3(x)).
\end{aligned}
\label{A4.30}
\end{equation}

Here we note the following elementary lemma, the proof of which is omitted.

\begin{lemma}{ Let $f(\xi) = \xi^4 + a_3\xi^3+a_2\xi^2+a_1 \xi + a_0$ and $ g(\xi) = b(\xi^4 + \beta_1 \xi^2 +\beta_0) $ be two polynomials, where $\beta_0$ and $\beta_1$ are positive. Then there exists a positive constant $b_0$ such that $f(\xi) \leq g(\xi)$ for all $\xi$ if $b \geq b_0$.
}
\end{lemma}

By applying this lemma to $|Q(\lambda)|^2$ and $\xi^4+(D_3+TT_3+S)\xi^2 +D_3S+DT_3$, we see that
\begin{equation*}
|Q(\lambda)|^2 \leq C_1\big(\xi^4 + (D_3(x)+T(x)T_3(x)+S(x))\xi^2 + D_3(x) S(x)+D(x)T_3(x) \big)
\end{equation*}
for some positive constant $C_1$. Consequently we obtain
\begin{equation}
{\rm Re}\, q(x,\lambda) \geq \frac{\xi+3\kappa_0}{C_1}
\label{A4.31}
\end{equation}
provided that ${\rm Re}\, \lambda \geq -3\kappa_0$ and $|{\rm Im}\, \lambda| \leq \mathit{ \Gamma}_{1}$. Hence,
\begin{equation}
{\rm Re}\, \mathcal{B}(\psi,\psi) 
 \geq D \Vert \psi^\prime \Vert_{L^2(I)}^2 + C_1^{-1}({\rm Re}\,\lambda +3\kappa_0) \Vert \psi \Vert_{L^2(I)}^2
\geq \min\{D, \kappa_0/C_1\} \Vert \psi \Vert_{H^1(I)}^2
\label{A4.32}
\end{equation}
whenever ${\rm Re}\, \lambda \geq -2\kappa_0$ and $|{\rm Im}\, \lambda| \leq \mathit{ \Gamma}_{1}$. Combining \eqref{A4.32} with \eqref{A4.24} completes the proof of the coerciveness estimate \eqref{A4.201}, hence Sublemma \ref{SL4.3}. \qed

{\it Step 3.}\quad
By the Lax-Milgram theorem, for any $h \in L^2(I)$ there exists a unique $\psi_h \in H^1(I)$ such that
\begin{equation}
\mathcal{B}(\psi_h, \psi) = (h, \psi)_{L^2(I)} \quad \hbox{for all}\ \psi \in H^1(I),
\label{A4.33}
\end{equation}
as long as ${\rm Re}\, \lambda + 2\kappa_0 \geq 0$ or $|{\rm Im}\, \lambda| \geq \mathit{\Gamma}_{1}$. 
We would like to prove that there exists a positive constant $M_1$ such that
\begin{equation*}
\Vert \psi_h \Vert_{L^2(I)} \leq \frac{M_1}{|\lambda+\kappa_0|+\kappa_0} \Vert h \Vert_{L^2(I)}
\end{equation*}
whenever ${\rm Re}\,\lambda \geq - \kappa_0$.

By \eqref{A4.33} with $\psi=\psi_h$, \eqref{A4.24} and \eqref{A4.32} we obtain
\begin{equation}
D \Vert \psi_h^\prime \Vert_{L^2(I)}^2 + \gamma_0 ({\rm Re}\,\lambda + 2 \kappa_0) \Vert \psi_h \Vert_{L^{2}(I)}^2 \leq {\rm Re}\, \mathcal{B}(\psi_h, \psi_h) \leq \Vert h \Vert_{L^2(I)} \Vert \psi_h \Vert_{L^2(I)}.
\label{A4.34}
\end{equation}
As in the computation of \eqref{A4.20}, we see that
\begin{equation}
\begin{aligned}
{\rm Im}\, P(\lambda)Q(\overline{\lambda}) & = - ({\rm Im}\, \lambda) \left\{ (D_3(x)S(x)-D(x)T_3(x)) + (D(x)-D_3(x)T(x)(\lambda+\overline{\lambda}) \right. \\
& \quad \left. - D_3(x)(\lambda^3-\overline{\lambda}^3) - (T(x)T_3(x)-S(x))|\lambda|^2 - T_3(x)|\lambda|^2(\lambda+\overline{\lambda}) + |\lambda|^4\right\},
\end{aligned}
\label{A4.35}
\end{equation}
which yields
\begin{equation} \label{A4.36}
{\rm Im}\, q(x,\lambda) = {\rm Im}\, \lambda (1+O(|\lambda|^{-1})) \quad \hbox{as}\ |\lambda| \to +\infty.
\end{equation}
Hence there exists a positive constant $\mathit{\Gamma}_{2} ( \geq \mathit{\Gamma}_{1})$ such that if $|{\rm Im}\,\lambda| \geq \mathit{\Gamma}_{2}$ then
\begin{equation*}
\left| {\rm Im}\, \mathcal{B}(\psi_h, \psi_h) \right| = \left| \int \limits_I {\rm Im}\, q(x,\lambda) |\psi_h|^2\,dx \right| 
\geq \frac{|{\rm Im}\,\lambda|}{2} \Vert \psi_h \Vert_{L^2(I)}^2.
\end{equation*}
Combining this estimate with \eqref{A4.33}, we conclude that
\begin{equation} \label{A4.37}
|{\rm Im}\,\lambda| \, \Vert \psi_h \Vert_{L^2(I)}^2 \leq 2 \Vert h \Vert_{L^2(I)} \Vert \psi_h \Vert_{L^2(I)}.
\end{equation}
From \eqref{A4.34} and \eqref{A4.37} we immediately obtain, as long as $|{\rm Im}\, \lambda| \geq \mathit{\Gamma}_{2}$,
\begin{equation} \label{A4.38}
(\gamma_0^2({\rm Re}\,\lambda+2\kappa_0)^2+({\rm Im}\,\lambda)^2) \Vert \psi_h \Vert_{L^2(I)}^2 \leq 5 \Vert h \Vert_{L^2(I)}^2.
\end{equation}
If ${\rm Re}\,\lambda + \kappa_0 \geq 0$, then $({\rm Re}\, \lambda+2\kappa_0)^2 \geq ({\rm Re}\, \lambda + \kappa_0)^2 +\kappa_0^2$. Hence, \eqref{A4.38} implies  
\begin{equation}
\Vert \psi_h \Vert_{L^2(I)} \leq \frac{m_1}{|\lambda+\kappa_0|+\kappa_0} \Vert h \Vert_{L^2(I)}
\label{A4.39}
\end{equation}
for ${\rm Re}\,\lambda +\kappa_0 \geq 0$ and $|{\rm Im}\, \lambda| \geq \mathit{\Gamma}_{2}$. Then from \eqref{A4.34} and \eqref{A4.39} it follows that
\begin{equation}
D \Vert \psi_h^\prime \Vert_{L^2(I)}^2 \leq \frac{m_1}{|\lambda+\kappa_0|+\kappa_0} \Vert h \Vert_{L^2(I)}^2.
\label{A4.40}
\end{equation}

Finally we treat the case $|{\rm Im}\, \lambda| \leq \mathit{ \Gamma}_{2}$ and ${\rm Re}\,\lambda +\kappa_0 \geq 0$. Then we have ${\rm Re}\, \lambda + 2\kappa_0 \geq \gamma_1(|\lambda +\kappa_0| + \kappa_0)$ with $2\gamma_1^2 = \kappa_0^2/(\mathit{\Gamma}_1^2+\kappa_0^2)$. Hence, \eqref{A4.34} yields the desired estimate
\begin{equation}
D \Vert \psi_h^\prime \Vert_{L^2(I)}^2 +
\gamma_0 \gamma_1(|\lambda+\kappa_0| + \kappa_0) \Vert \psi_h \Vert_{L^{2}(I)}^2 \leq \Vert h \Vert_{L^2(I)} \Vert \psi_h \Vert_{L^2(I)}.
\label{A4.41}
\end{equation}

{\it Step 4.}\quad Now we are ready to finish the proof of Lemma \ref{L4.2}. Clearly from \eqref{A4.13} it follows that
\begin{equation}
\sup_{x \in I} |p_1(x,\lambda)| \leq C_2
\quad \hbox{and} \quad
\sup_{x \in I} |p_{2}(x,\lambda)| \leq C_2 
\label{A4.42}
\end{equation}
for some positive constant $C_2$ independent of $\lambda  \not\in \mathcal{R}$. Therefore, the right-hand side of \eqref{A4.12} is bounded by $C_3 (\Vert r_1 \Vert_{L^2(I)} + \Vert r_2 \Vert_{L^2(I)}) + \Vert s \Vert_{L^2(I)}.$

Hence, the unique solution $\psi$ of \eqref{A4.12} under homogeneous Neumann boundary conditions satisfies the estimate
\begin{equation}
\Vert \psi \Vert_{L^2(I)} \leq \frac{M}{|\lambda+\kappa_0|+\kappa_0} \cdot \left\{ C_3  (\Vert r_1 \Vert_{L^2(I)} + \Vert r_2 \Vert_{L^2(I)}) +\Vert s \Vert_{L^2(I)} \right\}.
\label{A4.42a}
\end{equation}
From \eqref{A4.11} we see that
\begin{align*}
\Vert \phi_1 \Vert_{L^2(I)} 
& \leq c_{11}(\lambda) \Vert r_1 \Vert_{L^2(I)} +  c_{12}(\lambda) \Vert r_2 \Vert_{L^2(I)} + c_{13}(\lambda) \Vert \psi \Vert_{L^2(I)}, \\
\Vert \phi_2 \Vert_{L^2(I)} 
& \leq c_{21}(\lambda) \Vert r_1 \Vert_{L^2(I)} +  c_{22}(\lambda) \Vert r_2 \Vert_{L^2(I)} + c_{23}(\lambda) \Vert \psi \Vert_{L^2(I)}, \\
\end{align*}
where
\begin{align*}
c_{11}(\lambda) & =\left\Vert \frac{a_{22}-\lambda}{\det (A_{33}-\lambda)} \right\Vert_{L^\infty(I)}, \quad & 
c_{12}(\lambda) =\left\Vert \frac{a_{12}}{\det (A_{33}-\lambda)} \right\Vert_{L^\infty(I)}, \\
c_{13}(\lambda) & = \left\Vert \frac{ \det A_{31} + a_{13}\lambda }{ \det (A_{33}-\lambda)} \right\Vert_{L^\infty(I)}, \quad &
c_{21}(\lambda)  =\left\Vert \frac{a_{21}}{\det (A_{33}-\lambda)} \right\Vert_{L^\infty(I)}, \\
c_{22}(\lambda) & =\left\Vert \frac{a_{11}-\lambda}{\det (A_{33}-\lambda)} \right\Vert_{L^\infty(I)}, \quad &
c_{23}(\lambda) = \left\Vert \frac{ \det A_{32} - a_{23}\lambda }{ \det (A_{33}-\lambda)} \right\Vert_{L^\infty(I)}. 
\end{align*}
Since $| \det (A_{33}-\lambda) | \geq \gamma_* |\lambda+2\kappa|^2 $
 for all $\lambda \not\in \mathcal{R}$, there exists a positive constant $C_4$ such that
\begin{equation*}
c_{ij}(\lambda) \leq \frac{C_4}{|\lambda +\kappa|+\kappa} \quad \hbox{for}\ \lambda \not \in \mathcal{ R},\ i=1,\,2,\, 3; \ j=1,\,2.
\end{equation*}
Hence, by virtue of \eqref{A4.42a}, we see that
\begin{equation*}
\Vert \phi_1 \Vert_{L^2(I)} + \Vert \phi_2 \Vert_{L^2(I)} \leq \frac{C_5}{|\lambda+\kappa_0| + \kappa_0} \left( \Vert r_1 \Vert_{L^2(I)} + \Vert r_2 \Vert_{L^2(I)} + \Vert s \Vert_{L^2(I)} \right)
\end{equation*}
whenever ${\rm Re}\, \lambda + \kappa_0 \geq 0$. 
Putting these together, we conclude that if ${\rm Re}\,\lambda \geq -\kappa_0 $ then for any $(r_1, r_2, s) \in (L^2(I))^3$ equation \eqref{A4.7} has a unique solution $(\phi_1, \phi_2, \psi)$ and it satisfies the estimate
\begin{equation}
\Vert \phi_1 \Vert_{L^2(I)} + \Vert \phi_2 \Vert_{L^2(I)} + \Vert \psi \Vert_{L^2(I)} 
\leq \frac{C_5}{|\lambda + \kappa_0| + \kappa_0}.
\label{A4.888}
\end{equation}

Therefore, we have shown that $-\mathcal{ L}_2$ is sectorial in the sector $\{\lambda \in \mathbb{C} \mid |{\rm arg}\, (\lambda - \kappa_0)| \leq \pi/2 \}$. Once this is established, it is standard that the angle of $-\mathcal{ L}_2$ is less than $\pi/2$ (see, e.g., pp.~55--56 of \cite{Yagi}). This completes the proof of Lemma \ref{L4.2}.
\qed

\subsection{$(\varepsilon_0, A)$ stability} $\,$ \\
The initial-boundary value problem \eqref{full} is known to have a unique local-in-time classical solution $(u_1(t,x), u_2(t,x), v(t,x))$ if the initial data $(u_1^0, u_2^0, v^0)$ is sufficiently regular, e.g., $u_j^0 \in C^{\alpha}(\overline{I})$, $v^0 \in C^{2+\alpha}(\overline{I})$ and $dv^0/dx(0)=dv^0/dx(l)=0$ (see, e.g., \cite{R84} [Theorem 1, p.\ 111]). 
We may interpret this as a mild solution of the integral equation on the Hilbert space $\mathbb{X}=(L^2(I))^3$ as follows: Let $\Phi(t) = {}^T(\phi_1(t,\cdot), \phi_2(t,\cdot), \psi(t,\cdot))$. Then
\begin{equation}
\Phi(t) = e^{t \mathcal{L}_2}\Phi(0) + \int \limits_0^t e^{(t-\tau)\mathcal{L}_2} H(\Phi(s))\,d\tau, \qquad \Phi(0) \in (C(\overline{I}))^2 \times C^1(\overline{I}),
\label{A4.900}
\end{equation}
where we define the mapping $H : (L^\infty(I))^3 \to (L^2(I))^3 $ by
\begin{equation}
\left\{ \ 
\begin{aligned}
H(\Phi)&=\begin{pmatrix} r_1(\phi_1, \phi_2,\psi) \\ r_2(\phi_1,\phi_2,\psi) \\ s(\phi_1,\phi_2,\psi) \\ \end{pmatrix}, \\
r_1(\Phi) & = f_1(\tilde{u}_1+\phi_1,\tilde{u}_2+\phi_2,\tilde{v}+\psi, x) - a_{11}(x) \phi_1 - a_{12}(x) \phi_2 - a_{13}(x) \psi, \\
r_2(\Phi) & = f_2(\tilde{u}_1+\phi_1,\tilde{u}_2+\phi_2,\tilde{v}+\psi, x) - a_{21}(x) \phi_1 - a_{22}(x) \phi_2 - a_{23}(x) \psi, \\
s(\Phi) & =g(\tilde{u}_1+\phi_1,\tilde{u}_2+\phi_2,\tilde{v}+\psi, x) - a_{31}(x) \phi_1 - a_{32}(x) \phi_2 - a_{33}(x) \psi. \\
\end{aligned}
\right.
\label{A4.901}
\end{equation}
By virtue of \eqref{A4.6} of Lemma \ref{L4.2}, we have
\begin{equation*}
\Vert \Phi(t) \Vert_{\mathbb{X}} \leq M e^{-\kappa_1 t} \Vert \Phi(0) \Vert_{\mathbb{X}} + \int \limits_0^t M e^{-\kappa_1(t-\tau)} \Vert H(\Phi(\tau)) \Vert_{\mathbb{X}}\,d\tau.
\end{equation*}
Since $| H(\Phi) | \leq C (|\phi_1|+|\phi_2|+|\psi|)^2$, 
\begin{equation}
\Vert \Phi(t) \Vert_{\mathbb{X}} \leq M e^{-\kappa_1 t} \Vert \Phi(0) \Vert_{\mathbb{ X}} + CM \int \limits_0^t e^{-\kappa_1(t-\tau)} ( \Vert \phi_1(\tau) \Vert_{L^4(I)}^2 + \Vert \phi_2(\tau) \Vert_{L^4(I)}^2 + \Vert \psi(\tau) \Vert_{L^4(I)}^2)\,d\tau.
\label{A4.80}
\end{equation}

We consider the ODE subsystem, too:
\begin{equation}
\begin{aligned}
\pd{\phi_1}{t} & = a_{11}(x) \phi_1 + a_{12}(x) \phi_2 +  \rho_1(\phi_1, \phi_2, \psi),  \\
\pd{\phi_2}{t} & = a_{21}(x) \phi_1 + a_{22}(x) \phi_2 + \rho_2(\phi_1, \phi_2, \psi),
\end{aligned}
\label{A4.81}
\end{equation}
where
\begin{align*}
\rho_1(\phi_1,\phi_2,\psi) & = r_1(\phi_1, \phi_2, \psi) + a_{13}(x) \psi, \\
\rho_2(\phi_1, \phi_2, \psi) & = r_2(\phi_1, \phi_2,\psi) +a_{23}(x) \psi,
\end{align*}
Since this is a system of ordinary differential equations with constant coefficients for each $x$, the solution satisfies
\begin{align*}
\begin{pmatrix} \phi_1(t,x) \\ \phi_2(t,x) \\ \end{pmatrix} & = \exp \big(t A_{33}(x)\big) \begin{pmatrix} \phi_1(0,x) \\ \phi_2(0,x) \\ \end{pmatrix}  \\
& \quad + \int \limits_0^t \exp \big((t-\tau) A_{33}(x) \big) \begin{pmatrix} \rho_1(\phi_1(\tau,x),\phi_2(\tau,x), \psi(\tau,x)) \\ \rho_2(\phi_1(\tau,x), \phi_2(\tau,x), \psi(\tau,x)) \\ \end{pmatrix}\,d\tau.
\end{align*}
Moreover, by the assumptions on the matrix $A_{33}(x)$,
\begin{equation*}
\left\Vert \exp \big(tA_{33}(x)\big) \begin{pmatrix} \rho_1 \\ \rho_2 \\ \end{pmatrix} \right\Vert_{\R^2} \leq C e^{-\kappa t} \left\Vert \begin{pmatrix} \rho_1 \\ \rho_2 \\ \end{pmatrix} \right\Vert_{\R^2},
\end{equation*}
where $\Vert {}^T(\rho_1, \rho_2) \Vert_{\R^2} = |\rho_1|+|\rho_2|$. Hence, we have a pointwise estimate
\begin{equation}
\begin{aligned}
|\phi_1(t,x)|+|\phi_2(t,x)| & \leq M e^{-\kappa t} (|\phi_1(0,x)|+|\phi_2(0,x)|) \\
& \qquad + M \int \limits_0^t e^{-\kappa (t-\tau)} \sum_{j=1}^2 |\rho_j(\phi_1(\tau,x),\phi_2(\tau,x),\psi(\tau,x))|\,d\tau.
\end{aligned}
\label{A4.82}
\end{equation}

We now turn to the proof of Theorem \ref{thm:stabC-3c-gen}. In order to control the nonlinear terms, we first estimate the $H^1(I)$ norm of $\psi$, which gives us a bound on the $L^\infty(I)$ norm of $\psi$ by virtue of the Sobolev embedding theorem.

\begin{lemma} \label{L4.3}{  Let $\psi=\psi(t,x)$ be the solution of the initial-boundary value problem 
\begin{align}
& \pd{\psi}{t} - D \dpd{\psi}{x} = a_{31}(x) \phi_1 + a_{32}(x) \phi_2 + a_{33}(x) \psi + s(t,x) &\hbox{for} \ x \in I,\ t > 0, & \label{A4.50}\\
& \pd{\psi}{x}(t,x) = 0 & \hbox{for} \ x \in \partial I,\ t>0, & \label{A4.51}\\
& \psi(0,x) = \psi_0(x) &\hbox{for} \ x \in \overline{I}, & \label{A4.52}
\end{align}
where $s \in C([0, T] ; H^1(I)) \cap C^1((0,T] ; L^2(I))$. Then there exist positive constants $k$ and $M_1$ such that
\begin{equation}
\Vert \psi(t,\cdot) \Vert_{H^1(I)}^2 \leq M_1 \bigg\{ \Vert \psi_0 \Vert_{H^1(I)}^2 e^{- k t} + \int \limits_0^t  e^{-k (t-\tau)} \, d\tau \int \limits_I \bigg(\sum_{j=1}^2 \phi_j(\tau,x)^2 + s(\tau,x)^2 \bigg)\,dx \bigg\}.
\label{A4.53}
\end{equation}
}
\end{lemma}

{\it Proof}\/.  Let $k$ be a positive number to be determined later. We multiply both sides of \eqref{A4.50} by $e^{k t} \psi(t,x)$ and then integrate over $(0, t) \times I$. Observing that $\psi \psi_t = 2^{-1}(\psi^2)_t $ and $ \psi \psi_{xx} = (\psi \psi_x)_x - (\psi_x)^2$, we obtain after integration by parts:
\begin{equation}
\begin{aligned}
& \frac{e^{k t}}{2}  \int \limits_I \psi(t,x)^2\,dx - \frac{1}{2}\int \limits_I \psi_0(x)^2\,dx - \frac{k}{2} \int \limits_0^t e^{k \tau} \,d\tau \int \limits_I \psi(\tau,x)^2\,dx \\
& \quad 
+ D \int \limits_0^t e^{k \tau} \,d\tau \int \limits_I |\psi_x(\tau,x)|^2\,dx \\
& = \int \limits_0^t e^{k \tau} \,d\tau \int \limits_I \bigg( a_{33}(x) \psi(\tau,x)^2  + \sum_{j=1}^2 a_{3j}(x) \phi_j(\tau,x) \psi(\tau,x)  + s(\tau,x) \psi(\tau,x)\bigg)\,dx.
\end{aligned}
\label{A4.54}
\end{equation}
Similarly, we multiply both sides of \eqref{A4.50} by $e^{k t} \psi_t(t,x) $ and then integrate over $(0, t) \times I$. Noting that $\psi_{xx}\psi_t = (\psi_t \psi_x)_x - 2^{-1}(\psi_x^2)_t$, we obtain after integration by parts the following: 
\begin{equation}
\begin{aligned}
& \int \limits_0^t e^{k \tau} \,d\tau \int \limits_I \psi_t(\tau,x)^2 \,dx  + \frac{D}{2} \Big[ e^{k \tau} \int \limits_I \psi_x(\tau,x)^2\,dx \Big]_{\tau=0}^{\tau=t} - \frac{D k}{2} \int \limits_0^t e^{k \tau}\, d\tau \int \limits_I \psi_x(\tau,x)^2\,dx \\
& = - \frac{k}{2} \int \limits_0^t e^{k \tau} \,d\tau \int \limits_I a_{33}(x) \psi(\tau,x)^2\,dx + \frac{1}{2} \Big[ e^{k \tau} \int \limits_I a_{33}(x) \psi(\tau,x)^2\,dx \Big]_{\tau=0}^{\tau=t} \\
& \qquad + \int \limits_0^t e^{k \tau} \,d\tau \int \limits_I \bigg(\sum_{j=1}^2 a_{3j}(x) \phi_j(\tau,x)\psi_t(\tau,x) + s(\tau,x)\psi_t(\tau,x) \bigg)\,dx.
\end{aligned}
\label{A4.55}
\end{equation}

We use the well-known inequality $|ab| \leq \frac{1}{2}(\varepsilon a^2 + \varepsilon^{-1} b^2 )$ valid for arbitrary $\varepsilon>0$ to estimate the last two terms of \eqref{A4.54} and \eqref{A4.55} as follows:
\begin{equation}
\begin{aligned}
& \left| \int \limits_0^t e^{k\tau}\, d\tau \int \limits_I \bigg(\sum_{j=1}^2 a_{3j}(x) \phi_j(\tau,x) + s(\tau,x) \bigg)\psi(\tau,x) \,dx \right| \\
& \leq \frac{\varepsilon}{2} \int \limits_0^t e^{k\tau} \,d\tau \int \limits_I \psi(\tau,x)^2\, dx + \frac{1}{\varepsilon} R(t)
\end{aligned}
\label{A4.56}
\end{equation}
and
\begin{equation}
\begin{aligned}
& \left| \int \limits_0^t e^{k\tau} \,d \tau \int \limits_I \bigg(\sum_{j=1}^2 a_{3j}(x) \phi_j(\tau,x) + s (\tau,x) \bigg)\psi_t(\tau,x) \,dx \right| \\
& \leq \frac{\varepsilon^\prime}{2} \int \limits_0^t e^{k\tau} \, d\tau \int \limits_I \psi_t(\tau,x)^2\, dx + \frac{1}{\varepsilon^\prime} R(t), 
\end{aligned}
\label{A4.57}
\end{equation}
where
\begin{equation*}
R(t) = \int \limits_0^t e^{k\tau} \, d\tau \int \limits_I \bigg( \big( \sum_{j=1}^2 a_{3j}(x) \phi_j(\tau,x) \big)^2 + s(\tau,x)^2 \bigg) \,dx.
\end{equation*}
From \eqref{A4.54} and \eqref{A4.56} follows that
\begin{equation}
\begin{aligned}
& \frac{1}{2}e^{kt} \int \limits_I \psi(t,x)^2\,dx + D \int \limits_0^t e^{k\tau}\,d \tau \int \limits_I \psi_x(\tau,x)^2\,dx + \int \limits_0^t e^{k\tau}\, d\tau \int \limits_I |a_{33}(x) | \psi(\tau,x)^2\,dx  \\
& \leq \frac{1}{2} \int \limits_I \psi_0(x)^2\,dx + \frac{k}{2} \int \limits_0^t e^{k\tau}\,d\tau \int \limits_I \psi(\tau,x)^2\,dx  + \varepsilon \int \limits_0^t e^{k\tau}\,d\tau  \int \limits_I \psi(\tau,x)^2\,dx + \frac{1}{2\varepsilon} R(t). 
\end{aligned}
\label{A4.58}
\end{equation}
From \eqref{A4.55} and \eqref{A4.57} it follows that
\begin{equation}
\begin{aligned}
& \int \limits_0^t e^{k\tau}\,d\tau \int \limits_I \psi_t(\tau,x)^2\,dx + \frac{D}{2} e^{kt} \int \limits_I \psi_x(t,x)^2\,dx + \frac{e^{kt}}{2} \int \limits_I | a_{33}(x)| \psi(t,x)^2\,dx \\
& \leq \frac{D}{2} \int \limits_I \psi_0^\prime(x)^2\,dx + \frac{1}{2} \int \limits_I | a_{33}(x)|\psi_0(x)^2\,dx + \frac{Dk}{2} \int \limits_0^t e^{k\tau}\, d\tau \int \limits_I \psi_x(\tau,x)^2\,dx \\ 
& \quad +\frac{k}{2} \int \limits_0^t e^{k\tau}\, d\tau \int \limits_I | a_{33}(x)| \psi(\tau,x)^2 \,dx + \varepsilon^\prime \int \limits_0^t e^{k\tau}\, d\tau \int \limits_I \psi_t(\tau,x)^2\,dx + \frac{1}{2\varepsilon^\prime} R(t).
\end{aligned}
\label{A4.59}
\end{equation}
Adding \eqref{A4.58} and \eqref{A4.59} together results in
\begin{equation*}
\begin{aligned}
& \frac{1}{2} e^{kt} \int \limits_I \psi(t,x)^2\,dx + \frac{D}{2} e^{kt} \int \limits_I \psi_x(t,x)^2\,dx + (1-\varepsilon^\prime) \int \limits_0^t e^{k\tau}\, d\tau \int \limits_I \psi_t(\tau,x)^2\,dx \\
& \quad + D \left(1-\frac{k}{2}\right) \int \limits_0^t e^{k\tau}\,d \tau \int \limits_I \psi_x(\tau,x)^2\,dx + \left(1-\frac{k}{2}\right) \int \limits_0^t e^{k\tau}\,d \tau \int \limits_I |a_{33}(x)| \psi(\tau,x)^2\,dx  \\
& \qquad - \left(\frac{k}{2}+\varepsilon \right) \int \limits_0^t e^{k\tau}\, d\tau  \int \limits_I \psi(\tau,x)^2\,dx + \frac{e^{kt}}{2} \int \limits_I |a_{33}(x)| \psi(t,x)^2 \,dx  \\
& \leq \frac{1}{2} \left( \int \limits_I \psi_0(x)^2\,dx + \int \limits_I | a_{33}(x) | \psi_0(x)^2\,dx + D \int \limits_I \psi_0^\prime(x)^2\,dx \right)
+ \left( \frac{1}{2\varepsilon}+ \frac{1}{2\varepsilon^\prime} \right) R(t).
\end{aligned}
\end{equation*}
Recall from \eqref{A2.14} that $| a_{33}(x) | \geq 3\kappa $ for all $x \in \overline{I}$. Hence we choose $k=\varepsilon=\min\{\kappa,1\}$ to obtain
\begin{align*}
& \left(1-\frac{k}{2}\right) \int \limits_I | a_{33}(x)| \psi(\tau,x)^2\,dx - \left(\frac{k}{2} + \varepsilon \right) \int \limits_I \psi(\tau,x)^2\,dx \\
&
\geq \frac{3}{2}\left((2-k)\kappa -k\right) \int \limits_I \psi(\tau,x)^2\,dx \geq 0.
\end{align*}
Also we choose $\varepsilon^\prime = 1/2$ and conclude that
\begin{align*}
&\frac{1}{2} \int \limits_I \psi(t,x)^2\,dx + \frac{D}{2} \int \limits_I \psi_x(t,x)^2\,dx + \frac{1}{2} \int \limits_0^t e^{-k(t-\tau)}\, d\tau \int \limits_I \psi_t(\tau,x)^2\,dx \\
& \quad + \frac{D}{2} \int \limits_0^t e^{-k(t-\tau)} \,d\tau  \int \limits_I \psi_x(\tau,x)^2\,dx +0\\
& \leq \frac{1}{2} \left( \int \limits_I \psi_0(x)^2\,dx + \int \limits_I | a_{33}(x)| \psi_0(x)^2\,dx + D \int \limits_I \psi_0^\prime(x)^2\,dx \right) e^{-kt} \\
& \qquad + \left( \frac{1}{2k} + 1 \right) \int \limits_0^t e^{-k(t-\tau)}\, d\tau  \int \limits_I \bigg( \big(\sum_{j=1}^2 a_{3j}(x) \phi_j(\tau,x) \big)^2 + s(\tau,x)^2 \bigg)\, dx,
\end{align*}
which implies \eqref{A4.53}. 
\qed

\medskip
Put $X(t,x)=|\phi_1(t,x)|^2+|\phi_2(t,x)|^2$, $Y(t,x)=|\psi(t,x)|^2$ and $Z(t)=\Vert \psi(t,\cdot) \Vert_{H^1(I)}^2$. So far we have established the following inequalities:
\begin{align}
& X(t,x) \leq C X(0,x) + C \int \limits_0^t e^{-k(t-\tau)}\big( X(\tau,x)^2+Y(\tau,x) \big)\,d\tau, & \label{A4.100} \\
& \int \limits_I (X(t,x)+Y(t,x))\,dx \leq C \int \limits_I (X(0,x)+Y(0,x))\,dx e^{-kt} & \label{A4.101} \\
& \qquad \qquad \qquad \qquad \qquad + C \int \limits_0^t e^{-k(t-\tau)}\,d\tau \int \limits_I \big( X(\tau,x)^2 + Y(\tau,x)^2\big)\,dx, \nonumber \\
& Z(t) \leq CZ(0)e^{-kt} + C \int \limits_0^t e^{-k(t-\tau)}\,d\tau \int \limits_I \big(X(\tau,x)+ X(\tau,x)^2 + Y(\tau,x)^2 \big)\,dx. & \label{A4.102}
\end{align}
Also, by the Sobolev embedding theorem, 
\begin{equation}
 \sup_{x \in I} Y(t,x) \leq C Z(t).  
\label{A4.103}
\end{equation}
Here and in what follows, we use $C$ to denote various positive constants independent of $t$, $X(t)$, $Y(t)$ and $Z(t)$.
From \eqref{A4.101} we get
\begin{align*}
\int \limits_0^t e^{-k(t-\tau)}\, d\tau \int \limits_I X(\tau,x)\,dx & \leq C \int \limits_0^t e^{-k(t-\tau)} e^{-k\tau} \,d\tau \int \limits_I (X(0,x)+Y(0,x))\,dx \\
& + C \int \limits_0^t e^{-k(t-\tau)}\,d\tau \int \limits_0^\tau e^{-k(\tau-\sigma)}\,d\sigma \int \limits_I \big( X(\sigma,x)^2+Y(\sigma,x)^2\big)\,dx.
\end{align*}
It is straightforward to see that the right-hand side is equal to
\begin{align*}
 C t e^{-kt} \int \limits_I (X(0,x)+Y(0,x))\,dx + C \int \limits_0^t (t-\tau) e^{-k(t-\tau)}\,d\tau \int \limits_I \big( X(\tau,x)^2+Y(\tau,x)^2\big)\,dx
\end{align*}
  Since $(t-\tau) e^{-k(t-\tau)} \leq C_{k^\prime} e^{-k^\prime(t-\tau)} $ for $k^\prime \in (0, k)$, we obtain
\begin{align*}
& \int \limits_0^t e^{-k(t-\tau)}\,d\tau \int \limits_I X(\tau,x)\,dx  \\
& \leq C \int \limits_I\big(X(0,x)+Y(0,x)\big)\,dx + C \int \limits_0^t e^{-k^\prime(t-\tau)}\,d\tau \int \limits_I \big(X(\tau,x)^2 + Y(\tau,x)^2 \big)\,dx. 
\end{align*}

Hence, from \eqref{A4.102} we conclude that
\begin{equation}
\begin{aligned}
Z(t) & \leq C \big( Z(0) + \int \limits_I (X(0,x)+Y(0,x))\,dx \big) \\
& \qquad + C \int \limits_{I} e^{-k^\prime (t-\tau)}\,d\tau \int \limits_I \big( X(\tau,x)^2+Y(\tau,x)^2\big)\,dx
\end{aligned}
\label{A4.104}
\end{equation}

For any measurable set $R \subset I$, let
\begin{equation}
\xi_R(t) = \sup_{x\in R} X(t,x).
\label{A4.105}
\end{equation}
Then from \eqref{A4.100} and \eqref{A4.103}, we have
\begin{equation}
\xi_R(t) \leq C \xi_R(0) + \frac{C}{k} \xi_R(t)^2 + C \int \limits_0^t e^{-k(t-\tau)} Z(\tau)\,d\tau.
\label{A4.106}
\end{equation}
Decomposing the integral over $I$ into the sum of those over $R$ and $I \setminus R$, and denoting the Lebesgue measure of $R$ by $\mu(R)$, we see
\begin{equation}
\int \limits_I X(t,x)\,dx \leq \mu(R) \xi_R(t) + \mu(I \setminus R) K \  \hbox{and} \int \limits_I X(t,x)^2\,dx \leq \mu(R) \xi_R(t)^2 + \mu(I \setminus R) K^2,
\label{A4.107}
\end{equation}
where $K=\overline{U} - \underline{U}$.

From \eqref{A4.104} and \eqref{A4.107} it follows that
\begin{align*}
Z(t)  & \leq C (Z(0)+\mu(R)\xi_R(0) + \mu(I\setminus R)K)
 \\
 & \qquad + C \int \limits_0^t e^{-k(t-\tau)} \big(\mu(R)\xi_R(\tau)^2 +  \mu(I\setminus R) K^2 + Z(\tau)^2 \big)\,d\tau,
\end{align*}
which implies
\begin{equation}
\sup_{0\leq \tau \leq t} Z(\tau) \leq C (Z(0)+\mu(R)\xi_R(0) + \mu(I \setminus R)K) + C\sup_{0\leq \tau \leq t} \big(\mu(R)\xi_R(\tau)^2+ \mu(I \setminus R) K^2 + Z(\tau)^2 \big).
\label{A4.108}
\end{equation}
From \eqref{A4.106} we obtain
\begin{align*}
\xi_R(t) + Z(t)  & \leq C\xi_R(0) + C \xi_R(t)^2 + C (Z(0)+\mu(R)\xi_R(0) + \mu(I \setminus R) K) \\
& \qquad + C \sup_{0 \leq \tau \leq t} \big( \mu(R) \xi_R(\tau)^2 + \mu(I \setminus R) K^2 + Z(\tau)^2 \big),
\end{align*}
that is,
\begin{equation*}
\xi_R(t) +Z(t) \leq C(\xi_R(0)+Z(0)) + C \sup_{0\leq \tau \leq t} (\xi_R(\tau)^2+Z(\tau)^2) + C \mu(I \setminus R).
\end{equation*}
Therefore, $m(t) = \max_{0\leq \tau \leq t} ( \xi_R(\tau)+Z(\tau) ) $ satisfies
$ m(t) \leq C \big( m(0)+m(t)^2 + \mu(I\setminus R)\big)$.
Now it is easy to show that, for any $A$ satisfying $ A > \max\{ C, 1 \}$, if we choose $\varepsilon_0>0$ so small that
\begin{equation*}
1+A^2 \varepsilon_0^2 + \varepsilon_0^2 < \frac{A}{C},
\end{equation*}
then for any $0<\varepsilon<\varepsilon_0$, as long as $m(0) \leq \varepsilon^2 $ and $\mu(I\setminus R) < \varepsilon^4$, the inequality $m(t) < A\varepsilon^2$ holds for all $t \geq 0$.
This completes the proof of Theorem \ref{thm:stabC-3c-gen}.
\qed

\medskip

Theorem \ref{thm:stabC} on two-component systems is proved in exactly the same way and hence we omit it.

{\it Proof of Corollary \ref{thm:stabC-3c}\/.} Define a matrix $A^\delta(x)$ by
\begin{equation*}
A^\delta(x) = \begin{pmatrix} a_{11}(x) & a_{12}(x) & a_{13}(x) \\ a_{21}(x)/\delta & a_{22}(x)/\delta & a_{23}(x)/\delta \\ a_{31}(x) & a_{32}(x) & a_{33}(x) \\ \end{pmatrix}.
\end{equation*}
Then $\det A^\delta(x) = \delta^{-1} \det A(x)$, $\det A_{11}^\delta(x) = \delta^{-1} \det A_{11}(x)$, $\det A_{22}^\delta(x) = \det A_{22}$ and $ \det A_{33}^\delta(x) = \delta^{-1}\det A_{33}(x)$. Moreover,  as $\delta \downarrow 0$ we have ${\rm tr}\, A^\delta(x) = \delta^{-1} a_{22}+O(1)$ and ${\rm tr}\, A_{22}^\delta(x) = \delta^{-1} a_{22}(x) + O(1)$.

Therefore, it is easy to check that under the assumptions of the corollary, conditions \eqref{A2.11}--\eqref{A2.14} are satisfied for $\delta>0$ sufficiently small.
\qed

\bigskip

{\it Proof of Theorem \ref{lem:stabtrans}\/.} To prove Theorem \ref{lem:stabtrans}, we show the conditions of Corollary \ref{thm:stabC-3c} are satisfied. The only conditions left to check are \eqref{A2.18}--\eqref{A2.20}. 
First we differentiate the identity $f_2(u_1, u_2^*(u_1,v),v) = 0$ with respect to $u_1$ and $v$ at $(u_1,v)=(\tilde u_1, \tilde v)$, respectively, and obtain
\begin{equation}
a_{21} + a_{22}\pd{u_2^*}{u_1}=0 \quad \hbox{and} \quad a_{23}+a_{22}\pd{u_2^*}{v} =0.
\label{A4.800}
\end{equation}
Also, by differentiating $f_1(u_1, u_2^*(u_1,v),v)$ and $g(u_1, u_2^*(u_1,v),v)$ with respect to $u_1$ and $v$ at $(u_1,v)=(\tilde u_1, \tilde v)$ respectively, we get
\begin{equation}
 \eb_{11}=a_{11}+a_{12}\pd{u_2^*}{u_1}, \quad 
 \eb_{12} = a_{12}\pd{u_2^*}{v} + a_{13}, \quad 
 \eb_{21}=a_{31}+a_{32}\pd{u_2^*}{u_1}, \quad 
 \eb_{22} = a_{32}\pd{u_2^*}{v} + a_{33}.
\label{A4.801}
\end{equation}
Then we have
\begin{equation}
\begin{aligned}
\det A_{33} & = a_{11}a_{22}-a_{12}a_{21} \\
& =a_{22} \left( a_{11}+a_{12}\pd{u_2^*}{u_1} \right) \\
& = a_{22} \eb_{11} > - c_1 a_{22} \qquad (\hbox{by \eqref{A2.5}}).
\end{aligned}
\label{A4.802}
\end{equation}
Similarly, by using \eqref{A4.800} and \eqref{A4.801} we see
\begin{equation}
\det A_{11} = a_{22} \eb_{22} > 0 \qquad (\hbox{by (3) and \eqref{A2.7}}).
\label{A4.803}
\end{equation}

It is straightforward to verify the formula
\begin{equation}
\det A = a_{22} (\eb_{11} \eb_{22}-\eb_{12} \eb_{21}) = a_{22} \det \MB,
\label{A4.804}
\end{equation}
which is negative because of \eqref{A2.7} and (2).
\qed

\bigskip
\subsection{Quasi-steady state reduction of the ODE subsystem} $\,$\\
In this subsection we address the question of whether Diffusion Driven Instability can be investigated or not based on the quasi-steady state approximation. 
For $\delta>0$, consider a system of type \begin{align}
& \pd{\ueps}{t} = f_1(\ueps, \veps,\weps), & \label{A4.300}\\
& \delta \pd{\veps}{t} = f_2(\ueps, \veps, \weps), & \label{A4.301} \\
& \pd{\weps}{t} = D \frac{\partial^2 v}{\partial x^2} + g(\ueps, \veps, \weps), & \label{A4.302} \\
& \pd{\weps}{x} = 0 \quad (x \in \partial \text{$I$}).\label{uredBC}
\end{align}
Let $(\overline{u}_1,\overline{u}_2,\overline{v})$ satisfy
\begin{align}
 f_1(\overline{u}_1,\overline{u}_2,\overline{v})&=f_2(\overline{u}_1,\overline{u}_2,\overline{v})=g(\overline{u}_1,\overline{u}_2,\overline{v})=0, \label{ubarstst}\\
 \del{f_2}{u_2}(\overline{u}_1,\overline{u}_2,\overline{v}) & <0.\label{ImplExist}
\end{align}
Then by the implicit function theorem there exists a smooth function $u_2^*(u_1,v)$ defined in the neighborhood $\mathcal{N}=\{(u_1,v) \, | \, |u_1-\overline{u}_1| < \eta, \, |v-\overline{v}|<\eta\}$ of
$(\overline{u}_1,\overline{v})$ such that
\begin{align}
 f_2(\ueps,u_2^*(\ueps,\weps),\weps) = 0 \quad \text{for } (\ueps,\weps) \in \mathcal{N} \quad \text{and} \quad u_2^*(\overline{u}_1,\overline{v})=\overline{u}_2.\nonumber
\end{align}
Hence, $(\overline{u}_1,\overline{u}_2,\overline{v})$ is a constant steady state of \eqref{A4.300}-\eqref{A4.302}, and $(\overline{u}_1,\overline{v})$ is a constant steady state of its quasi-steady state reduction:
\begin{align}
\pd{u_1^0}{t} & = f_1(u_1^0, u_2^*(u_1^0,v^0), v^0), & \label{A4.304} \\
\pd{v^0}{t} & = D \frac{\partial^2 v}{\partial x^2} + g(u_1^0, u_2^*(u_1^0,v^0),v^0), & \label{A4.305}\\
\pd{v^0}{x} & = 0 \quad (x \in \partial \text{$I$}). \label{redBC}
\end{align}
\begin{prop}\label{propDDIImply}
 Assume that \eqref{ubarstst} and \eqref{ImplExist} are satisfied. If the reduced system \eqref{A4.304}-\eqref{redBC} exhibits DDI at the constant steady state $(\overline{u}_1,\overline{v})$, then there
 exists a positive number $\delta^*$ such that the unreduced system \eqref{A4.300}-\eqref{uredBC} exhibits DDI at $(\overline{u}_1,\overline{u}_2,\overline{v})$ as long as $0 < \delta < \delta^*$.
\end{prop}
Let $\MB=(\eb_{ij})_{1 \leq i,j\leq 2}$ denote the Jacobian matrix around the equilibrium $(\overline{u}_1,\overline{v})$ of the kinetic system for the reduced system, 
i.e. \eqref{A4.304}--\eqref{A4.305} with $D=0$. Let $A^{\delta}=(a_{ij}^{\delta})_{1 \leq i,j\leq 3}$ denote the Jacobian matrix around $(\overline{u}_1,\overline{u}_2,\overline{v})$ of the kinetic
system for the unreduced system, i.e., \eqref{A4.300}-\eqref{A4.302} with $D=0$. We prove this proposition by demonstrating the following two lemmas, both of which are algebraic in nature.
First, in Lemma \ref{L4.7} we prove that if $\delta$ is sufficiently small then two of the three eigenvalues of $A^{\delta}$ remain in the neighborhood of the two eigenvalues of $\MB$, while one
of the eigenvalues of $A^{\delta}$ tends to $-\infty$ as $\delta \rightarrow 0$, which in particular shows that $(\overline{u}_1,\overline{u}_2,\overline{v})$ is stable under
spatially homogeneous disturbance, provided that $\delta$ is sufficiently small. Second, in Lemma \ref{L4.9} we show that DDI of the reduced system \eqref{A4.304}-\eqref{redBC} implies
instabilty of $(\overline{u}_1,\overline{u}_2,\overline{v})$ for all $D>0$. Combining these lemma proves Proposition \ref{propDDIImply}.

\begin{lemma} \label{L4.7}
 Let $\MB$ and $A^{\delta}$ be as above. Assume that $\partial f_2 / \partial u_2 (\overline{u}_1,\overline{u}_2,\overline{v})<0$. Let $\lambda_{0,1}$ and $\lambda_{0,2}$ be two eigenvalues of $B$. Then
 \begin{enumerate}[{\rm (1)}]
\setlength{\leftskip}{2pc}
   \item if $\lambda_{0,1} \neq \lambda_{0,2}$ then for each $j=1,\,2$ there exists an eigenvalue $\lambda_{\delta,j} $ of  $A^\delta$ such that $\lim_{\delta \to 0} | \lambda_{\delta,j}-\lambda_{0,j} | = 0$ for $j=1,\,2$;
  \item if $\lambda_{0,1} = \lambda_{0,2} $ then there exist two eigenvalues $\lambda_{\delta,j}$, ($j=1,\,2$), of $A^\delta$ such that $\lim_{\delta \to 0} | \lambda_{\delta,j} - \lambda_{0,1} | = 0$ for $j=1,\,2$;
    \item there exists an eigenvalue $\lambda_{\delta,3}$ of $A^{\delta}$ such that, $\lambda_{\delta,3}$ is real and $\lim_{\delta \rightarrow 0}\lambda_{\delta,3}=-\infty$.
   \end{enumerate}
\end{lemma} 

\proof
Assertions (1) and (2) are proved by applying Rouch\'e's theorem. We put $a_{ij}=a_{ij}^1$ and $A=(a_{ij})_{1 \leq i,j \leq 3}$. Then $a_{ij}^{\delta}=a_{ij}$ for $i = 1,3$ and $j=1,2,3$,
whereas $a_{2j}^{\delta}=a_{2j}/\delta$ for $j=1,2,3$. Let $\Phi(\lambda)$ and $\Psi(\lambda)$ be the characteristic polynomials for $\MB$ and $A^{\delta}$, respectively:
\begin{align*}
 \Phi(\lambda) =& \lambda^2 - (\eb_{11} + \eb_{22})\lambda + \det \MB,\\
 \Psi(\lambda) =& -\lambda^3 + (a_{11}+\delta^{-1} a_{22} + a_{33})\lambda^2 - (\delta^{-1} \det A_{11} + \det A_{22} + \delta^{-1} \det A_{33}) \lambda + \delta^{-1} \det A.
\end{align*}
We put $\gamma_1 =\partial u_2^*/\partial u_1 (\overline{u}_1,\overline{v})$ and 
$\gamma_2 =\partial u_2^*/\partial v (\overline{u}_1,\overline{v})$ to see that
\begin{align*}
 & \eb_{11} = a_{11} + \gamma_1 a_{12}, \quad \eb_{12} = \gamma_2 a_{12} + a_{13}, \quad \eb_{21} = a_{31} + \gamma_1 a_{32}, \quad \eb_{22} = \gamma_2 a_{32} + a_{33},\\
 & a_{21}+\gamma_1 a_{22} = 0 \quad \text{and} \quad \gamma_2 a_{22} + a_{23} = 0.
\end{align*}
Using these relations, we obtain \eqref{A4.801}, \eqref{A4.802} and \eqref{A4.803}, which yield $\delta \Psi(\lambda) = a_{22} \Phi(\lambda) - \delta (\lambda^3 - (a_{11} + a_{33}) \lambda^2 + \det A_{22} \lambda)$.
Recall that $a_{22}<0$ by assumption \eqref{ImplExist}. Hence, the characteristic equation for $A^{\delta}$ reads
\begin{equation}\label{charPola}
 a_{22} \Phi(\lambda) - \delta \psi(\lambda) = 0, \qquad \text{where } \psi(\lambda) = \lambda^3 - (a_{11} + a_{33}) \lambda^2 + \det A_{22} \lambda.
\end{equation}
Let $\lambda_{0,1}$ and $\lambda_{0,2}$ be the roots of $\Phi(\lambda)=0$. Fix a disc $D_{\rho} = \{ |\lambda|<\rho \}$ which contains both $\lambda_{0,1}$ and $\lambda_{0,2}$.
We apply Rouch\'e's theorem on $\overline{D}_{\rho}$ to show that \eqref{charPola} has exactly two roots $\lambda_{\delta,1}$ and $\lambda_{\delta,2}$ inside $D_{\rho}$ for 
sufficiently small $\delta$. Moreover, if $\lambda_{0,1} \neq \lambda_{0,2}$ then, by passing to a small disc around $\lambda_{0,j}$, we can prove that $\lambda_{\delta,j} \rightarrow \lambda_{0,j}$ as $\delta \rightarrow 0$ for each $j=1, 2$. This verifies (1) and (2).

In addition, one can easily find that 
\eqref{charPola} has a real root $\lambda_{\delta,3}=\delta^{-1} a_{22} + O(\delta^{-2/3})$ as $\delta \rightarrow 0$, proving (3). Therefore, the assertions of the lemma hold true.
\qed

\begin{rem}
 Since $\det A^{\delta}=\delta^{-1} \det A$, if $A^{\delta}$ is a stable matrix for $\delta>0$ sufficiently small, then $\det A < 0$  and hence $\det A^\delta < 0$ for all $\delta>0$.
\end{rem}

\begin{lemma} \label{L4.9}
 Under assumption \eqref{ImplExist}, if the reduced system \eqref{A4.304}-\eqref{redBC} exhibits DDI at $(\overline{u}_1,\overline{v})$, i.e.,
 \begin{equation}
  \eb_{11}>0, \quad \operatorname{tr} \MB < 0 \text{  and  } \det \MB > 0,
 \end{equation}
then $(\overline{u}_1,\overline{u}_2,\overline{v})$ is unstable as a steady state of the unreduced system \eqref{A4.300}-\eqref{uredBC} for all $D>0$ and $\delta>0$.
\end{lemma}

\proof
 In a way analogous to the proof of \eqref{A4.802}, we obtain 
 \begin{equation}
  \det A_{33}^{\delta} = \delta^{-1} \begin{pmatrix} a_{11} & a_{12} \\ a_{21} & a_{22} \end{pmatrix} = \delta^{-1} a_{22} \eb_{11} < 0,
 \end{equation}
by virtue of our assumptions. Hence $A_{33}^\delta$  has a positive eigenvalue $\lambda^*$ for any $\delta>0$. Therefore $(\overline{u}_1,\overline{u}_2,\overline{v})$ is unstable if $D>0$ by  the result in \S 2.1.2 of \cite{Klika}. 
\qed

\medskip

\noindent
{\it Remark.} \ If we assume, instead of \eqref{A4.304}--\eqref{A4.305}, that both species diffuse and consider the system
$$
(u_1^0)_t = D_1 (u_1^0)_{xx} + f_1(u_1^0, u_2^*(u_1^0,v^0), v^0), \quad \quad (v^0)_t = D (v^0)_{xx} + g(u_1^0, u_2^*(u_1^0,v^0), v^0)
$$
then the constant steady state $(\overline{u}_1, \overline{v})$ is destabilized only for $D>D_{\rm cr}$ with $D_{\rm cr}$ being a positive number. This is quite different from the case $D_1=0$ in which $(\overline{u}_1, \overline{v})$ is always unstable for any $D>0$. The same observation applies to the three component system \eqref{A4.300}--\eqref{A4.302}. See \cite{Sakamoto12} for the case where all three species diffuse.

\section{Proof of statements concerning the example system}\label{sec:exproof}
First, we prove uniform boundedness of solutions and classification of different branches.

\proof[Proof of Lemma \ref{lem:trivstab}] To obtain bounds on $u$ and $v$, we see that for $v \geq 0$ \begin{equation} \label{eq5.1}
 -(1+v)u\leq f_r(u,v) \leq \left(-1+m_1 \frac{u}{1+k u^2}\right)u.
\end{equation}
The first inequality of \eqref{eq5.1} proves that $u$ remains positive if $v\geq 0$.

Define
\begin{equation*}
 h(u):=-1+m_1 \frac{u}{1+k u^2}.
\end{equation*}
If $m_1 \geq 2 \sqrt{k}$ then $h(u)=0$ has two real roots, and
the larger one is $u=U_2 := \frac{1}{2k}(m_1+\sqrt{(m_1)^2-4k}) $. Then $u(0,x) < U_2$ implies $ u(t,x) < U_2 $ for all $t>0$, since $ f_r(u,v) < 0$ for $u > U_2$. 
To obtain boundedness of $v$, note
\begin{equation}
 D\dell{v}{x} -(\mu_3+u)v \leq D \dell{v}{x}+g_r(u,v) \leq D \dell{v}{x} -\mu_3 v + \frac{m_2}{k}
\end{equation}
for $u \geq 0$, proving the second inequality due to the maximum principle.

If $m_1 < 2 \sqrt{k}$, then $\max h(u)  < 0$. Choose a $\gamma_h$ satisfying $0<\gamma_h < -\max h(u)$. Then $f_r(u,v) \leq - \gamma_h u $  for all $u \geq 0$ and $v \geq 0$. Therefore $\partial u / \partial t \leq - \gamma_h u$, proving that $u(t,x) \to 0$ uniformly as $t \to \infty$. Since $g_r \leq - \mu_3 v + m_2 u^2/(1+u^2)$ and $m_2 u(t,x)^2/(1+u(t,x)^2) \to 0$ as $t\to\infty$, we conclude that $ \lim_{t\to\infty} v(t,x) = 0$, uniformly on $\overline{I}$.

\qed 

\proof[Proof of Lemma \ref{lem:3-ub}]
To obtain uniform boundedness of solutions, note
\begin{align}
 \del{\ueps}{t}                  &\geq \big( \theta_1 \frac{\ueps}{1+\kappa (\ueps)^2}-\nu_1-\beta \weps \big) \ueps,\\
 \delta \del{\veps}{t} &\geq -(\nu_2+\alpha)\veps,\\
 \del{\weps}{t}                  &\geq \gamma \dell{\weps}{x} - (\nu_3+ \beta \ueps)\weps,
\end{align}
and for $\delta \leq 1$
\begin{equation} \label{eq5.6}
\begin{aligned}
 \del{(\ueps+\delta\veps)}{t} &= -(\nu_1\ueps +\nu_2 \veps)+\theta_1 \frac{(\ueps)^2}{1+\kappa (\ueps)^2} \\
                         &\leq - \min(\nu_1,\nu_2)(\ueps+\delta\veps) + \frac{\theta_1}{\kappa}. 
\end{aligned}
\end{equation}
Consequently,
\begin{equation}
 \limsup_{t \rightarrow \infty}  \big( u_1^\delta + \delta u_2^\delta \big)  \leq \frac{\theta_1}{\min(\nu_1,\nu_2) \kappa}.
\end{equation}
In particular, if $ u_1^\delta(0,x) + \delta u_2^\delta(0,x) < \theta_1/ (\kappa \min(\nu_1, \nu_2)) $, then $u_1^\delta(t,x) + \delta u_2^\delta(t,x) < \theta_1/ (\kappa \min(\nu_1, \nu_2)) $ for all $t > 0$.
Therefore, if $u_2^\delta(0,x) < \theta_1 / (\delta \kappa \min (\nu_1, \nu_2))$ then we have
\begin{equation}
 \del{\weps}{t} \leq \gamma \dell{\weps}{x} - (\nu_3+\beta \ueps)\weps + \frac{\alpha \theta_1}{\kappa \delta \min(\nu_1,\nu_2)}+\frac{\theta_2}{\kappa},
\end{equation}
yielding the result for $\weps$.

Now assume $\theta_1 < 2\nu_1 \sqrt{\kappa} $ and put $h_1(\xi) = - \nu_1 + \theta_1 \xi / (1+\kappa \xi^2) $. Then $\max h_1(\xi)$ is negative, and hence we can choose $\gamma_1$ satisfying $\max h_1(\xi) < - \gamma_1 < 0$. Therefore, $-\nu_1 u_1^\delta + \theta_1 (u_1^\delta)^2/(1+\kappa (u_1^\delta)^2) \leq - \gamma_1 u_1^\delta $ for $ u_1^\delta \geq 0$. We thus find from \eqref{eq5.6} that $\partial (u_1^\delta + \delta u_2^\delta) / \partial t \leq - (\gamma_1 u_1^\delta + \nu_2 u_2^\delta) $, which yields the uniform convergence $u_1^\delta(t,x) + \delta u_2^\delta(t,x) \to 0$ as $t\to\infty$. Once this is established, we have $\partial v^\delta / \partial t \leq \gamma \partial^2 v^\delta / \partial x^2 - \nu_3 v^\delta + c(t)$ with  $\lim_{t\to\infty} c(t) = 0$ as in the proof of Lemma \ref{lem:trivstab}, yielding the result for $v^\delta$.
\qed

\proof[Proof of Lemma \ref{stst:f=0}] If $(u,v)$ is a steady state of system \eqref{eq1}-\eqref{bc}, it satisfies
\begin{equation}\label{f=0}
 \del{u}{t} = f_r(u,v)= -(1+v)u+m_1\frac{u^2}{1+ku^2} = 0.
\end{equation}
This is satisfied if and only if either $u=0$ or
\begin{equation}\label{f=0-2}
 -(1+v)+m_1 \frac{{u}}{1+k u^2}=0.
\end{equation}
Solving \eqref{f=0-2} for $u$ yields
\begin{align}
 u_{\pm}(v) =   \frac{1}{2k(1+v)} \left( m_1 \pm \sqrt{m_1^2-4k(1+v)^2} \right)
\end{align}
proving Lemma \ref{stst:f=0}. \qed \par

The spatially homogeneous steady states are defined as common roots of two polynomials $(1+k u^2)f_r(u,v)=0=(1+ku^2)g_r(u,v)$. We therefore want to distinguish the different branches by a criterion which is easier to handle and note that the branches of steady states can be characterized in larger generality:
\begin{lemma}\label{lem:char}  The solution $u(v)$ of $f_r(u,v)=0$ has three branches:
\begin{align}
 u_0(v) &= 0, \\
 u_{\pm}(v) &=  \frac{1}{2k(1+v)} \left( m_1 \pm \sqrt{m_1^2-4k(1+v)^2} \right).
\end{align}
For $u\neq 0$ and $0 \leq v < v_r=m_1/(2 \sqrt{k})-1$, the following are true:
\begin{enumerate} [{\rm (1)}]
\setlength{\leftskip}{2pc}
 \item $u(v)=u_+(v)$ if and only if $\frac{d}{dv}u < 0$;
 \item $u(v)=u_-(v)$ if and only if $\frac{d}{dv}u > 0$.
\end{enumerate}
\end{lemma}
\proof
 Note that
\begin{equation}
 u_+(v)= \frac{1}{2k(1+v)} \left( m_1 + \sqrt{m_1^2-4k(1+v)^2} \right) 
\end{equation}
implies $\frac{d}{dv}u_+<0$.
Combining this with 
\begin{equation}
 u_+ u_- = \frac{1}{k},
\end{equation}
yields
\begin{equation}
 \frac{d}{dv}u_->0.
\end{equation}
\qed 

Using this characterization, we can prove existence of spatially homogeneous steady states and classify them:

\proof[Proof of Lemma \ref{lem:homStSt}] Clearly, $(\overline{u}_0,\overline{v}_0)=(0,0)$ is a spatially homogeneous steady state. 

To obtain nontrivial spatially homogeneous steady states, observe that a solution $(u,v)$ is a spatially homogeneous steady state only if
\begin{equation}
 \del{u}{t} = f_r(u,v) = -(1+v)u+m_1 \frac{u^2}{1+ku^2}=0.
\end{equation}
Solving $f_r(u,v)=0$ yields either $u=0$ or
\begin{equation}\label{Isocl-eq1}
 \frac{u}{1+k u^2}=\frac{1}{m_1}(1+v).
\end{equation}
Inserting \eqref{Isocl-eq1} into 
\begin{equation}
\del{v}{t}=g_r(u,v)=-(\mu_3+u)v+m_2 \frac{u^2}{1+ku^2}=0
\end{equation}
and solving for $v$ yield
\begin{equation}\label{w:f=g}
 v_{f_r,g_r}(u)=\frac{m_2}{m_1}\frac{u}{\mu_3+(1-\frac{m_2}{m_1})u}.
\end{equation}
Note that $v_{f_r,g_r}(u)$ is $k$-independent.

\begin{figure}[ht]
  \centering
  \begin{tikzpicture}[domain=0:3]
\newcommand{\cm}{2}
\newcommand{\cn}{1.44}
\newcommand{\cmu}{4.2}
\newcommand{\ck}{0.1}
    \draw[very thin,color=gray] (0.0,0.0) grid (0,0);
    \draw[->] (0,0) -- (12,0) node[below] {$u$};
    \draw[->] (0,-1) -- (0,5) node[left] {$v$};
    \draw[samples=100,domain=0:6, dashed] plot(\x,   {min((\cm/\cn)*\x/(\cmu + (1-(\cm/\cn))*\x),5)});     \draw[samples=100,domain=0:11]         plot(\x,   {-1+\cn*\x/(1+\ck*\x*\x)});     \draw[samples=100,domain=0:11, dotted] plot(\x,   {\cm/(\cmu+\x)*\x*\x/(1+\ck*\x*\x)});     \draw (6,4.5) node[right] {$v_{f,g}$};
    \draw (5,1.8) node[right] {$v_{g=0}$};
    \draw (4.5,0.7) node[right] {$v_{f=0}$};
    \draw (10.8,0) node[below] {$\frac{m_1 \mu_3}{m_2-m_1}$};
    \draw (3.063,0) node[below] {$\frac{1}{\sqrt{k}}$};
\end{tikzpicture}
   \label{fig:ws}
\caption{Nullclines of $f_r$ and $g_r$ and $v_{f_r,g_r}$ for parameters $m_1 = 1.44, m_2 = 2, \mu_3 = 4.2, k=0.1$.}
\end{figure}
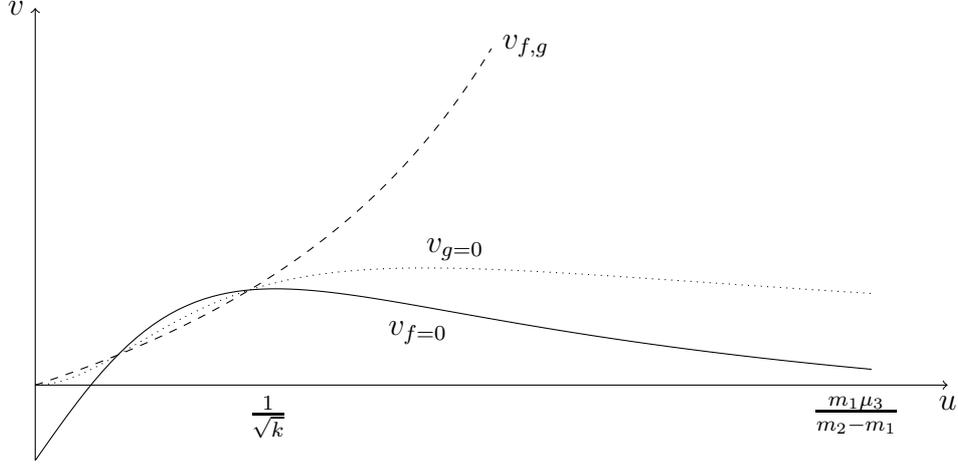

The nullclines of $f_r$ and $g_r$ for $u\neq0$ are described by
\begin{align}
 v_{f_r=0}(u):= -1 + m_1 \frac{u}{1+ku^2},\label{nc:f}\\
 v_{g_r=0}(u):= \frac{m_2}{\mu_3+u} \frac{u^2}{1+ku^2}.\label{nc:g}
\end{align}
If there exists $u>0$ such that $v_{g_r=0}(u)=v_{f_r,g_r}(u)$, then $(u,v_{f_r,g_r}(u))$ is a homogeneous steady state.

Note that $v_{f_r=0}(u)$ has a unique positive maximum at $u=1/\sqrt{k}$, is strictly concave on $[0,\sqrt{3}/\sqrt{k})$, strictly convex on $(\sqrt{3}/\sqrt{k}, \infty)$ and satisfies
\begin{equation}
 \lim_{u \searrow 0} v_{f_r=0}(u)=\lim_{u \nearrow \infty} v_{f_r=0}(u)=-1.
\end{equation}

Recall $m_1<m_2$. Defining $l(u):=m_2u^2/(\mu_3+u)$, we see that
\begin{align}
 l(u)-v_{g_r=0}(u)&=\frac{m_2}{\mu_3+u} \left(1-\frac{1}{1+ku^2}\right)u^2,\\
 l(u)-v_{g_r=0}(u)&>0,\\
 \lim_{k \rightarrow 0}\|l-v_{g_r=0}\|_{L^{\infty}([0,c])}&=0,
\end{align}
 on any finite interval $[0,c]$. Note also that
$v_{f_r,g_r}(u)$ is strictly increasing and continuous on $\mathbb{R}_{\geq 0}\setminus\{(m_1 \mu_3)/(m_2-m_1)\}$, nonnegative and convex on $[0,(m_1 \mu_3)/(m_2-m_1))$, negative on $((m_1 \mu_3)/(m_2-m_1),\infty)$ and satisfies 
\begin{align}
 \lim_{u \nearrow \frac{m_1 \mu_3}{m_2-m_1}} v_{f_r,g_r}(u) &= \infty,\label{asymp:f,g}\\
 \lim_{u \searrow \frac{m_1 \mu_3}{m_2-m_1}} v_{f_r,g_r}(u) &= -\infty.
\end{align}
Furthermore, it holds that
\begin{align}
 v_{f_r,g_r}(0)=v_{g_r=0}(0)=0,\\
 \frac{dv_{g_r=0}}{du}(0)=0<\frac{dv_{f_r,g_r}}{du}(0).
\end{align}
Hence, $v_{g_r=0}(\varepsilon)<v_{f_r,g_r}(\varepsilon)$ for $\varepsilon >0$ small.

Combining this with \eqref{asymp:f,g}, the uniform boundedness of $v_{g_r=0}$ and the strict convexity of $v_{f_r,g_r}$ on $(0,(m_1\mu_3)/(m_2-m_1))$, we see that if $(u,l(u))$ and $ u_, v_{f_r,g_r}(u))$ intersect twice in the region $\frak{R} = (0, (m_1 \mu_3)/(m_2-m_1)) \times (0, \infty)$, then there exists $k_1^*$ such that $(u,v_{g_r=0}(u))$ and $(u, v_{f_r,g_r}(u))$ intersect twice in the same region $\frak{R}$ for all $k<k_1^*$. Moreover, if $(u,l(u))$ and $(u, v_{f_r,g_r}(u))$ do not intersect twice, neither do $(u,v_{g_r=0}(u))$ and $ (u, v_{f_r,g_r}(u))$ for any $k \geq 0$.
Now we solve the $k$-independent equation $l(u)=v_{f_r,g_r}(u)$, i.e.,
\begin{equation}
 \frac{m_2 u^2}{\mu_3+u}=\frac{m_2}{m_1} \frac{u}{\mu_3+(1-\frac{m_2}{m_1})u} 
\end{equation}
for $u$ and obtain for $u\neq 0$:
\begin{equation}
 u_{\pm} = \frac{m_1\mu_3-1}{2(m_2-m_1)}\pm \sqrt{\left(\frac{m_1\mu_3-1}{2(m_2-m_1)}\right)^2 -\frac{\mu_3}{m_2-m_1}}.
\end{equation}
Clearly,
$u_{\pm}>0$ holds if and only if 
\begin{align}
 \mu_3&>\frac{1}{m_1}, \label{nonneg}\\
 \left(\frac{m_1\mu_3-1}{2(m_2-m_1)}\right)^2 &> \frac{\mu_3}{m_2-m_1}. \label{real}
\end{align}
Observe that \eqref{real} is equivalent to
\begin{equation} \label{eq5.35}
 \mu_3^2+2\frac{m_1-2 m_2}{m_1^2} \mu_3 + \frac{1}{m_1^2} > 0.
\end{equation}
Recall $m_2>m_1$ and
note that 
\begin{equation*}
 \left(\frac{m_2}{m_1}\right)^2-\frac{m_2}{m_1}>0
\end{equation*}
Hence, \eqref{eq5.35} is satisfied if and only if
\begin{equation}\label{mu:ub}
 \mu_3 < \frac{1}{m_1}\left(\frac{2 m_2 - m_1}{m_1} - 2 \sqrt{\left(\frac{m_2}{m_1}\right)^2-\frac{m_2}{m_1}}\right),
\end{equation}
or
\begin{equation}\label{mu:lb}
\mu_3 > \frac{1}{m_1}\left(\frac{2 m_2 - m_1}{m_1} + 2 \sqrt{\left(\frac{m_2}{m_1}\right)^2-\frac{m_2}{m_1}}\right).
\end{equation}

Inequality \eqref{mu:ub} is never satisfied if \eqref{nonneg} holds, because it holds that
\begin{equation}
 \frac{1}{m_1}\left(\frac{2 m_2 - m_1}{m_1} - 2 \sqrt{\left(\frac{m_2}{m_1}\right)^2-\frac{m_2}{m_1}}\right) < \frac{1}{m_1},
\end{equation}
implying $\mu_3 < 1/m_1$ if \eqref{mu:ub} holds, which contradicts \eqref{nonneg}.

On the other hand, since $m_2>m_1$, it holds that
\begin{equation}
\frac{2 m_2 - m_1}{m_1} + 2 \sqrt{\left(\frac{m_2}{m_1}\right)^2-\frac{m_2}{m_1}}>1.
\end{equation}
Therefore, if
\begin{align}
\mu_3 &> \frac{1}{m_1}\left(\frac{2 m_2 - m_1}{m_1} + 2 \sqrt{\left(\frac{m_2}{m_1}\right)^2-\frac{m_2}{m_1}}\right),\\
m_2 &> m_1,
\end{align}
both (5.33) and (5.34) are satisfied, and hence
 $(u,l(u))$ and $(u,v_{f,g}(u))$ intersect twice in $\frak{R}$.
This establishes the existence of two homogeneous steady states for sufficiently small $k$.

To prove that both steady states are of type $u_-$, we use the fact that the steady states component $u$ is less than $(m_1 \mu_3)/(m_2-m_1)$.
Since $v_{f_r=0}(u)=-1 + m_1 u/(1+k u^2)$ is strictly increasing on $(0,1/\sqrt{k})$, it holds that at any steady state $u^*$
\begin{equation}
 \frac{dv(u)}{du} (u^*)>0
\end{equation}
for $\sqrt{k}<(m_2-m_1)/(m_1 \mu_3)$.
By the characterization Lemma \ref{lem:char}, both spatially homogeneous steady states are of type $(u_-(v),v)$ for $k<((m_2-m_1)/(m_1 \mu_3))^2$.
\qed 

\medskip

\begin{figure}[ht]
  \centering
  \begin{tikzpicture}[domain=0:3,scale=1.5]
\newcommand{\cm}{2}
\newcommand{\cn}{1.44}
\newcommand{\cmu}{4.2}
\newcommand{\ck}{0.1}
\newcommand{\Cx}{(\cn/(2*(1+0.5*\x)))}
\newcommand{\upl}{(1/\ck *(\Cx + sqrt(\Cx*\Cx-\ck)) )}
\newcommand{\uminus}{(1/\ck *(\Cx - sqrt(\Cx*\Cx-\ck) ))}
\newcommand{\guplus}{(-(\cmu+\upl)*0.5*\x + \cm * \upl*\upl/(1+\ck*\upl*\upl))}
\newcommand{\guminus}{(-(\cmu+\uminus)*0.5*\x + \cm * \uminus*\uminus/(1+\ck*\uminus*\uminus))}
\newcommand{\gunull}{(-\cmu *\x)}
    \draw[very thin,color=gray] (0.0,0.0) grid (0,0);
    \draw[->] (0,0) -- (2.6,0) node[below] {$v$};
    \draw[->] (0,-1.5) -- (0,2.5);
    \draw[samples=100,domain=0:2.552] plot(\x,   {0.1*\guplus});     \draw[samples=100,domain=0:2.552, dashed] plot(\x,   {0.1*\guminus});     \draw[samples=100,domain=0:2.552, dotted] plot(\x,   {0.1*\gunull});
    \draw (1,-1) node[right]  {$g_r(u_0,v)$};
    \draw (1.4,-0.3)  node[right]  {$g_r(u_-(v),v)$};
    \draw (0.75,1.5) node[right] {$g_r(u_+(v),v)$};
\end{tikzpicture}
 \caption{Illustration of the right-hand side of $-\partial^2 v/\partial x^2 = g_r(u,v)$ for different branches of the solution $u(v)$ of $u_t=f_r(u,v)=0$. The parameters for illustration are $D=1,m_1=1.44,m_2=2,\mu_3=4.2$. We can observe that all nontrivial homogeneous steady states are of type $u_-$.}\label{fig:g}
\end{figure}
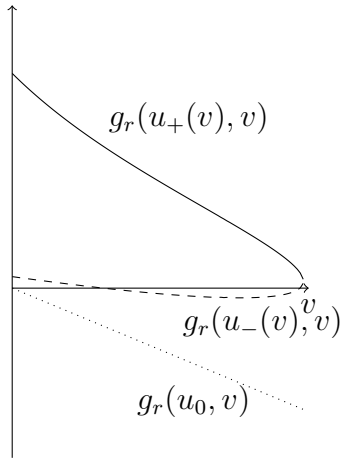

Existence of spatially inhomogeneous steady states can be proved similarly to \cite{MTH} or \cite{MCNT}.
To apply this work, we solve $f_r(u,v)=0$ in $u$ and substitute this into
\begin{equation} 
 g_r(u,v)= -(\mu_3+u)v+m_2 \frac{u^2}{1+k u^2}.
\end{equation}
The shape of $g_r(u(v),v)$ for different branches of the solution $u(v)$ of $f_r(u,v)=0$ is shown in Fig. \ref{fig:g}.

\begin{lemma}\label{lem:g} Assume $2 \sqrt{k}<m_1<m_2$. Then the following inequalities hold:
\begin{align} 
g_r(u_{\pm}(0),0)&>0, \label{g0}\\
g_r(u_+(v),v)&> g_r(u_{+}(v_r),v_r)=g_r(u_{-}(v_r),v_r) \quad \text{for $0<v<v_r$}, \label{g+>gwc}\\
g_r(u_+(v),v)&<g_r(u_+(0),0) \quad \text{for $0<v<v_r$}. \label{g+<g0}
\end{align}
If there exist two positive spatially homogeneous steady states of type $u_-$, it holds that
\begin{align}
g_r(u_{\pm}(v_r),v_r)&>0,\label{gwc}\\
g_r(u_+(0),0)&<\infty,\label{g0bdd}
\end{align}
and
\begin{equation}\label{derivg:unif}
\frac{d}{dv}g_r(u_-(v),v) \rightarrow \frac{1}{m_1} -\mu_3+\frac{2(m_2-m_1)}{m_1^2} (1+v),
\end{equation}
uniformly on any finite interval $v \in (0,c)$ as $k \searrow 0$.
\end{lemma}

\proof
Inserting $v=0$ into the right-hand side  (5.43) leads to
\begin{equation}\label{gpm0}
 g_r(u_{\pm}(0),0)=m_2 \frac{u_{\pm}(0)^2}{1+ku_{\pm}(0)^2}>0,
\end{equation}
yielding \eqref{g0}.

We differentiate $g_r(u_+(v),v)$ with respect to $v$ and obtain
\begin{equation}\label{d:g+}
\begin{aligned}
 \frac{d}{dv}g_r(u_+(v),v) &= -\frac{(1+v) m_1^2+2 k (1+v)^3 m_1 \mu _3}{2 k (1+v)^3 m_1}-\frac{m_1^3+4 k (1+v)^3 \left(m_2- m_1\right)}{2 k (1+v)^3 m_1\sqrt{\frac{-4 k (1+v)^2+m_1^2}{(1+v)^2}}}\\
                         &<0
\end{aligned}
\end{equation}
for $v \geq 0$, hence \eqref{g+>gwc} and \eqref{g+<g0} hold true.

If there exist two spatially homogeneous steady states of type $u_-$, then $g_r(u_{\pm}(v_r),v_r)$ must be positive because $g_r(u_-(0),0)>0$ and $g_r(u_-(v),v)$ has exactly two roots of order one in the interval $0<v<v_r$. Hence \eqref{gwc} is proved. 
Inequality \eqref{g0bdd} follows immediately from $u_+(0)\leq m_1/k$ and \eqref{gpm0}.

To see \eqref{derivg:unif}, we differentiate $g_r(u_-(v),v)$ with respect to $v$:
\begin{equation*}
\begin{aligned}
 \frac{d}{dv}g_r(u_-(v),v) &= \frac{m_1^3+4k(1+v)^3 (m_2-m_1)}{m_1 2 k (1+v)^2 \sqrt{-4k (1+v)^2+m_1^2}}+\frac{-m_1^2-2k(1+v)^2m_1 \mu_3}{2 k (1+v)^2 m_1},\\
                         &= -\mu_3+\frac{m_1^3 - m_1^2 \sqrt{-4k(1+v)^2+m_1^2}}{2 m_1 k (1+v)^2 \sqrt{-4k(1+v)^2+m_1^2}}-\frac{2(1+v) (m_2-m_1)}{m_1 \sqrt{-4k(1+v)^2+m_1^2}}. 
\end{aligned}
\end{equation*}
The limits for $k \rightarrow 0$ of the first and the last term are clear. To obtain the limit of the second term, we use l'H\^opital's rule and obtain
\begin{equation}
\frac{\frac{d}{dk}(m_1^2-m_1 \sqrt{-4k(1+v)^2+m_1^2})}{\frac{d}{dk}(2k(1+v)^2 \sqrt{-4k(1+v)^2+m_1^2})} = \frac{m_1}{-6k(1+v)^2+m_1^2},
\end{equation}
leading to 
\begin{equation}
\lim_{k \searrow 0}\frac{d}{dv}g_r(u_-(v),v)=\frac{1}{m_1} - \mu_3 + \frac{2(m_2-m_1)}{m_1^2} (1+v),
\end{equation}
which proves \eqref{derivg:unif}.
\qed 

\medskip

Using this knowledge, especially about the sign of $g_r$, we can prove existence of steady states with jump discontinuity. The proof is oriented on the shooting method and follows the principle presented in \cite{MTH}, \cite{Weinberger:87} and \cite{MCNT}. However, for completeness we state

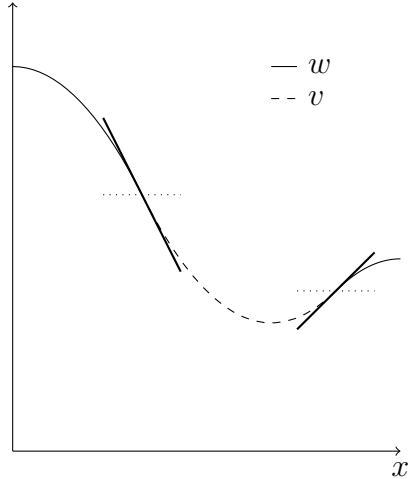
\begin{figure}[ht]
  \centering
  \begin{tikzpicture}[domain=-1:1,scale=1.7]
    \draw[->] (0,0) -- (0,3.5) node[left] {};
    \draw[->] (0,0) -- (3,0) node[below] {$x$};

    \draw[samples=20, domain=0:1] plot(\x, {2-\x*\x+1});
    \draw[dashed, samples=20, domain=1:2.5] plot(\x, {\x*\x-4*\x+4+1});
    \draw[samples=20, domain=2.5:3] plot(\x, {-\x*\x+6*\x-34/4+1});
    
    \draw[thick, samples=5, domain=0.7:1.3] plot(\x, {-2*\x+3+1});
    \draw[thick, samples=5, domain=2.2:2.8] plot(\x, {\x-9/4+1});
    
    \draw[dotted, samples=5, domain=0.7:1.3] plot(\x, {1+1});
    \draw[dotted, samples=5, domain=2.2:2.8] plot(\x, {1/4+1});

    \draw[-] (2,3)--(2.2,3) node[right] {$w$};
    \draw[-,dashed] (2,2.75)--(2.2,2.75) node[right] {$v$};
  \end{tikzpicture}
 \caption{Illustration of the construction of weak steady states in the proof of Lemma \ref{lem:disc}.}
  \label{fig:illu_shoot}
\end{figure}

\proof[Proof of Lemma \ref{lem:disc}] The concept of the proof is illustrated in Figure \ref{fig:illu_shoot}.
{\it Step 1\/.} 
We construct steady-states of \eqref{eq1}-\eqref{bc} by using two building blocks. 
The first one is the solution of the initial value problem
\begin{equation}
\begin{aligned}
 v^{\prime\prime} (x) &= - g_r(0,v(x)) = \mu_3 v(x) \qquad \hbox{for}\ x>0,\\
 v^{\prime}(0) &= 0, \\
 v(0) &= a,
\end{aligned}
\end{equation}
which is solved explicitly:
\begin{equation}
v(x) = a \cosh (\sqrt{\mu_3}x).
\end{equation}We denote this solution by $v(x;a)$.
The second one is the solution of the initial value problem
\begin{equation}\label{5.56}
\begin{aligned}
 w^{\prime\prime}(x) &= - g_r(u_+(w(x)),w(x)) \qquad \hbox{for}\ x>0,\\
 w^{\prime}(0) &= 0, \\
 w(0) &= b.
\end{aligned}
\end{equation}
We denote this solution by $w(x;b)$, where
the initial data $w_0$ is taken from the interval $(0, v_r)$. Denote $-g_r(u_+(w),w) $ by $h(w)$, which is 
defined for $w \in [0, v_r]$. Then we know from Lemma \ref{lem:g} that
\begin{align}
 &h(0) < h(v_r) < 0,\\
 &\od{h}{w}(w) >0 \quad \hbox{for}\ \ w \in [0, v_r).
\end{align}
Clearly $w(x)$ is positive and monotone decreasing in a certain interval $0 \leq x < x_{b}$ where $w(x_{b})=0$. 
Note that $x_b < \sqrt{2b/|h(b)|} $ since $w(x)=b+h(w(\theta x))x^2/2$ with $0<\theta<1$.

{\it Step 2\/.} We would like to connect $w(x;b)$ with $v(x;a)$, or $v(x;a)$ with $w(x;b)$, so that the boundary 
conditions are satisfied. Fix a number $c$ in the interval $(0, b)$ and let $x_c =x_c(c,b) \in (0, x_{b})$ be such 
that $w(x_c;b)=c$. Now we connect $w(\cdot;b)$ to $v(\cdot;a)$ at $w=c$ as follows: We would like to choose an 
appropriate $a \in (0, c)$ and a $y_c=y_c(c,a) > 0$ so that the matching conditions
\begin{equation} \label{5.99}
v(- y_c ;a)=c \quad \hbox{and} \quad
v^\prime(-y_c;a) = w^\prime(x_c;b)
\end{equation}
are satisfied.
This is possible if and only if 
\begin{equation}\label{5.100}
| w^\prime(x_c(c,b);b) | < \sqrt{\mu_3} c
\end{equation}
is satisfied. Indeed, if \eqref{5.100} holds, then there is a unique $\xi_c < 0$ such that 
$\tanh \xi_c = w^\prime(x_c;b)/ (\sqrt{\mu_3} c) $. Then we obtain $ y_c = - \xi_c/\sqrt{\mu_3}$ and 
$ a = c / \cosh \xi_c$. We say that {\it $w(x;b)$ is switchable to $v(x)$ at $w=c$} if \eqref{5.99} is satisfied. 
We remark that for each $b$ there exists a unique $c^* =c^*(b) \in (0, b)$ such that $w(x;b)$ is switchable to 
$v(x)$ at $c$ if and only if $c^*(b)<c<b$. This is proved easily by using the following three facts: 
$(w^\prime/w)^\prime = h(w)/w - (w^\prime/w)^2 < 0$, $w^\prime(0)/w(0) = 0$ and 
$\lim_{x \downarrow x_{w_0}} w^\prime(x)/w(x) = -\infty$. 

Next, fix a number $\gamma \in (a, v_r)$ and let $y_c =y_c(\gamma,a) >0$ be such that $v(y_c;a)=\gamma$. 
We extend $w(x;b)$ to $-x_{b} \leq x \leq 0$ by defining $ w(x;b)=w(-x;b)$ for $x \in [-x_{b}, 0]$. 
This time we would like to find $b \in (a, v_r)$ and $x_c=x_c(\gamma,b) \in  (0, x_{b})$ such that
\begin{equation}\label{5.100a}
w(-x_c; b)= \gamma \quad \hbox{and} \quad
w^\prime(-x_c; b) = v^\prime(y_c; a),
\end{equation}
which is equivalent to the following conditions:
\begin{equation}\label{5.101}
w(x_c; b) = \gamma \quad \hbox{and} \quad w^\prime(x_c; b)= - v^\prime(y_c;\gamma).
\end{equation}
We say that {\it $v(x;a)$ is switchable to $w(x)$ at $v=\gamma$} if there exist $b$ and $x_c(\gamma,b)$ satisfying 
the matching condition \eqref{5.100a}. Likewise for $w(x;b)$, we shall prove that there exists a unique 
$c_*={c_*}(a)$ such that $v(x;a)$ is switchable to $w(x)$ at $v=\gamma$ if and only if $a <\gamma< c_*(a)$.

{\it Step 3\/.} Assuming the existence of $c_*(a)$, we describe how to construct a steady state of \eqref{eq1}-\eqref{bc}. 
First we choose $b \in (0, v_r)$ arbitrarily and then choose any $c \in (c^*(b), b)$. We now obtain $a$ and $y_c(c,a) $ 
which satisfy \eqref{5.99}. Let us put
\begin{equation*}
W(x)=\begin{cases}w(x;b) & \text{for } x \in [0, x_{c}(c,b) ), \\ 
      v(x-x_c(c,b)-y_c(c,a); a) & \text{for } x \in [x_{c}(c,b), x_{c}(c,b)+y_c(c,a)].
      \end{cases}
\end{equation*}
Clearly $W_1(x) = W((x_c(c,b)+y_c(c,a) )x) $ is a weak solution of \eqref{eq1}-\eqref{bc} for $D=1/(x_c+y_c)^2$. 
Next, choose $\gamma \in (a, {c_*}(a))$ arbitrarily and find a pair $(b_{1}, x_c(\gamma,b_1))$ satisfying \eqref{5.101}. We extend $W(x)$ to the interval $[0, x_c(c,b)+y_c(c,a)+y_{c}(\gamma,a) + y_{c}(\gamma,b_1)]$ by putting
\begin{equation*}
W(x) = \begin{cases}W(x) & \text{for } x \in [0, L_1], \\
              v(x-L_1; b_{1}) & \text{for } x \in [L_1, L_2], \\ 
              w(x-L_3;b_{1}) & \text{for } x \in [L_2, L_3],
       \end{cases}
\end{equation*}
where $L_1=x_c(c,b)+y_c(c,a)$, $L_2=L_1+y_{c}(\gamma,a)$ and $L_3=L_2+x_c(\gamma,b_1)$. 
Actually, $W(x)$ is extended to $[0, L_3 +x_{b_1}]$.  If we put $W_2(x)=W(L_3 x)$, then this gives rise to a weak 
solution of \eqref{eq1}--\eqref{bc}  for $D=1/L_3^2$, which has two switching points. 
We can repeat these procedures to obtain steady states of \eqref{eq1}--\eqref{bc} with an arbitrary number of 
switching points. In the same manner, starting with $v(x;a)$, one obtains another type of steady states with any number 
of switching points.

{\it Step 4\/.} Now we prove the existence of $(b,x_c(\gamma,b))$ satisfying \eqref{5.101}. Since $w(x;b)$ is not an elementary 
function, we investigate its qualitative properties. First, let $0<b_1< b_2 < v_r$. Then $w(x;b_1) < w(x;b_2) $ in the common 
interval of existence. For, $\phi(x) = w(x;b_2)-w(x;b_1)$ satisfies
$\phi^{\prime \prime} = h^\prime(w(x;b_1)+\theta \phi(x)) \phi, \ \phi^\prime(0) = 0$ and $\phi(0) =b_2 - b_1 > 0 $
with $0<\theta<1$. Since $h^\prime(w)>0$, this implies $\phi(x) > 0$ for $x >0$ as long as it exists. 
Next, fix $\gamma \in (0, v_r)$ and let $y_c =y_c(\gamma,a)>0$ be such that $v(y_c;a)=\gamma$. 
Then for each $ b \in (\gamma, v_r)$ there exists a unique $ x(\gamma;b) > 0$ such that $w(x(\gamma;b);b) = \gamma$. 
Notice that $x(\gamma,b)$ is a continuously differentiable function of $b$ in $(\gamma, v_r)$ due to the implicit function 
theorem, since $w(x(\gamma,b);b)$ $=\gamma$ and $ w^\prime(x(\gamma,b);b) < 0$. Moreover,
we claim that if $ \gamma < b_1 < b_2 < v_r$ then $ w^\prime(x(\gamma,b_1);b_1) > w^\prime(x(\gamma,b_2);b_2)$. 
To verify this assertion, we define
\begin{equation}
H(w) = \int \limits_0^w h(v)\,dv.
\end{equation}
Then from \eqref{5.56} it follows that
\begin{equation}
(w^\prime(x;b))^2 = 2\left( H(w(x;b)) - H(b) \right),
\end{equation}
and hence
\begin{equation}
(w^\prime(x(\gamma;b_2);b_2))^2 - (w^\prime(x(\gamma;b_1); b_1))^2 = - 2\left( H(b_2)-H(b_1) \right).
\end{equation}
Recalling (5.57) and $b_2> b_1$, we find that $(w^\prime(x(\gamma;b_2);b_2))^2 > (w^\prime(x(\gamma;b_1);b_1))^2$. 
Since $w^\prime(x;b_j) < 0$ for $x \in (0, x_{b_j})$, $j=1,\,2$, we conclude that 
$ w^\prime(x(\gamma;b_2);b_2) < w^\prime(x(\gamma; b_1); b_1)$.  

We have proved that $b \mapsto w^\prime(x(\gamma,b);b)$ is a strictly decreasing continuous function. 
It is easy to check that (i) $w^\prime(x(\gamma,b); b) \to 0$ as $b \downarrow 0$, 
(ii) $w^\prime(x(\gamma,b); b) \to - \sqrt{2(H(\gamma)-H(v_r))} $ as $ b \uparrow v_r$, and 
(iii) if we put $p_*(\gamma)=\sqrt{2(H(\gamma)-H(v_r))} $ then 
$p_*^\prime(\gamma)=h(\gamma)/p_*(\gamma) < 0$, $p_*(0)=\sqrt{-2H(v_r)} >0$, and $p_*(v_r)=0$. 
Therefore, given $q \in (0, p_*(\gamma))$, there exists a unique $b \in (\gamma, v_r)$ such that $w^\prime(x(\gamma,b),b)=-q$.

Hence, our problem reduces to asking whether $|v^\prime(y_c(\gamma,a);a)| \in (0, p_*(\gamma))$ or not, 
given $a$ and $\gamma$.  Since $v(x;a)=a \cosh (\sqrt{\mu_3} x)$, we see by elementary computations that 
$\sqrt{\mu_3} y_c(\gamma,a)=\cosh^{-1}(\gamma/a)$ and $$
v^\prime(y_c(\gamma,a);a) = - \sqrt{\mu_3}a \sinh^{-1}(\cosh^{-1}(\gamma/a)) = - \sqrt{\mu_3} \sqrt{\gamma^2-a^2}.
$$ 
These facts imply that there exists a unique $c_*(a) \in (a, v_r)$ such that 
$\sqrt{\mu_3(c_*(a)^2-a^2)} = p_*(c_*(a))$ and $|v^\prime(y_c(\gamma,a);a)| < p_*(\gamma)$ if $a<\gamma<c_*(a)$, 
while $|v^\prime(y_c(\gamma,a);a)| > p_*(\gamma)$ if $a \in (c_*(a), v_r)$. This proves the existence of $c_*(a)$.

{\it Step 5\/.}  Finally we examine the range of $D$ which has  monotone decreasing solutions. Note that $c^*$ is characterized by the property that $|w^\prime(x_c(c^*,b);b) | = \sqrt{\mu_3} c^*$. Therefore, $\tanh \xi_c \to 1$, i.e., $\xi_c \to -\infty$ as $c\downarrow c^*(b)$, so that $y_c \to +\infty$ and $a \to 0$. This means that $D=1/(x_c+y_c)^2 \to 0$ as $c \downarrow c^*(b)$. On the other hand, as $c \uparrow b$, we have (i) $x_c \to 0$ and (ii) $w^\prime(x_c;b) \to 0$ and hence $\xi_c \to 0$. This implies that $y_c \to 0$ and $a \to c$ as $c \uparrow 0$. In particular, $D=1/(x_c+y_c)^2 \to +\infty$ as $c\uparrow b$. Observing that $x_c$ and $y_c$ are continuous function of $c$, we conclude that for each $D>0$ and for each $b \in (0, v_r)$ there exists a monotone decreasing solution $W_1(x;b)$ of \eqref{eq1}-\eqref{bc} as stated in Step 3. Clearly $W_1(x;b_1) \neq W_1(x;b_2)$ if $b_1 \neq b_2$, we hence obtain infinitely many steady states for each $D>0$. (The arguments above also prove the existence of infinitely many steady states which are not monotone in $x$ for each $D>0$.)
This completes the proof of Lemma \ref{lem:disc}.
\qed

Now that we have established the existence of spatially homogeneous steady states and spatially inhomogeneous steady states,  we turn to the proof of the lemmas stating that the conditions of Theorem \ref{thm:stabC} are satisfied:

\proof[Proof of Lemma \ref{lem:entries}] First, we calculate the Jacobian matrix of the kinetic system for given arbitrary $(u,v)$:
\begin{equation}
 \MB(x):=
\begin{pmatrix}
 \eb_{11} & \eb_{12} \\
 \eb_{21} & \eb_{22}
\end{pmatrix}
:=
\begin{pmatrix}
  -(1+v)+m_1 \frac{2u}{(1+ku^2)^2} & -u \\
  -v + m_2 \frac{2u}{(1+ku^2)^2} & -(\mu_3+u)
\end{pmatrix}.
\end{equation}
The signs of $b_{ij}$ for $u=u_0=0$ follow immediately by inserting $u=0$; also the signs of $\eb_{12}$ and $\eb_{22}$ are obvious for arbitrary $u$.

For $f_r(u,v)=-(1+v)u+m_1 u^2/(1+k u^2)$ it follows
\begin{equation}\label{a11}
 \begin{aligned}
  \eb_{11} &= -(1+v)+ \frac{2 m_1 u}{(1+ku^2)^2}\\
         &= -(1+v)+ \frac{2}{1+ku^2} (1+v)\\
         &= \frac{1-ku^2}{1+k u^2}(1+v).
 \end{aligned}
\end{equation}
Combining \eqref{a11} with $u_-(v) \leq 1/\sqrt{k} \leq u_+(v)$ yields the result for $\eb_{11}$.

We investigate $\eb_{21}$:
\begin{equation}
\begin{aligned}
 \eb_{21} &=     -v + m_2 \frac{2u}{(1+ku^2)^2} \\
        &=     -v + \frac{m_2}{m_1}(1+v)+\frac{m_2}{m_1} \left(-1-v+m_1 \frac{2u}{(1+ku^2)^2}\right)\\
        &=     \frac{m_2}{m_1}+\left(\frac{m_2}{m_1}-1\right)v+\frac{m_2}{m_1} \eb_{11}.
\end{aligned}
\end{equation}
The continuity of $\eb_{21}$ and $\eb_{11}(v_r)=0$ imply the result for $\eb_{21}$.

In the proof of Lemma \ref{lem:g}, we saw in \eqref{d:g+} that
\begin{equation}
 \frac{d}{dv}g_r(u_+(v),v)<0,
\end{equation}
and in Lemma \ref{lem:char} that
\begin{equation}
 \frac{d}{dv}u_+(v)<0.
\end{equation}

It is left to prove the dependence of the determinant of $\MB(x)$ on $dg(u_{f_r=0}(v),v)/dv$. To see this, recall from the identity $f_r(u(v),v)=0$ that
\begin{equation} \label{eq.5.59}
\begin{aligned}
 & b_{11} \frac{du}{dv} + b_{12} =0,  \quad \text{and} \quad \frac{d}{dv} g_r(u(v),v) = \frac{d u}{dv} \eb_{21} + \eb_{22}. \\
\end{aligned}
\end{equation}
The first equation yields $du/dv = - b_{12}/b_{11}$. Hence, from the second equation we obtain
\begin{equation*}
\begin{aligned}
\frac{d}{dv} g_r(u(v),v) & = -\frac{b_{12}}{b_{11}} \eb_{21} + \eb_{22} =
 \frac{\det \MB}{b_{11}}.
\end{aligned}
\end{equation*}
From the first equation of \eqref{eq.5.59} we get the expression $b_{11}=-b_{12}/(du/dv)$. This in turn give us
\begin{equation*}
\begin{aligned}
\det \MB & = b_{11} \frac{d}{dv} g_r(u(v),v) \\
& = - b_{12}\left(\frac{du}{dv}\right)^{-1} \frac{d}{dv} g_r(u(v),v),
\end{aligned}
\end{equation*}
as desired.
\qed 

\proof[Proof of Lemma \ref{lem:stabKin}] First, we prove that spatially homogeneous steady states of type $(u_-(v),v)$ are unstable as steady states of system \eqref{eq1}-\eqref{bc} for any $D>0$. Lemma \ref{lem:entries} states
\begin{equation}
 \partial_u f_r|_{(u_-(v),v)}>0,
\end{equation}
proving instability since autocatalysis implies instability, see \cite{MCKS13}.

Letting $\MB$ again denote the Jacobian matrix of the kinetic system at the steady state, we recall that a steady state of the kinetic system is stable if and only if $\text{tr}(\MB)<0$ and $\operatorname{det}{\MB}>0$. 
Hence, we calculate the trace:
\begin{equation}
\begin{aligned}
 \text{tr}\, \MB &= \left(-1+\frac{2}{1+ku^2}\right)(1+v)-(\mu_3+u)\\
                &= \left(-1+\frac{2}{1+ku^2}-\frac{m_2}{m_1}\frac{u}{v}\right)(1+v)\\
\end{aligned}
\end{equation}
since $v(\mu_3+u)= m_2 u^2/(1+ku^2) = m_2 u(1+v)/m_1$ due to (5.18) and (5.19).
From 
\begin{equation}
 v = -1 + m_1 \frac{u}{1+ku^2} \leq -1+m_1 u,
\end{equation}
it follows that
\begin{equation}
\begin{aligned}
\eb_{11}+\eb_{22}           &\leq  \left(-1+2-\frac{m_2}{m_1}\frac{u}{-1+m_1 u}\right) (1+v).\\
\end{aligned}
\end{equation}
Therefore, it is easy to conclude that $\operatorname{tr} \MB \leq 0$ if $m_1 \leq \sqrt{m_2}$.

To prove $\operatorname{det} \MB >0$, we recall that $g_r(u_-(0) ,0) > 0$, $g_r(u_-(v_r),v_r)>0$  and $g_r(u_-(v),v)$ has exactly two roots $v_-<v_+ \in (0,v_r)$, see Lemma \ref{lem:g}.
It follows from the proof of Lemma \ref{lem:g} that
\begin{equation}\label{der:g}
\begin{aligned}
 \frac{dg_r(u_-(v),v)}{dv}(v_+)&>c_1(k),\\
 \frac{dg_r(u_-(v),v)}{dv}(v_-)&<-c_2(k).
\end{aligned}
\end{equation}

Note 
\begin{enumerate} [{\rm (1)}]
\setlength{\leftskip}{2pc}
 \item $\frac{d}{dv}u_-(v)>0$, see Lemma \ref{lem:char};
 \item $u_-(v)>1/m_1$, see proof of Lemma \ref{lem:char};
 \item $\operatorname{det} \MB = u \frac{d}{dv}g_r(u_-(v),v) / \frac{du_-}{dv}(v)$, see Lemma \ref{lem:entries}.
\end{enumerate}
Assuming $dg_r(u_-(v),v)/dv(v_{\pm}) \rightarrow \pm |c_{\pm}|$ as $k \searrow 0$,
\eqref{der:g} therefore implies
\begin{equation} \label{eq.5.65}
\begin{aligned}
\operatorname{det} \MB(u_-(v_+),v_+) >0,\\
\operatorname{det} \MB(u_-(v_-),v_-) <0.
\end{aligned}
\end{equation}
It is left to prove $dg_r(u_-(v),v)/dv(v_{\pm}) \not\rightarrow 0$ as $k \rightarrow 0$.
To see this, recall that $(u_-(v_-),v_-)$ and $(u_-(v_+),v_+)$ converge towards two different, finite limits as $k \searrow 0$, see the proof of Lemma \ref{lem:homStSt}. Since \eqref{derivg:unif} states that $dg_r(u_-(v),v)/dv$ converges uniformly towards an affine-linear function on any finite interval, we obtain the bound on the derivative.
\qed 

Finally we come to the \proof[Proof of Lemma \ref{lem:3er-sys}]
Assertion $(1)$ follows directly from Lemma \ref{lem:homStSt} since $f_2(\ueps,\veps,\weps)=0$ is uniquely solvable and the rescaling is linear with positive constants. 

To investigate stability of spatially homogeneous steady states, we notice
\begin{align}
 \partial_2 f_2(\ueps,\veps,\weps) = -(\mu_3+\delta)<0.\label{eq_final_lemma}
\end{align}
First we verify the stability or instability for the kinetic system in assertions (2), (3).
Due to \eqref{eq_final_lemma} and Lemma \ref{L4.7},
for all $\lambda^{\delta} \in \sigma(A^{\delta})$ it holds that either $ \lambda^{\delta} \in \bigcup_{\lambda \in \sigma(\MB)} B_{c(\delta)}(\lambda)$ with $\lim_{\delta \rightarrow 0} c(\delta) = 0$ or $\lambda^{\delta} \in \R \rightarrow -\infty$, where $B_{c(\delta)}(\lambda)$ denotes the open disk $\{ \mu \in  \mathbb{C} \mid |\mu - \lambda| < c(\delta)\}$.
By Lemma \ref{lem:stabKin}, the Jacobi matrix $\MB$ for $(u_-(v_-), v_-)$ has a positive eigenvalue, whereas both eigenvalues of $\MB$ for $(u_-(v_+), v_+)$ have negative real part. Therefore assertions for the kinetic system in (2) and (3) hold true.
Second, we obtain the destabilization of spatially homogeneous steady state $(u_-(v_+), v_+)$ for $D>0$ by applying Lemma 4.9, completing the proof of assertion (3). Since $(u_-(v_-), v_-)$ is an unstable equilibrium of the kinetic system, it is unstable also as a steady state of \eqref{3er-sys}. 
To obtain stability of steady states with jump discontinuity and also of $(0,0,0)$, we show that the conditions of Theorem \ref{lem:stabtrans} are satisfied. Condition $(1)$ is satisfied, see Lemma \ref{lem:entries} and Corollary \ref{lem:stabRD}. Condition $(2)$ is satisfied due to Lemma \ref{lem:3-ub}, and condition $(3)$ is satisfied due to \eqref{eq_final_lemma} and 
\begin{equation}
 a_{33}=-\nu_3-\beta \ueps < 0.
\end{equation}
Assertion $(5)$ is a consequence of Corollary \ref{lem:stabRD} and Theorem \ref{lem:stabtrans}. Assertion (4) follows from Lemma \ref{lem:entries} and Theorem \ref{lem:stabtrans}. This finishes the proof.
\qed

\section{Acknowledgments}
This work was undertaken in the framework of German-Japanese University Partnership Program (HeKKSaGOn Alliance).Anna Marciniak-Czochra and Steffen H\"arting were supported by European Research Council Starting Grant No 210680 
`Multiscale mathematical modelling of dynamics of structure formation in cell systems' and Emmy Noether Program of DFG. Steffen  H\"arting was a member of the Heidelberg Graduate School of Mathematical and Computational Methods for the Sciences. He also wants to express his gratitude for `Baden-W\"urttemberg Stipendium plus' scholarship of Baden-W\"urttemberg Stiftung .
Izumi Takagi was supported in part by JSPS Grant-in-Aid for Scientific Research (A) \#22244010 `Theory of Differential Equations Applied to Biological Pattern Formation---from Analysis to Synthesis' and \#26610027 `Control of Patterns by Multi-component  Reaction-Diffusion Systems of Degenerate Type'. 

The authors would like to express their sincere gratitude to the referees for their careful reading of the manuscript and their helpful comments.

\end{document}